\input amssym.def    \input amssym.tex



\magnification 1000

\catcode`\@=11

\hsize=150 mm   \vsize =210mm
\hoffset=4mm    \voffset=10mm
\pretolerance=500 \tolerance=1000 \brokenpenalty=5000

\catcode`\;=\active
\def;{\relax\ifhmode\ifdim\lastskip>\z@
\unskip\fi\kern.2em\fi\string;}

\catcode`\:=\active
\def:{\relax\ifhmode\ifdim\lastskip>\z@\unskip\fi
\penalty\@M\ \fi\string:}

\catcode`\!=\active
\def!{\relax\ifhmode\ifdim\lastskip>\z@
\unskip\fi\kern.2em\fi\string!}

\catcode`\?=\active
\def?{\relax\ifhmode\ifdim\lastskip>\z@
\unskip\fi\kern.2em\fi\string?}

\def\^#1{\if#1i{\accent"5E\i}\else{\accent"5E #1}\fi}
\def\"#1{\if#1i{\accent"7F\i}\else{\accent"7F #1}\fi}


\catcode`\@=12

\newif\ifpagetitre      \pagetitretrue
\newtoks\hautpagetitre  \hautpagetitre={\hfil}
\newtoks\baspagetitre   \baspagetitre={\hfil}

\newtoks\auteurcourant  \auteurcourant={\hfil}
\newtoks\titrecourant   \titrecourant={\hfil}

\newtoks\hautpagegauche \newtoks\hautpagedroite
\hautpagegauche={\hfil\the\auteurcourant\hfil}
\hautpagedroite={\hfil\the\titrecourant\hfil}

\newtoks\baspagegauche  \baspagegauche={\hfil\tenrm\folio\hfil}
\newtoks\baspagedroite  \baspagedroite={\hfil\tenrm\folio\hfil}

\headline={\ifpagetitre\the\hautpagetitre
\else\ifodd\pageno\the\hautpagedroite
\else\the\hautpagegauche\fi\fi}

\footline={\ifpagetitre\the\baspagetitre
\global\pagetitrefalse
\else\ifodd\pageno\the\baspagedroite
\else\the\baspagegauche\fi\fi}

\hautpagetitre={\hfill\tenbf Preliminary version: Not for diffusion\hfill}
\hautpagegauche={\tenbf\folio\hfill\tenrm\the\auteurcourant}
\hautpagedroite={\tenrm\the\titrecourant\hfill\tenbf\folio}
\baspagegauche={\hfil} \baspagedroite={\hfil}
\auteurcourant{Marmi and Sauzin}
\titrecourant{Quasianalytic monogenic solutions $\ldots$}
\def\mois{\ifcase\month\or January\or February\or March\or April\or
May\or June\or July\or August\or September\or October\or November\or
December\fi}
\def\Date{\rightline{\mois\ /\ \the\day\ /\/ \the\year}}
\hfuzz=1pt
\font\tit=cmb10 scaled \magstep1


\catcode`\@=11

\hoffset=4mm    \voffset=10mm
\pretolerance=500 \tolerance=1000 \brokenpenalty=5000

\def\^#1{\if#1i{\accent"5E\i}\else{\accent"5E #1}\fi}
\def\"#1{\if#1i{\accent"7F\i}\else{\accent"7F #1}\fi}

\catcode`\@=12

\newif\ifpagetitre      \pagetitretrue
\newtoks\hautpagetitre  \hautpagetitre={\hfil}
\newtoks\baspagetitre   \baspagetitre={\hfil}

\newtoks\auteurcourant  \auteurcourant={\hfil}
\newtoks\titrecourant   \titrecourant={\hfil}

\newtoks\hautpagegauche \newtoks\hautpagedroite
\hautpagegauche={\hfil\the\auteurcourant\hfil}
\hautpagedroite={\hfil\the\titrecourant\hfil}

\newtoks\baspagegauche  \baspagegauche={\hfil\tenrm\folio\hfil}
\newtoks\baspagedroite  \baspagedroite={\hfil\tenrm\folio\hfil}

\headline={\ifpagetitre\the\hautpagetitre
\else\ifodd\pageno\the\hautpagedroite
\else\the\hautpagegauche\fi\fi}

\footline={\ifpagetitre\the\baspagetitre
\global\pagetitrefalse
\else\ifodd\pageno\the\baspagedroite
\else\the\baspagegauche\fi\fi}

\def\mois{\ifcase\month\or January\or February\or March\or April\or
May\or June\or July\or August\or September\or October\or November\or
December\fi}
\def\Date{\rightline{\mois\ /\ \the\day\ /\/ \the\year}}

\hfuzz=1pt 
\overfullrule=3pt
\font\tiny=cmr8

\font\eightrm=cmr8
\font\tit=cmb10 scaled \magstep1



\def\al{\alpha}
\def\be{\beta}
\def\ga{\gamma}
\def\Ga{{\Gamma}}
\def\de{\delta}
\def\De{\Delta}
\def\la{\lambda}
\def\La{\Lambda}
\def\om{\omega}
\def\Om{\Omega}
\def\th{\theta}
\def\eps{\varepsilon}
\def\sig{\sigma}
\def\ze{{\zeta}}
\def\ka{\kappa}
\def\ph{\varphi}


\def\fS{{\Sigma}}


\def\B{\Bbb B}
\def\Bd#1{\B_{#1}}
\def\Bda#1#2{\Bd#1({#2})}

\def\R{{\Bbb R}}

\def\Z{{\Bbb Z}}
\def\Q{{\Bbb Q}}
\def\C{{\Bbb C}}
\def\N{{\Bbb N}}
\def\D{{\Bbb D}}
\def\E{{\Bbb E}}
\def\S{{\Bbb S}}
\def\L{{\Bbb L}}


\def\cA{{\cal A}}
\def\cB{{\cal B}}
\def\cC{{\cal C}}
\def\cD{{\cal D}}

\def\cL{{\cal L}}
\def\cN{{\cal N}}
\def\cR{{\cal R}}

\def\cT{{\cal T}}

\def\cE{{\cal E}}
\def\cF{{\cal F}}
\def\cG{{\cal G}}
\def\cK{{\cal K}}
\def\cM{{\cal M}}
\def\cO{{\cal O}}
\def\cS{{\cal S}}

\def\JJ{{{\cal J}_\psi}}


\def\resp{{resp.}~}
\def\cf{{cf.}~}
\def\ie{{i.e.}\ }
\def\eg{{e.g.}\ }



\def\nor{\Vert}         
\def\ao{\{\,}          
\def\af{\,\}}          
\def\dist{\mathop{\hbox{\rm dist}}\nolimits}
\def\length{\mathop{\hbox{\rm length}}\nolimits}
\def\INT{\mathop{\hbox{\rm int}}\nolimits}
\def\id{\mathop{\hbox{\rm Id}}\nolimits}
\def\demi{\frac{1}{2}}
\def\ii{^{-1}}
\def\pa{\partial}
\def\IM{\mathop{\Im m}\nolimits}
\def\RE{\mathop{\Re e}\nolimits}

\def\limproj{\mathop{\oalign{lim\cr\hidewidth$\longleftarrow$\hidewidth\cr}}}

\def\modZ{\, (\hbox{\rm mod}\,\Z)}


\def\pppar{\smallskip\par}
\def\ppar{\medskip\par}
\def\ens{\enspace}

\def\dst{\displaystyle}
\def\sst{\scriptstyle}
\def\txt{\textstyle}

\def\tst{\textstyle}
\def\ti{\tilde}
\def\ov{\overline}


\def\mbox{\hbox}
\def\text#1{\;\hbox{#1}\;}
\def\em{\sl}
\def\:{\>}
\def\frac#1#2{{#1\over #2}}


\outer\def\beginsection#1\par{\vskip0pt plus.3\vsize\penalty-50
  \vskip0pt plus-.3\vsize\bigskip\vskip\parskip
  \noindent{\bf #1 }\nobreak\bigbreak\noindent}

\outer\def\Def#1#2{\bigbreak \noindent {\bf #1\enspace }{\sl #2\par}%
\ifdim\lastskip<\bigskipamount \removelastskip\penalty55\medskip\fi}


\outer\long\def\Proc#1#2{\bigbreak \noindent {\bf #1\enspace }{\sl #2\par}%
\ifdim\lastskip<\bigskipamount \removelastskip\penalty55\medskip\fi}


\def\remark#1#2{\bigbreak \noindent {\sl Remark\ #1\enspace }{#2\par}%
\ifdim\lastskip<\bigskipamount \removelastskip\penalty55\medskip\fi}

\outer\long\def\longremark#1#2{\bigbreak \noindent {\sl Remark\ #1\enspace }{#2\par}%
\ifdim\lastskip<\bigskipamount \removelastskip\penalty55\medskip\fi}

\def\proof{\bigbreak \noindent{\sl Proof:\ }}
\def\Pf#1{\bigbreak \noindent{\sl #1:\ }}

\def\qed{\hfill$\square$\par\bigbreak} 


\def\BWD{Borel-Wolff-Denjoy}

\def\Eisen#1{\mathop{\sum\nolimits^e}\limits_{#1\hphantom{e}}}
\def\bsum#1#2{\sum_{\hbox{$ {#1}\atop{#2} $}}}
\def\BigEisen#1#2{\mathop{\sum\nolimits^e}\limits_{\hbox{$ {#1}\atop{#2} $}
\hphantom{e}}}

\def\DC{\hbox{\rm DC}}
\def\oDC{\underline{\DC}}
\def\CH{{\cal C}^{\infty}_{hol}}
\def\IQ{\hbox{\rm QI}}

\def\REC{\hbox{\rm REC}_\la}
\def\rec{\hbox{\rm rec}_\la}
\def\CM#1{\cC^#1(\la,\{M_n\},B)}
\def\CMr#1{\cC^#1(\la,\{M_n\},B_r)}
\def\CMs{\cC(\la,\{M_n\},B)}
\def\LF#1#2{\ti\cL_{(#1\rightarrow #2)}}

\def\Sg{\S_\la^{<}}
\def\Sd{\S_\la^{>}}
\def\Sgd{\S_\la^{\SupInf}}
\def\fSg{\fS_\la^{<}}
\def\fSd{\fS_\la^{>}}
\def\fSgd{\fS_\la^{\SupInf}}
\def\Gag{\Ga_\la^{<}}
\def\Gad{\Ga_\la^{>}}

\def\degnl{\de_{n,l}^{<}(\la)}
\def\dednl{\de_{n,l}^{>}(\la)}
\def\degdnl{\de_{n,l}^{\SupInf}(\la)}
\def\degdlrl{\de_{l+1+r,l}^{\SupInf}(\la)}
\def\Agdn{A_{n,\la}^{\SupInf}}
\def\deg{\de_\la^{<}}
\def\ded{\de_\la^{>}}
\def\degd{\de_\la^{\SupInf}}
\def\Egd{E^{\SupInf}}
\def\cGg#1#2#3{\cG_{#1}^{<}(#2,#3)}
\def\cGd#1#2#3{\cG_{#1}^{>}(#2,#3)}
\def\cGgd#1#2#3{\cG_{#1}^{\SupInf}(#2,#3)}

\def\SupInf{{\raise4pt\hbox{$\scriptscriptstyle >$}
\hskip-.55em\raise.5pt\hbox{$\scriptscriptstyle <$}}}
\def\Deg{\De^{\!^<}}
\def\Ded{\De^{\!^>}}
\def\Degd{\De^{\!\SupInf}}


\def\Eqno#1{
\eqno#1}

\def\Hfill#1{
\hfill#1}

\def\LLap#1{
\llap#1}

\def\EspNo#1{
&#1}


\def\equaCohom{(1.1)}
\def\equadefinfg{(1.2)}
\def\equaCohomCercle{(1.3)}
\def\equadefinF{(1.4)}
\def\equadecompf{(1.5)}

\def\equadefKn{(2.1)}
\def\ineqCondQA{(2.2)}
\def\equabwdR{(2.3)}
\def\equadefCpsi{(2.4)}
\def\excltoocl{(2.5)}
\def\propJJ{(2.6)}
\def\defJJi{(2.7)}
\def\defJJii{(2.8)}
\def\eqdecCpk{(2.9)}
\def\equadefSRa{(2.10)}
\def\equadefSrB{(2.11)}
\def\ineqDistLa{(2.12)}
\def\equacondarlin{(2.13)}
\def\equaSerLengj{(2.14)}
\def\equaDerivj{(2.15)}
\def\ineqDistRac{(2.16)}
\def\ineqconcllem{(2.17)}
\def\ineqpflemun{(2.18)}
\def\ineqpflemdeux{(2.19)}

\def\ineqlem{(3.2)}

\def\ineqtruncated{(3.3)}

\def\defichin{(3.4)}
\def\defipsi{(3.5)}
\def\defiGun{(3.6)}
\def\defiGdeux{(3.7)}
\def\defiG{(3.8)}
\def\defipsip{(3.11)}
\def\defipsim{(3.12)}
\def\defiPsip{(3.13)}
\def\defiPsim{(3.14)}

\def\lemmongev{3.2}

\def\ineqdiophS{(4.1)}
\def\PhiPsi{(4.2)}
\def\eqFE{(4.3)}
\def\eqphimr{{(4.4)}}
\def\eqGa{{(4.5)}}
\def\eqSF{{(4.6)}}
\def\defiPsiun{(4.7)}
\def\eqBT{{(4.8)}}
\def\eqTB{{(4.9)}}

\def\eqGaser{(5.1)}
\def\eqidexp{(5.2)}
\def\afbrjuno{(5.3)}
\def\Rieslermap{(5.4)}
\def\Rieslermapz{(5.5)}
\def\lF{(5.6)}
\def\eqSSM{(5.7)}

\def\thmConjGam{5.1}
\def\thmris{5.2}
\def\thmmongev{3.3}
\def\thmrisgev{5.3}

\def\lineqcircle{\equaCohomCercle}

\def\eqxnm{\hbox{$(\hbox{A}\,3.1)$}}
\def\eqAlt{\hbox{$(\hbox{A}\,3.2)$}}
\def\ineqfc{\hbox{$(\hbox{A}\,3.3)$}}
\def\defrpm{\hbox{$(\hbox{A}\,4.1)$}}


\def\noteHerPoin{$^1$}
\def\noteBranch{$^2$}
\def\noteCar{$^3$}
\def\noteRiemsurf{$^4$}
\def\noteCDD{$^5$}
\def\noteGam{$^6$}
\def\noteGN{$^7$}
\def\noteRam{$^8$}


\input epsf


\titrecourant{Quasianalytic monogenic solutions of a cohomological equation}
\auteurcourant{S. Marmi and D. Sauzin}
\centerline{\tit Quasianalytic monogenic solutions of a cohomological equation}
\par
\vskip 1. truecm 
\rightline{8 December 2000}
\vskip 1.5 truecm 
\centerline
{Stefano Marmi\footnote{$^1$}{Di\-par\-ti\-men\-to di Ma\-te\-ma\-ti\-ca 
e In\-for\-ma\-ti\-ca, Universit\`a di Udine, Via delle Scienze 206, 
Loc.\ Rizzi, 33100 Udine, Italy; e-mail: 
marmi@dimi.uniud.it}
and David Sauzin\footnote{$^2$}{
CNRS - Institut de M\'ecanique C\'eleste (UMR 8028) 77, avenue 
Denfert-Rochereau, 75014 Paris, France;  e-mail: sauzin@bdl.fr}
}
\vskip 2. truecm


\centerline{\bf Abstract} 
\vskip .5 truecm

We prove that the solutions of a cohomological equation of complex dimension one and in the
analytic category have a monogenic dependence on the parameter, and we investigate the
question of their quasianalyticity.
This cohomological equation is the standard linearized conjugacy equation for
germs of holomorphic maps in a neighborhood of a fixed point.
The parameter is the eigenvalue of the linear part, denoted by~$q$.

Borel's theory of non-analytic monogenic functions has been first investigated by Arnol'd and
Herman in the related context of the problem of linearization of analytic diffeomorphisms
of the circle close to a rotation. 
Herman raised the question whether the solutions of the cohomological equation
had a quasianalytic dependence on the parameter~$q$.
Indeed they are analytic for $q\in\C\setminus\S^1$,
the unit circle~$\S^1$ appears as a natural boundary (because of resonances,
\ie roots of unity),
but the solutions are still defined at points of~$\S^1$ which lie ``far enough from resonances''.
We adapt to our case Herman's construction of an increasing sequence of compacts which avoid 
resonances and prove that the solutions of our equation
belong to the associated space of monogenic functions;
some general properties of these monogenic functions and
particular properties of the solutions are then studied.

For instance the solutions are defined and admit asymptotic expansions at the points of~$\S^1$
which satisfy some arithmetical condition, and the classical Carleman Theorem allows us
to answer negatively to the question of quasianalyticity at these points.
But resonances (roots of unity) also lead to asymptotic expansions, for which
quasianalyticity is obtained as a particular case of \'Ecalle's theory of resurgent functions.
And at constant-type points, where no quasianalytic Carleman class contains the solutions,
one can still recover the solutions from their asymptotic expansions
and obtain a special kind of quasianalyticity.

Our results are  obtained by reducing the problem, 
by means of Hadamard's product, to the study of a fundamental solution
(which turns out to be the so-called $q$-logarithm or ``quantum logarithm'').
We deduce as a corollary of our work the proof of a conjecture of 
Gammel on the monogenic and quasianalytic properties of a 
certain number-theoretical Borel-Wolff-Denjoy series.


\vfill\eject

\centerline{CONTENTS}

\vskip .5 truecm

\item{\bf 1.}{\bf Introduction}

\medskip

\item{\bf 2.}{\bf Monogenic properties of the solutions of the cohomological equation}

\medskip

\item{2.1} ${\cal C}^{1}$-holomorphic and ${\cal C}^\infty$-holomorphic functions

\medskip

\item{2.2} Borel's monogenic functions

\medskip

\item{2.3} Domains of monogenic regularity: The sequence~$(K_j)$

\medskip

\item{2.4} Monogenic regularity of the solutions

\medskip

\item{2.5} Whitney smoothness of monogenic functions

\medskip

\item{\bf 3.}{\bf Carleman classes at Diophantine points}

\medskip

\item{3.1} Carleman and Gevrey classes

\medskip

\item{3.2} Gevrey asymptotics at Diophantine points for monogenic functions 

\medskip

\item{3.3} Borel transform at quadratic irrationals for the fundamental solution

\medskip

\item{3.4} Deduction of Theorem~3.4 from Theorem~3.5

\medskip

\item{3.5} Proof of Theorem~3.5

\medskip

\item{\bf 4.}{\bf Resummation at resonances and constant-type points}

\medskip

\item{4.1} Asymptotic expansions at resonances

\medskip

\item{4.2} Resurgence of the fundamental solution at resonances

\medskip

\item{4.3} Proofs of Theorems 4.2 and 4.3

\medskip

\item{4.4} A property of quasianalyticity at constant-type points

\medskip

\item{\bf 5.}{\bf Conclusions and applications}

\medskip

\item{5.1} Gammel's series 

\medskip

\item{5.2} An application to the problem of linearization of analytic diffeomorphisms of
the circle

\medskip

\item{5.3} An application to a nonlinear small divisor problem (semi-standard map)

\medskip

\item{}{\bf Appendix}

\medskip

\item{A.1} Hadamard's product

\medskip

\item{A.2} Some elementary properties of the fundamental solution

\medskip

\item{A.3} Some arithmetical results. Continued fractions

\medskip

\item{A.4} Proof of Lemma~3.3

\medskip

\item{A.5} Reminder about Borel-Laplace summation


\vfill\eject


\beginsection{1. Introduction}

{\bf 1.1}\ens
%
%
Let $q$ a complex number, $g(z)$ a germ of holomorphic function which vanishes at~0, 
and consider the {\sl one-dimensional cohomological equation}
$$
f(qz)-f(z) = g(z), \Eqno\equaCohom
$$
where the unknown function~$f$ is required to vanish at~0.
If $|q|\neq1$ there is a unique solution, which
can be obtained directly 
by iterating the equation forwards or backwards: 
$$
f(z) = f^{-}_{g}(q,z)=-\sum_{m\ge 0}g(q^mz)   \ens\text{if}\; |q|<1 , \quad
f(z) = f^{+}_{g}(q,z)=\sum_{m\ge 1}g(q^{-m}z) \ens\text{if}\; |q|>1.
$$
These two series are uniformly convergent in each compact 
subset of $\D\times\D_{r}$ or $\E\times\D_{r}$ respectively,
where the factor~$\D_r$ denotes the disk of convergence of~$g$
and the first factor corresponds to the parameter~$q$, with
$$
\D=\{q\in\C \mid \; |q|<1 \},\quad 
\E=\{q\in\C \mid \; |q|>1 \}.
$$
Thus we get two holomorphic functions of $q$ and~$z$.
We will be particularly interested in their dependence on~$q$,
and specifically in the relationship
between these two functions of~$q$: Is it possible to cross the unit circle which
separates  one domain of analyticity from the other?


At a formal level, we obviously obtain from the
Taylor expansion of $g(z)=\sum_{n=1}^\infty g_{n}z^n$
a unique power series satisfying~\equaCohom:
$$
f(z)=
f_{g}(q,z) = \sum_{n\ge 1}g_{n}{z^n\over q^n-1} \Eqno\equadefinfg
$$
which, as a series of functions of~$q$ and~$z$, converges 
towards $f^{-}_{g}$ in $\D\times\D_{r}$ 
and towards $f^{+}_{g}$ in $\E\times\D_{r}$. 
The case where $|q|=1$ gives rise to the simplest non-trivial small divisor problem.
Each root of the unity appears indeed as a ``resonance'', \ie a pole for some terms of
this series, and it is easy to define by an appropriate arithmetical condition a subset
of full measure of~$\S^1=\{|q|=1\}$ for which the serie converges. 
Our purpose will be to investigate the behaviour of~$f$ in the neighborhood of this set 
but also near the roots of unity,
from the point of view of regularity and asymptotic expansions.

\bigbreak\noindent
{\bf 1.2}\ens
Equation~\equaCohom\ arises naturally in the 
study of the existence of analytic conjugacies of germs of 
holomorphic diffeomorphisms of $({\C}, 0)$ with their linear part $z\mapsto qz$;
it is called cohomological because
it is the linearization of the conjugacy equation.
The study of the $q$-dependence is needed to investigate the dependence on parameters of Fatou 
components (more specifically Siegel disks)
in the dynamics of families of rational maps on the Riemann sphere [Ris].
The conformal change of variables $z=e^{2\pi i w}$, $q=e^{2\pi i h}$ transforms~\equaCohom\ into 
$$
{\cal F}(w+h)-{\cal F}(w)={\cal G}(w), 
\Eqno\equaCohomCercle
$$
where the given function~$\cG(w) = g(e^{2\pi i w})$ is 1-periodic, 
analytic in the infinite semi-cylinder $\IM w >-\de$ for some $\de\in\R$ 
and tends to zero at infinity, and the unknown function $\cF$ is required to have 
the same properties.
In this form, but under the assumption that ${\cal G}$ be
1-periodic and analytic in the complex strip 
$|\IM w |<\delta$,
the cohomological equation has been studied in detail by many authors, 
especially Wintner [Wi], Arnol'd [Ar] and Herman [He], 
since it is the linearization of the conjugacy equation of an analytic 
circle diffeomorphism to the rotation $w\mapsto w+h$. 
If $h$ is real a small divisor problem occurs once again. 

\bigbreak\noindent
{\bf 1.3}\ens
Let us return to the solutions of~\equaCohom.
We will call {\sl fundamental solution} the function
$$
f_{\delta}(q,z) = \sum_{n\ge 1}{z^n\over q^n-1}
$$ 
which is obtained in the particular case where $g(z)=\de(z)=\frac{z}{1-z}$.
In view of~\equadefinfg, we recover the general solution~$f_g$ by using the Hadamard
product with respect to~$z$:
$f_g = f_\de\odot g$.
Here, the Hadamard product of two formal series 
$A=\sum A_n z^n$ and $B=\sum B_n z^n$ is defined to be
$A\odot B = \sum A_n B_n z^n$ (see Appendix~A.1).
The formula 
$$
F(q)g= f_g(q,\cdot) = f_{\delta}(q,\cdot)\odot g(\cdot)
$$
defines a mapping~$F$ from~$\D\cup\E$ to some space of linear operators.
For all $r>0$ we denote by $H^\infty (\D_{r})$ the Banach algebra of the functions 
which are holomorphic and bounded in~$\D_r=\{|z|<r\}$ (equipped with the norm of the
supremum  over~$\D_r$), and we consider the subspace
$
B_r=zH^\infty (\D_{r})
$
of the functions which vanish at the origin.
We can now consider $F$ as a mapping
$$
F=F_{{r_1},{r_2}} \,:\ens
\D\cup\E \; \rightarrow \; \cL(B_{r_{1}}, B_{r_{2}}) \Eqno\equadefinF
$$
for $r_{1}>0$ and $r_{2}\in\left]0,r_{1}\right[$.
This allows one to describe in a compact way all the solutions of~\equaCohom\ and
to reduce most of the questions to the study of the fundamental solution.

\bigbreak\noindent
{\bf 1.4}\ens
To investigate the behaviour of the solutions for $q$ near the unit circle, we introduce a few
notations in connection with the roots of unity which appear as simple poles in~\equadefinfg.
For $m\in\N^{*}$, we set 
$\cR_{m} = \{\La\in\C\,|\ens \La^m=1\}$ 
(roots of unity of order $m$) 
and  
$\cR_{m}^{*} = \{\La= e^{2\pi i n/m},\; (n|m)=1 \}$ 
(primitive roots of order~$m$). 
We will denote by
$$\cR = \bigcup_{m\ge1} \cR_{m} = \bigsqcup_{m\ge1} \cR_{m}^{*}$$
the set of all roots of unity. 
To each $\La\in\cR$ we associate its order 
$m(\La) = \min\{m\in\N^{*}\,\mid\,\La\in\cR_{m}\}$
so that $\La\in\cR_{m(\La)}^{*}$.
Considered as an analytic function in~$(\D\cup\E)\times\D$,
the fundamental solution satisfies the following easy but important identity:
$$
f_{\delta}(q,z) = 
\sum_{\La\in\cR}{\Lambda\over q-\Lambda} \cL_{m(\La)}(z),    
\quad\text{with}\ens
\forall m\ge 1,\ens {\cal L}_{m}(z) = -{1\over m}\log (1-z^m)
\Eqno\equadecompf
$$
(see Appendix~A.2, Lemma~A2.1).
This formula,
which may be viewed as a ``decomposition into simple elements'',
is in fact an example of {\sl Borel-Wolff-Denjoy series} (see Section 2.2).
By using the Hadamard product we immediately obtain an analogous formula for the
general solution~$f_g$, or more globally for the mapping~$F_{r_1,r_2}$.

Such a formula suggests an analogy with meromorphic functions. Indeed, for each $\La\in\cR$,
we will see that 
$(q-\La)f_\de(q,z)$ tends to $\La\cL_{m(\La)}$ as $q$ tends to~$\La$ non-tangentially with respect
to the unit circle (uniformly in~$z$),
\ie $f_\de$ behaves as a function with a simple pole at~$\La$.
There is even a ``Laurent series'' at~$\La$: an asymptotic expansion which is valid
near~$\La$, inside or outside the unit circle.  But this asymptotic series must be
divergent, since there are singularities infinitely close to~$\La$: the unit circle is a
natural boundary of analyticity for~$f_\de(.,z)$, and the same is true for~$F_{r_1,r_2}$.

\bigbreak\noindent
{\bf 1.5}\ens
On the other hand, we already mentioned that $f_\de$ or~$F_{r_1,r_2}$ are defined when $q$ lies in
a special subset of~$\S^1$. There too, restricting ourselves to Diophantine points, we will find
asymptotic expansions.
We will study the Gevrey properties of those various series, and discuss the question of
{\em quasianalyticity} in the sense of Hadamard at the corresponding base-points:
we say that a space~$\cF$ of functions is quasianalytic at a point~$q_0$ if all its members admit
an asymptotic expansion at~$q_0$ and if any two functions in~$\cF$ with the same asymptotic
expansion at~$q_0$ coincide (\ie the functions of~$\cF$ are determined by their asymptotics
at~$q_0$). 
The question of quasianalyticity is a classical one for the {\em Carleman
classes}, but other spaces of functions are conceivable.

We wish also to investigate the regularity of~$f_\de$ or~$F_{r_1,r_2}$ in closed sets which
intersect the unit circle. This naturally leads to study {\em monogenic} functions in
domains which avoid the roots of unity:
in spite of the natural boundary $\{|q|=1\}$, we try to connect the function in~$\D$
and the function in~$\E$ by some monogenic continuation which would replace analytic
continuation.

Notice that, when we say that we wish to connect these two functions, our concern is
not a relationship like $f^-_g(q,z)+f^+_g(q\ii,z)=-g(z)$ 
(easy consequence of the definition of~$f_g^\pm$)
which is not ``local'' with respect to~$q$.

\bigbreak\noindent
{\bf 1.6}\ens
Section~2 deals with the definition and properties of monogenic functions; 
it gives a framework in which the solutions of the cohomological equations fall,
as shown in Section~2.4.

Section~3 is concerned with asymptotic expansions at those points of the unit circle
which satisfy Diophantine inequalities.
The question of quasianalyticity is answered negatively as far as one chooses a Diophantine
base-point associated to a quadratic irrational and considers only the classical Carleman classes.
This is in agreement with M. Herman's comment 
``The (solution of the) linearized equation
does not seem to belong to any quasianalytic class'' [He, p.\  82]. 

Section~4 proposes a constructive way to recover any solution from its asymptotic expansion 
at some particular points: roots of unity (resonances) but also constant-type points
display such a quasianalyticity property.
The {\em resurgent}  structure which appears at resonances allows one to elucidate
completely the local behaviour of the solutions
and to pass directly from the Laurent series at a given root of the unity to the whole
\BWD\ series~\equadecompf.
At constant-type points we use the Hadamard product to define a quasianalytic space which
contains the solutions.

Section~5 discusses some applications and generalizations of our work.

\bigbreak\noindent
{\bf 1.7}\ens
To conclude this introduction, let us add that
the fundamental solution~$f_\de$ is known as {\sl $q$-logarithmic series}
([Du]) but is perhaps more popular under the name of {\sl ``quantum logarithm''}. 
It is also related to  {\sl Weierstrass' $\zeta$ function}. 
The identities
$$
\eqalign
{f_{\delta}^-(q,z) &= -\sum_{n\ge1,m\ge0} z^n q^{nm} 
= \sum_{m\ge0} {zq^m\over zq^m-1}
= z{\partial\over\partial z}\log\prod_{m\ge 0} (1-zq^m) 
\; \hbox{if }\,q\in \D, \cr
f_{\delta}^+(q,z) &= \sum_{n\ge1,m\ge1} z^n q^{-nm} 
= \sum_{m\ge1} {zq^{-m}\over 
1-zq^{-m}}= -z{\partial\over\partial z}\log\prod_{m\ge1} (1-zq^{-m}) 
\; \hbox{if }\,q\in \E \cr
}
$$
show that the fundamental solution is related to the logarithmic 
derivative of Jacobi's infinite product 
([HL], [Tr]). For fixed 
$q\in\D\setminus\{0\}$, $f_{\delta}^-$ is meromorphic over~$\C$
with respect to~$z$,  with only simple poles 
at $z=q^{-m},m\ge0$. For fixed $q\in \E$, $f_{\delta}^{+}$ is meromorphic over~$\C$
with respect to~$z$,  with only simple poles 
at $z=q^{m},m\ge1$. 
On the other hand if $q$ lies on the unit circle and satisfies some arithmetical condition,
$\{|z|=1\}$ is a  natural boundary of analyticity as one can immediately check directly 
using \equaCohom\ and the fact that the r.h.s. has a pole at $z=1$
(see [Sim] for more details). 

>>From the relation with Jacobi's infinite product it immediately 
follows that Weierstrass'~$\ze$ function relative to the lattice~$\Z\oplus h\Z$
can be expressed in terms of~$f_{\de}^-$,~$f_{\de}^+$ and the corresponding Eisenstein series
$$
e_2 = \Eisen{(n,m)\in\Z^2\setminus\{(0,0)\}}{(n+mh)^{-2}},
$$
where the symbol~$\sum^e$ denotes {\sl Eisenstein summation} [We, p.\ 14].
Indeed, if $q=e^{2\pi i h}$ and $z=e^{2\pi i w}$,
$$
\ze(w) = \frac{1}{w} + e_2 w + 
         \Eisen{\om\in\Z\oplus h\Z}{\frac{1}{w+\om}}
       = e_2 w - \pi i + 2\pi i[f_{\de}^-(q,z) + f_{\de}^+(q\ii,z\ii)],
$$
where the last equality holds for $|q|<|z|<|q|\ii$
([We, p.\ 21] and [La, p.\ 248]).

\vfill \eject

\beginsection{2. Monogenic properties of the solutions of the cohomological equation}

The importance of Borel's monogenic functions in parameter-dependent 
small divisor problems was emphasized by Kolmogorov [Ko]. In his address to the 1954 
International Congress of Mathematicians (the same where he first stated 
the theorem on invariant tori in the analytic case) he considers parameter-dependent vector
fields on the two-dimensional torus and comments: ``It is possible that 
the dependence $\ldots$ on the parameter $\ldots$ is related to the class 
of functions of the type of monogenic Borel functions $\ldots$''

In his work~[Ar] on the local linearization problem of analytic diffeomorphisms 
of the circle, Arnol'd discussed in detail this issue; he 
complexified the rotation number but he did not prove that the 
dependence of the conjugacy on it is monogenic.
This point was dealt with by M. Herman [He]. Later, Risler [Ris]
extended considerably some parts of Herman's work showing that 
the parameter-dependence is Whitney-smooth also if one assumes less 
restrictive arithmetical conditions (\ie the Brjuno condition used 
by Yoccoz in [Y1, Y2, Y3]). However he did not investigate monogenic 
properties.
One should also mention that Whitney smooth dependence on parameters 
has been established also in the more general framework of KAM theory 
by P\"oschel [P\"o] who did not however consider complex frequencies. 

Borel [Bo] wanted to extend the notion of 
holomorphic function so as to allow, in certain situations, analytic continuation
through what is considered as a natural boundary of analyticity in Weierstrass' theory.
One of his goals was apparently to determine, with the help of Cauchy's formula, not too
restrictive conditions which would have ensured uniqueness of the continuation, \ie a
quasianalyticity property (see~[Th]).

Extending the presentation given in [He, III.16], 
we recall in Section~2.1 some properties of ${\cal C}^{1}$ (and ${\cal 
C}^\infty$)-holomorphic mappings on a compact subset $K$ of $\C$
with values in an arbitrary complex Banach space $B$. These are 
${\cal C}^{1}$ maps in the sense of Whitney [Wh] 
which satisfy the Cauchy-Riemann condition.
Being the uniform limits of $B$-valued rational functions with poles outside~$K$,
${\cal C}^{1}$-holomorphic maps on~$K$ share many 
properties of holomorphic functions. In particular Cauchy's Theorem
and Cauchy's Formula hold,
and they are automatically ${\cal C}^\infty$-holomorphic on a subdomain of~$K$. 

Following Borel's memoir [Bo], we define in Section~2.2 the space of {$B$-valued 
monogenic functions} associated to an increasing sequence of compact subsets of~$\C$
as the projective limit of the corresponding sequence 
of spaces of ${\cal C}^{1}$-holomorphic functions.
Borel's quasianalyticity theorem for monogenic functions is then 
recalled, in a refined form extracted from [Wk]. 

In Section 2.3 we construct an increasing sequence $K_{j}$ of compact sets 
whose union has a full-measure intersection with the unit circle.
We prove in Section~2.4 that the map $q\mapsto F_{r_1,r_2}(q)$ belongs to the associated space
of monogenic functions.
This implies that there exists an increasing sequence of smaller compact sets~$K^*_{A,j}$
on which our map is ${\cal C}^\infty$-holomorphic (Section 2.5).

Unfortunately the assumptions of Borel's 
quasianalyticity theorem are too restrictive to be applied to~$F_{r_1,r_2}$. 
This is not too surprising since Borel's result 
is much more general and includes also 
monogenic functions with singularities which are dense in an open subset of~$\C$. 
%
The problem of the quasianalyticity of $q\mapsto F_{r_1,r_2}(q)$ is 
addressed in Sections 3 and 4.


\beginsection{2.1 ${\cal C}^{1}$-holomorphic and ${\cal C}^{\infty}$-holomorphic functions}

Let $(B,\Vert\;\Vert )$ be a complex Banach space. 
%
%
The following definition is taken from [He] and makes use of 
the generalization of the notion of smoothness of a function to 
a closed set due to Whitney ([St], [Wh]). 

\Def{Definition 2.1}
{Let $C$ a closed subset of~$\C$.
A continuous function $f:\,C\rightarrow B$ is said to be ${\cal C}^{1}$-holomorphic
if there exists a continuous map $f^{(1)} :\,C\rightarrow B$ such that
$$\eqalign
{\forall z\in C,\;\forall\varepsilon >0,\ens \exists\delta >0 \;/\quad
\forall z_{1},z_{2}\in C,\ens 
&|z_1-z|<\delta,\; |z_2-z|<\delta
\cr
&\qquad\Rightarrow
\Vert f(z_{2})-f(z_{1})-f^{(1)}(z_{1})(z_{2}-z_{1})
\Vert \le \varepsilon |z_{1}-z_{2}|.
\cr}
$$
}

Notice that $f^{(1)}$ in the above definition is a complex derivative:
$\bar{\partial}f =0$,
$\pa f = f^{(1)}$
and $f$ is holomorphic in the interior of~$C$.

If $C$ is compact then ${\cal C}^{1}_{hol}(C,B)$
will denote the Banach space obtained by taking as norm  
$$
||| f |||= \max\Bigl(\sup_{z\in C}\Vert f(z)\Vert\, , \, 
\sup_{z\in C}\Vert f^{(1)}(z)\Vert\, , \, 
\sup_{z_{1},z_{2}\in C, \, z_{1}\not=z_{2}}{\Vert
f(z_{2})-f(z_{1})-f^{(1)}(z_{1})(z_{2}-z_{1})\Vert\over |z_{1}-z_{2}|}
\Bigr)
$$
(see [ALG], Remark III.4 and Proposition 
III.8: in their terminology our functions define W-Taylorian 
$1$-fields; see also [Gl], pp.\  65--66).

Let $K$ be a  compact non-empty subset of~${\C}$ 
and let ${\cal C}(K,B)$  denote the uniform algebra of continuous 
$B$-valued functions on $K$.
Let ${\cal R}(K,B)$
denote the uniform algebra of continuous functions
from $K$ to $B$ which are uniformly approximated
by rational functions with all the poles outside $K$. Let 
${\cal O}(K,B)$ denote the uniform algebra of functions of 
${\cal C}(K,B)$ which are holomorphic in the interior of $K$.
Notice that $f$ belongs to one of these uniform algebras
if and only if $\ell\circ f$ belongs to the corresponding ${\C}$-valued algebra 
for all $\ell\in B^{*}$. 

The inclusions
$$
{\cal R}(K,B)
\subset {\cal O}(K,B)
\subset {\cal C}(K,B)
$$
are in general proper;
it is not too difficult to construct examples (``swiss cheeses'')
of compacts~$K$ with empty interior such that 
${\cal R}(K,{\C}) \not= {\cal O}(K,B) = {\cal C}(K,{\C})$
(see [Ga] and the construction of monogenic functions below for more details).

\Proc{Proposition 2.1}{${\cal C}^{1}_{hol}(K,B)\subset {\cal R}(K,B)$.}

\proof
Let $f \in {\cal C}^{1}_{hol}(K,B)$. By Whitney's extension theorem
([Wh], Theorem I, see also [ALG], Theorem III.5)
$f$ admits a continuously differentiable extension $F$ to a neighborhood of~$K$. 
But according to Theorem 1.1 of [Ga],
for all $\ell\in B^{*}$, 
the function $g = \ell\circ f$ which admits a continuously differentiable 
extension to a neighborhood of~$K$
and satisfies $\bar{\partial} g\equiv 0$ on $K$
necessarily belongs to~${\cal R}(K,{\C})$. Hence $f\in {\cal R}(K,B)$. 
\qed

\remark{2.1}
{As noticed by Herman, functions in ${\cal C}^{1}_{hol}(K,B)$
share some of the properties of holomorphic functions.
Let $(U_{\ell})_{\ell\ge 1}$ be the connected components of 
${\C}\setminus K $ and assume that each $\partial U_{\ell}$
is a piecewise smooth Jordan curve. If $\sum_{\ell\ge 1}\length
(\partial U_\ell)<+\infty$, Cauchy's theorem holds:
$$
\sum_{\ell=1}^\infty \int_{\partial U_\ell} f(z)\,dz=0.
$$
Indeed, since $f\in {\cal R}(K,B)$, 
one can approximate $f$ by a sequence $(r_{k})_{k\in N}$ 
of $B$-valued 
rational functions with poles off $K$. Cauchy's theorem applies to 
these rational functions and one can pass to the limit since the 
convergence is uniform. 
Moreover, if $z\in K$ satisfies
$$
\sum_{\ell=1}^\infty \int_{\partial U_\ell} {|d\zeta|\over |\zeta-z|}<+\infty,
$$
Cauchy's formula also holds:
$$
f(z) = {1\over 2\pi i} \sum_{\ell=1}^\infty \int_{\partial U_\ell} 
{f(\zeta)\over \zeta-z}\,d\zeta.
$$
However to define higher order derivatives by means of Cauchy's formula 
one needs further assumptions on~$z$ 
(namely $\sum_{\ell=1}^\infty \int_{\partial U_\ell} {|d\zeta|\over 
|\zeta-z|^{n+1}}<+\infty $ to obtain a derivative of order $n$). 
}


The following definition is taken from [Ri]; it generalizes 
Whitney ${\cal C}^{\infty}$-smoothness to the complex case.

\Def{Definition 2.2}
{Let $C$ a closed subset of~$\C$.
A function $f:\,C\rightarrow B$ is said to be
${\cal C}^{\infty}$-holomorphic
if there exist an infinite  sequence of continuous functions 
$(f^{(n)})_{n\in {\N}}:\,C\rightarrow B$ 
such that $f^{(0)}= f$ 
and, for all $n,m\ge0$, the function~$R^{(n,m)}$ defined by
$$
f^{(n)}(z_{2}) = 
\sum_{h=0}^{m}
\frac{f^{(n+h)}(z_1)}{h!}(z_{2}-z_{1})^h + R^{(n,m)}(z_{1},z_{2}), 
\qquad
z_{1},z_{2}\in C,
$$
satisfies the following property: 
$$
\forall z\in C,\;\forall\varepsilon >0,\; \exists\delta >0 \;/\;
\forall z_{1},z_{2}\in C,\;
|z_1-z|<\delta,\; |z_2-z|<\delta
\;\Rightarrow\;
\Vert R^{(n,m)}(z_1,z_2)\Vert\le \varepsilon |z_1-z_2|^{m}.
$$
}

Clearly ${\cal C}^{\infty}$-holomorphic $B$-valued functions 
on a compact set form a Fr\'echet space. 
Once again the derivatives are taken in a complex sense, 
thus $\bar{\partial} f^{(n)}=0$ for all $n\in {\N}$. 
The functions~$f^{(n)}$ are some generalized ``weak derivatives for 
$f$''; clearly 
$f$ must be analytic in the interior of $C$ and
$$
\forall n,m\in {\N},\quad
\forall z\in \INT(C),
\quad
f^{(n+m)}(z) = \pa^m f^{(n)}(z).
$$
Whitney's extension theorem applies again: 
any $f\in {\cal C}^\infty_{hol}(C,B)$ admits an infinitely differentiable extension~$F$
to~$\C\simeq {\R}^{2}$. 
Moreover for any $n\in {\N}$,
$\pa^n F$ extends~$f^{(n)}$,
but of course $F$ is not unique and $\bar{\partial}F$ need not  
vanish outside $C$. 

\beginsection{2.2 Borel's monogenic functions}


\Def{Definition 2.3}
{Let $B$ a complex Banach space and $(K_{j})_{j\in {\N}}$ an increasing sequence of compact
subsets of~${\C}$. The associated space of {\sl $B$-valued monogenic functions} 
is defined to be the projective limit
$$
{\cal M}((K_{j}),B) = \limproj {\cal C}^{1}_{hol}(K_{j},B).
$$
}

The restrictions ${\cal C}^{1}_{hol}(K_{j+1},B)\rightarrow 
{\cal C}^1_{hol}(K_{j},B)$ are continuous linear operators between Banach spaces,
thus ${\cal M}((K_{j}),B)$ is a Fr\'echet space with seminorms 
$\Vert \,.\, \Vert_{{\cal C}^{1}_{hol}(K_{j},B)}$. 

The above definition is inspired by the work of Borel [Bo]
(see also [He], p.\  81). Borel considered the case $B={\C}$ and wanted to extend 
the notions of holomorphic function and analytic continuation.
In the usual process of analytic continuation 
(defined by means of couples $([f], D(z_{0},r))$ where $[f]$ is the 
germ at~$z_0$ of a function analytic in the open disk $D(z_{0},r)$),
the domain of holomorphy of a function is necessarily open
and one cannot distinguish between the points on a natural boundary of analyticity
(see the discussion in [Re], Chapter V, 
for a nice elementary introduction, which is also related to \BWD\ series defined below). 
Borel's idea was to allow monogenic continuation through natural boundaries of 
analyticity\footnote{\noteHerPoin}
{M. Herman pointed out to us that Poincar\'e himself investigated the possibility
of generalizing Weierstrass' process of analytic continuation so as to consider
functions whose singular points are dense on an open set or a Jordan curve [P1,
P2].} 
by selecting points at which
the function is $\cC^1$-holomorphic. If the function is moreover
$\cC^\infty$-holomorphic at such a point, the question of quasianalyticity may
be raised: Is the function determined by its Taylor series?
Such a uniqueness property could depend on the choice of the sequence~$(K_j)$ which defines the
monogenic class (and not only on the union of the~$K_j$'s),
and the Cauchy formula could help to establish it.


In the rest of Section~2.2, we illustrate the previous definition by a construction due to Borel
of a special sequence~$(K_j)$ which is adapted to the case of {\sl\BWD\ series}
[Gou, Bo, Wo, De, Si]. They are the most studied examples of monogenic functions,
and quasianalyticity can be proved in their case under suitable assumptions.

Let ${\om} = (\om_{\nu})_{\nu\ge 1}$ a bounded sequence of points in~$\C$
and $\Om = \{\om_\nu\}$.
We will exclude smaller and smaller disks around these points;
the open disk of center~$\om_\nu$ and radius~$\rho$
will be denoted by~$D(\om_\nu,\rho)$.
Let $G$ be an open bounded Jordan domain which contains~$\Om$.
We fix a sequence $(r_{\nu})_{\nu\in\N^*}\in \ell^{1}({\R}^{+})$ and define
$$
K_{j}= \ov G\setminus\bigcup_{\nu\ge1} D(\om_{\nu}, 2^{-j}r_{\nu}), 
\quad C=\bigcup_{j\ge1} K_{j}.
\Eqno\equadefKn
$$
Notice the inclusions
$$
\ov G\setminus\ov\Om \subset C \subset \ov G\setminus\Om,
$$
which are in general proper.

For each each sequence $a=(a_{\nu })_{\nu \ge 1}\in \ell^{1}(B)$, we can define a
function
$$
\fS_\om(a):\; q \mapsto \bigl(\fS_{\om}(a)\bigr)(q) = \sum_{\nu =1}^\infty {a_{\nu}\over q-\om_{\nu}}
$$
which is holomorphic in~$\C\setminus\ov\Om$.
We get a linear operator 
$\fS_\om : \ell^{1}(B)\rightarrow {\cal O}({\C} \setminus \ov\Om, B)$ 
which is generally not injective (see [Wo] for some examples). 
%
But we have also the following

\Proc{Lemma 2.1}
{The operator $\fS_{\om}$ induces an injective operator from the space
$$
\ell^{1}_{r}(B)=\{a=(a_{\nu })_{\nu \ge 1}\in \ell^{1}(B)\, \mid\; 
\forall \nu\ge 1,\;\Vert a_{\nu}\Vert^{1/4} <r_{\nu }\}.
$$
into ${\cal M}((K_j),B)$.}

\proof 
Since for all $q\in K_{j}$ and $\nu\ge1$, $|q-\om_{\nu}|\ge 2^{-j}r_{\nu}
\ge 2^{-j}\Vert a_{\nu}\Vert^{1/4}$, 
it is easy to check that $\fS_{\om}(a)\vert_{K_{j}}\in {\cal C}^{1}_{hol}
(K_{j},B)$ for all $j\ge 1$. 

To prove injectivity we make use of a residue computation. 
Let $f_{j}= \fS_{\om}(a)\vert_{K_{j}}$. Let $\gamma_{j}^{(\mu)}=
\partial D(\om_{\mu},2^{-j}r_{\mu})$ with positive orientation and let 
$\Gamma_{j}^{(\mu)}$ denote the curve obtained from 
$\gamma_{j}^{(\mu)}$ replacing those parts which are covered 
by disks $D(\om_{\nu},2^{-j}r_{\nu})$ with $\nu\not=\mu$ by the 
corresponding arcs of circles $\partial D(\om_{\nu},2^{-j}r_{\nu})$
which are contained in $K_{j}$. Clearly $\Gamma_{j}^{(\mu)}$
is a countable union of arcs of circle, all positively 
oriented, and the length of $\Gamma_{j}^{(\mu)}$ is bounded by 
$2^{-j}\sum_{\nu=1}^\infty r_{\nu}$. 
If $G_{j}^{(\mu)}$ denotes the domain of ${\C}$ enclosed by~$\Gamma_{j}^{(\mu)}$,
$$
{1\over 2\pi i}\int_{\Gamma_{j}^{(\mu)}}f_{j}(q)dq = 
\sum_{\om_{\nu}\in G_{j}^{(\mu)}} a_{\nu}. 
$$
The sequence $\nu (\mu ,j)=\inf \{\nu\in\N^*\, \mid \, \om_{\nu}\in 
G_{j}^{(\mu)}, \, \om_{\nu}\not= \om_{\mu}\}$
tends to infinity as $j\rightarrow \infty$, thus
$$
\Vert 
{1\over 2\pi i}\int_{\Gamma_{j}^{(\mu)}}f_{j}(q)dq - a_{\mu}
\Vert \le \sum_{\nu (\mu ,j)}^\infty \Vert a_{\nu}\Vert 
\;\rightarrow 0\quad \hbox{as }\; j\rightarrow \infty. 
$$
This implies injectivity. 
\qed

Of course, if none of the coefficients~$a_\nu$ vanishes,
$\fS_\om(a)$ is not analytic at any point of $C$ which is an accumulation point of
the sequence~$\om$.
Borel's example 
([Bo], p.\  144) is $B={\C}$, $\{\om_\nu\}=
\{ \frac{r+si}{n};\, 1\le r,s\le n, \, (r,n)=1, \, (s,n)=1\}$, 
$a_\nu=\exp (-\exp (n^{4}))$ and
$G=\{ q\in {\C}\, \mid \, 0<\RE q< 1\, , \, 0<\IM q< 1\}$. 

A remarkable result of Borel and Winkler is the following (see also~[Tj])

\Proc{Theorem 2.1}
{We still use the notations~\equadefKn\  and
assume furthermore that $r_{\nu}<1$ for all $\nu\in\N^*$ and 
$$
\sum_{\nu=1}^\infty \Bigl(\log {1\over r_{\nu}}\Bigr)^{-1} <+\infty.
\Eqno\ineqCondQA
$$
Let
$$
K_{j}^* = G \setminus \bigcup_{\nu=1}^\infty 
                      D\Bigl(\om_{\nu},
                             2^{-j}\bigl(\log {\tst\frac{1}{r_{\nu}}}\bigr)\ii
                       \Bigr),
\quad C^{*}=\bigcup_{j=1}^\infty K_{j}^{*}.
$$
$C^*$ is included in~$C$ and
if $f\in {\cal M}((K_j),{B})$, the restriction 
$f\vert_{K_{j}^{*}}$ is ${\cal C}^\infty$-holomorphic for all 
$j\ge 1$. Moreover, if there exist $q_{0}\in C^{*}$ and $j\in\N^*$ such that 
\item{(i)} there exists a straight line $s$ such that $q_{0}\in s\cap G\subset K_{j}^{*}$,
\item{(ii)} $f^{(n)}(q_{0})=0$ for all $n\ge 0$,

\noindent
the function~$f$ vanishes identically on~$K_{j}^{*}$.
}

In particular, according to the definition of quasianalyticity given in Section~1.5,
${\cal M}((K_j),{B})$ is quasianalytic at all points of~$C^*$ which satisfy the
condition~{\sl (i)}.  
We refer to [Wk] for a proof of Theorem~2.1 (in the case where $B={\C}$,
but this restriction is not essential).

\remark{2.2} 
{Borel (without using Whitney's extension theorem)
also proves that Cauchy's formula holds: 
let $\gamma$ a simple positively oriented closed curve bounding a 
simply connected region~$D$ of~$G$. 
Let $\gamma_{j}$ denote the curve 
obtained from $\gamma$ by replacing those parts of $\gamma$ which are 
covered by disks $D(\om_{\nu},2^{-j}r_{\nu})$ by the corresponding 
parts of the circles $\pa D(\om_{\nu},2^{-j}r_{\nu})$ which are contained in~$K_{j}\cap D$ 
(see [Wk] and [Ar], section 7, for more details). 
Let $\Ga_j$ denote the union of those parts of the circles $\pa
D(\om_{\nu},2^{-j}r_{\nu})$ which are contained in~$K_{j}\cap D$ 
and not part of~$\gamma_{j}$. Then
$$
f^{(n)}(q) = {n!\over 2\pi i} \Bigl(\int_{\gamma_{j}}{f(w)\over(w-q)^{n+1}}dw 
- \int_{\Ga_j}{f(w)\over(w-q)^{n+1}}dw\Bigr),
\qquad
q\in K_j^*\cap D,\;n\in\N.
$$
}

\remark{2.3}{The previous 
theorem was proved by Winkler under less restrictive 
assumptions than those originally required by Borel, using
Carleman's Theorem (see [Ca] and Theorem 3.1 below). Note that it holds 
{\sl without any 
further assumption on the distribution of the singular points }
$(\om_{\nu})_{\nu \ge 1}$,
while for the problem we are interested in roots of unity will play a role in the
sequel.
The quasianalyticity properties of Borel-Wolff-Denjoy series are 
studied also in [Be1], [Be2] and [Si]
(which focus in fact on the broader question of the injectivity of~$\fS_{\om}$).
}

\remark{2.4}{Unfortunately one cannot apply the previous theorem 
to the solutions of cohomological equations since the condition~\ineqCondQA\ is 
too restrictive for that situation. Let $0<\rho_{2}<\rho_{1}$, 
$B={\cal L}(B_{\rho_{1}}, B_{\rho_{2}})$ and consider the mapping~\equadefinF. Ordering
the primitive roots  of unity by increasing order (\ie following the Farey ordering 
of rational numbers), one can write it as a \BWD\ series
$$
F(q)=\fS_\cR(a)(q)=
\sum_{\nu=1}^\infty {\Lambda_{\nu}\over q-\Lambda_{\nu}}
{\cal L}_{m(\nu )}\odot, 
\qquad
\cR = \{\La_1,\La_2,\ldots\},
\Eqno\equabwdR
$$
setting $a_{\nu }= \Lambda_{\nu}{\cal L}_{m(\nu )}\odot$.
Since the number of terms in the Farey series of order $m$ is 
approximately~${3m^{2}\over \pi^{2}}$ ([HW], Theorem 331, p.\  268) one 
has $m(\nu ) \simeq {\pi\over\sqrt{3}}\sqrt{\nu}$. 
On the other hand, one checks easily that
$\Vert a_\nu \Vert \simeq \frac{1}{m(\nu)} {(\frac{\rho_2}{\rho_1})}^{m(\nu)}$.
The requirement $(a_{\nu})_{\nu \ge 1}\in \ell^{1}_{r}(B)$ 
leads to a lower bound 
$r_{\nu}\ge c_1\, c_2^{m(\nu)/4}$ and the condition~\ineqCondQA\ is violated. 
}

\beginsection{2.3 Domains of monogenic regularity: The sequence~$(K_j)$}

The goal of this section is to specify a 
sequence of compact sets $(K_j)_{j\in {\N}}$
so as to be able to prove (in Section~2.4) that
${\cal M}((K_j),B)$ contains the solutions of the cohomological 
equation. 
In the definitions of the domains $C_{\psi, \ka, d}$ and 
$W^{A}_{\gamma ,\ka, d}$ given below (Definitions~2.4 and~2.6),
we will follow a construction given by M. Herman [He]
for Diophantine numbers (see also [Ris] for a similar construction for Brjuno numbers).
We adapt it slightly so as to deal with more general irrational numbers. 

\bigbreak\noindent {\bf a)}
The conformal change of variable $q=e^{2\pi ix}$ maps ${\C}^{*}$  
biholomorphically on 
${\C}/{\Z}$, 
the circle $\{|q|=1\}$ on ${\R}/{\Z}$
and $\cR^*_m$ on $\{ \frac{n}{m} \mid m\in {\N}^{*},\: 0\le n\le m-1,\: (n,m)=1\}$.
We will use the notations of Appendix~A.3 for continued fractions:
if $x\in {\R}\setminus {\Q} \modZ$,
we will denote by
$[0,a_1(x),a_2(x),\ldots]$ its continued fraction expansion
and by
${({n_{k}(x)\over m_{k}(x)})}_{k\ge 0}$ the corresponding sequence of convergents, 
omitting sometimes the dependence on~$x$.
Note that $n_0/m_0=0/1$. 

\Def{Definition 2.4}
{We call an approximation function any decreasing function~$\psi$ on~$\N^*$ such that 
$$ 
2 \sum_{m=1}^{+\infty}\psi(m)<1
\quad\text{and}\quad
\forall m\ge1,\ens 0 < \psi(m)\le {1\over 2m}.
$$
We associate with it a subset of~$\R\setminus\Q \modZ$:
$$
C_{\psi} = \Bigl\{ x\in {\R}\setminus {\Q} \modZ \mid \; \forall k\ge 0,\;
m_{k+1}(x)\le {1\over \psi (m_{k}(x))} \Bigr\},
\Eqno\equadefCpsi
$$
and some subsets of~$\C/\Z$ whose traces on~$\R/\Z$ are~$C_\psi$:
$$
C_{\psi,\ka} = \bigcup_{y \in  C_\psi}
\bigl\{ x\in\C/\Z\, \mid \, |\IM \ti x| \ge \ka |\RE (\ti x-\ti y )| \bigr\},
\quad
C_{\psi,\ka,d} = C_{\psi,\ka} \cap \bigl\{|\IM x| \le d\bigr\},
$$
for $\ka\in\,]0,1[$ and $d>0$,
where $\ti x$ and~$\ti y $ denote some lifts in~$\C$ of~$x$ and~$y$.
}

Notice that $C_\psi$ consists of points which are ``far enough from
the rationals'', as measured by~$\psi$; namely, according to~\ineqfc\ and 
Proposition~A3.2, 
$$
\bigcap_{n/m} \{ x\,|\; |x-\frac{n}{m}|\ge \frac{\psi (m)}{m}\} 
\;\subset\;
C_{\psi}
\;\subset\;
\bigcap_{n/m} \{ x\,|\; |x-\frac{n}{m}|> \frac{\psi (m)}{2m}\}.
\Eqno\excltoocl
$$
The most interesting case for our purposes will be  $\psi(m) = \ga\,e^{-\al m}$
with fixed $\al >0$ and
$\ga \in\,]0,\inf ({\al e\over 2},{e^{\al}-1\over 2})[$.
The classical Diophantine condition of exponent~$\tau>2$ 
(see Section~3.2)
would correspond to $\bigcup C_{\phi_{\ga,\tau}}$, where 
$\phi_{\ga,\tau}(m)=\ga\,m^{1-\tau}$ and the union is taken over those $\ga>0$
such that~$\phi_{\ga,\tau}$ is an approximation function
(\ie $\ga\le\frac{1}{2\ze(\tau-1)}$, denoting by~$\ze$ the Riemann zeta function).

The Diophantine exponent $\tau=2$ (which is associated to constant-type
points) was not considered here, only because the
corresponding functions~$\phi_{\ga,2}$ do not satisfy the condition of
summability in Definition~2.4. This condition is used in the next lemma to
ensure positive measure for~$C_\psi$, and indeed the set of constant-type points
has measure zero.

\Proc{Lemma 2.2}
{If $\psi$ is an approximation function,
$C_{\psi}$ has positive Lebesgue measure.
For all $\varepsilon >0$ there exists an approximation function 
$\psi$  such that $| C_{\psi} \modZ |>1 - \varepsilon$.
}

\proof  
According to~\excltoocl,
the one-dimensional Lebesgue measure of~$(\R/\Z)\setminus C_{\psi}$ is less than
$$
2\sum_{m=1}^\infty\sum_{n=0}^{m-1}{\psi (m)\over m} \le 
2\sum_{m=1}^\infty \psi (m) <1.
$$
Given $\eps>0$, we choose $\psi (m) = {\varepsilon\over 2\zeta (2)m^{2}}$ 
to make the previous quantity less than~$\eps$.
\qed

In order to investigate the structure of this kind of set, it is useful 
to refer to a suitable partition of~$\R/\Z$
obtained by considering a finite number of iterations of the Gauss map $A$
(see Appendix~A.3 for the definition of the Gauss map; 
the intervals defined below are called ``intervals of rank~$k$'' in~[Khi]).

Let $a_1,\ldots,a_k$ positive integers ($k\in\N^*$).
We associate with them the finite continued fractions 
$[0,a_1,\ldots,a_{k-1}] = \frac{n_{k-1}}{m_{k-1}}$
and
$[0,a_1,\ldots,a_{k}] = \frac{n_{k}}{m_{k}}$,
and define an interval
$$
\tst
I(a_{1},\ldots ,a_{k}) = 
\{ x=\frac{n_{k}+n_{k-1}y}{m_{k}+m_{k-1}y}, \; y\in\,]0,1[\,\}
=\cases{
] \frac{n_{k}}{m_{k}},{n_{k}+n_{k-1}\over m_{k}+m_{k-1}}[
&if $k$ is even\cr
] \frac{n_{k}+n_{k-1}}{m_{k}+m_{k-1}}, {n_{k}\over m_{k}}[
&if $k$ is odd\cr}
$$
%
(the alternative stems from~\eqAlt).
Each such interval is a branch 
of the $k$-th iterate $A^{k}$ of the Gauss map,
precisely the branch which is determined 
by the fact that all points 
$x\in I(a_{1},\ldots ,a_{k}) $ have $\{0,a_{1},\ldots ,a_{k}\}$ as 
first $k+1$ partial quotients
(see Formula~\eqxnm).
For a given $k\ge1$, the union of all branches of~$A^k$ yields a partition of~$\R/\Z$
$$
\forall k\ge1,\quad
\R/\Z = \cF_k \cup \bigcup_{a_{1},\ldots ,a_{k}\ge 1}
{I(a_{1},\ldots ,a_{k})},
$$
where\footnote{\noteBranch}{Notice that
any rational number $n/m\in\,]0,1[$ is the endpoint 
of exactly two branches of the iterated Gauss map. Indeed $n/m$ can 
be written in a unique way as $n/m=[0,\bar{a}_{1},\ldots ,\bar{a}_{\ell}]$ 
for some $\ell\ge 1$, with $\bar a_1,\ldots,\bar a_{\ell-1}\ge1$
and $\bar{a}_{\ell}\ge 2$; it is the left (right) 
endpoint of $I(\bar{a}_{1},\ldots ,\bar{a}_{\ell})$ and 
the right (left) endpoint of $I(\bar{a}_{1},\ldots,\bar{a}_{\ell-1},\bar{a}_{\ell}-1,1)$ 
if $\ell$ is even (odd). 
} 
$\cF_k = \bigl\{ [0,a_1,\ldots,a_\ell],\; 1\le \ell\le k,\; a_i\ge1 \bigr\} 
\subset\Q/\Z$.
%
The previous definition allows for a convenient rephrasing of~\equadefCpsi:
$$
C_{\psi}= \bigcap_{k\ge1}\tst\bigsqcup_{\psi}I(a_{1},\ldots ,a_{k}),
$$
where for each $k\ge1$,
$\tst\bigsqcup_{\psi}$ denotes the disjoint union over those $(a_{1},\ldots,a_{k})$
such that $m_{i+1}\le 1/\psi (m_{i})$ for $i=0,\ldots ,k-1$ 
(here, of course, $m_{i}$ is the denominator of $[0,a_{1},\ldots,a_{i}]$).

\bigbreak\noindent {\bf b)}
We will indicate some more properties of the set~$C_\psi$ associated to an approximation
function~$\psi$.
As a preliminary, to each rational number $n/m\in\Q/\Z$ we attach 
an open interval $\JJ(n/m)$ such that
$$
n/m \in \JJ(n/m) \subset (\R/\Z)\setminus C_\psi.
\Eqno\propJJ
$$
To define it we proceed as follows:
\item{(i)} if $n/m=0/1$, we set
$$
\JJ(0/1) := 
\INT\bigl(\,
\bigcup_{a_{1}> 1/\psi (1)} \ov{I(a_{1})}
\;\cup
\dst\bigcup_{a_{2}+1> 1/\psi (1)} \ov{I(1,a_{2})}
\,\bigr);
\Eqno\defJJi
$$
\item{(ii)} if $n/m\not= 0/1$ and $(n,m)=1$, we write 
$n/m = [0,\bar{a}_{1},\ldots ,\bar{a}_{k}]$, with $k\ge 1$, $\bar a_1,\ldots,\bar a_{k-1}\ge1$ and 
$\bar{a}_{k}\ge 2$, we set
$n_{-}/m_{-}=[0,\bar{a}_{1},\ldots ,\bar{a}_{k-1}]$
(if $k\ge 2$; otherwise $n_{-}/m_{-}=0/1$) and 
$$
\displaylines
{\JJ(n/m) := 
\Hfill\defJJii\cr
\hfill
\INT\bigl(\,
\bigcup_{a_{k+1}m+m_{-}>1/\psi (m)} \ov{I(\bar a_{1},\ldots ,\bar{a}_{k},a_{k+1})}\;\cup
\dst\bigcup_{(a_{k+2}+1)m-m_{-}> 1/\psi (m)} \ov{I(\bar{a}_{1},\ldots ,\bar{a}_{k}-1,1,a_{k+2})}
\,\bigr).
\cr}
$$

\noindent
This definition is motivated by the relations~\hbox{$(\hbox{A}\,3.7)$}.
For instance the points in the first union of~\defJJii\
have continued fraction expansions such that
$a_{k+1}m+m_{-}=m_{k+1}$ since $m_{k-1}=m_-$ and $m_k=m$ for them,
and in the second one,
$(a_{k+2}+1)m-m_{-}=m_{k+2}$ since $m_{k}=m-m_-$ and $m_{k+1}=m$
(except at one of the boundary-points of each interval:
$m_{k+1} = (a_{k+1}+1) m + m_-$ and $m_{k+2} = (a_{k+2}+2) m - m_-$
respectively, for these exceptional rational points).
We thus have $m_{k+1}>1/\psi(m_k)$ or $m_{k+2}>1/\psi(m_{k+1})$ respectively,
hence $\JJ(n/m)$ is contained in the complement of~$C_\psi$.
To check that it is an open interval, consider for instance the case of
odd~$k$:
using~\eqxnm\ one can write the first union as 
$\{\frac{n+n_-y}{m+m_-y}\,;\; 0< y\le 1/M\}$, where $M$ is the minimum value of~$a_{k+1}$
(\ie $M = [\frac{1}{m}(\frac{1}{\psi(m)}-{m_-})]+1$), 
this union is thus a non-empty interval whose right endpoint is~$n/m$; 
similarly the second union is a non-empty closed interval whose left endpoint is~$n/m$.
Notice that, in case~(i), $\JJ(0/1)=\INT([0,\frac{1}{M}]\cup[1-\frac{1}{M},1])$
must be identified with~$]-\frac{1}{M},\frac{1}{M}[$ (where $M=[\frac{1}{\psi(1)}]+1$).

\Proc{Lemma 2.3}
{The set~$C_{\psi}$ associated to any approximation function is totally disconnected,
closed and perfect.
}

\proof
Since $C_{\psi}\cap(\Q/\Z)=\emptyset$, $C_{\psi}$ is totally disconnected.
To see that $C_{\psi}$ is closed observe that,
if $x\in (\R/\Z)\setminus C_\psi$,
one can exhibit an open neighborhood of~$x$ which is contained in the
complement of~$C_\psi$:
either 
$x\in\Q$ and $\JJ(x)$ is such a neighborhood,
or
$x\notin\Q$ and $m_{k+1}(x)>1/\psi (m_{k}(x))$ for some $k\ge 0$,
hence $I(a_{1}(x),\ldots,a_{k+1}(x))$ will do.

We now prove that any $x \in C_{\psi}$ is an accumulation point of~$C_\psi$.
For each $j\in\N^*$ we define a linear fractional map
$$
T_{x,j}\,:
\;
y\in\,]0,1[
\;\mapsto\;
T_{x,j}(y) = {n_{j}(x)+n_{j-1}(x)y\over m_{j}(x)+m_{j-1}(x)y}
=
[0,a_1(x),\ldots,a_j(x), a_1(y), a_2(y), \ldots].
$$
Let us use $y^*={\sqrt{5}-1\over 2}=[0,1,1,\ldots]$.
The sequence
$x^{(j)} = T_{x,j}(y^*)$
converges to~$x$ as $j\to\infty$, and one can check that each $x^{(j)}\in C_\psi$:
\item{--}
If $k\le j$,
${n_{k}(x^{(j)})\over m_{k}(x^{(j)})} = {n_{k}(x)\over m_{k}(x)}$;
thus $m_{k+1}(x^{(j)}) \le {1\over\psi (m_{k}(x^{(j)}))}$ for all $k\le j-1$ and
$$
m_{j+1}(x^{(j)}) = m_{j}(x)+m_{j-1}(x)\le a_{j+1}(x)m_{j}(x)+m_{j-1}(x)=m_{j+1}(x)\le 
{1\over\psi (m_{j}(x^{(j)}))}. 
$$
\item{--}
If $k\ge j+1$, we use $\psi (m)\le 1/2m$:
$$  
m_{k+1}(x^{(j)}) = m_{k}(x^{(j)})+m_{k-1}(x^{(j)})\le 2m_{k}(x^{(j)})\le 
{1\over\psi (m_{k}(x^{(j)}))}.
$$
\qed

The intervals~$\JJ(n/m)$ defined above will also help us in the proof of the next proposition
which describes the connected components of~$(\R/\Z)\setminus C_\psi$.

\Proc{Proposition 2.2}
{Let
$$
{\Q}_{\psi} = 
\bigl\{0/1\bigr\}
\; \cup \;
\bigl\{n/m\in \Q/\Z\, \mid \;n/m\not=0/1 \;\text{and}
\; m_{j+1}\le 1/\psi (m_{j}) \;\text{for}\; j=0,\ldots,k-1\bigr\},
$$
with the usual notations and conventions:
the~$m_j$'s ($0\le j\le k$) are the denominators of 
the convergents $[0,a_1,\ldots,a_j]$ of~$n/m=[0,a_1,\ldots,a_k]$, with $a_k\ge2$.
%
\item{(1)} Each connected component of~$(\R/\Z)\setminus C_\psi$
contains one and only one point of~$\Q_\psi$,
which is a convergent of both of its endpoints.
We denote the connected component of~$n/m\in\Q_\psi$ 
by $]x_{n/m},x_{n/m}'[\subset\R/\Z$ 
(which must be identified to an open subinterval of~$]0,1[$ if $n/m\not=0/1$,
or of~$]-1/2,1/2[$ if $n/m=0/1$).
\item{(2)} ${\psi(m)\over 2m} \le |x-\frac{n}{m}| < \frac{2\psi(m)}{m}$
if $x=x_{n/m}$ or~$x_{n/m}'$.
\item{(3)} If $r/s\in \,]x_{n/m},x_{n/m}'[$ and $r/s\not= n/m$, 
$s>\frac{1}{\psi(m)}\ge 2m$. 
}

\proof 
Any connected component of $U=(\R/\Z)\setminus C_\psi$ contains at least
a rational~$r/s$. Suppose this rational does not belong to~$\Q_\psi$ and write it as
$r/s = [0,a_1,\ldots,a_\ell]$ with $a_\ell\ge2$:
we must have $m_{k+1}>1/\psi(m_k)$ for some $k\ge0$.
Choosing $k$ to be minimal, we obtain
$n/m=[0,a_1,\ldots,a_k]\in\Q_\psi$ and $r/s\in\JJ(n/m)$
(note that $n/m=0/1$ if $k=0$).
Thus the connected component of~$r/s$ contains~$\JJ(n/m)$, and in particular~$n/m$.
We notice in passing that $s>1/\psi(m)\ge 2m$.

Let us now suppose that $]x,x'[$ is the connected component of $n/m\in\Q_\psi$ in~$U$
and check that $n/m$ is a convergent of~$x$ and~$x'$.
We may suppose that $n/m\not=0/1$.
Let us write $n/m = [0,a_1,\ldots,a_k]$ with $a_k\ge2$.
Denoting by~$m_-$ the denominator of~$[0,a_1,\ldots,a_{k-1}]$, 
we choose positive integers~$a$ and~$b$ such that
$$
a m + m_- \le 1/\psi(m), \quad
(b+1) m - m_- \le 1/\psi(m)
$$
(this is possible since $1/\psi(m)\ge 2m > m+m_-$).
By the same kind of argument as at the end of the proof of Lemma~2.3, 
one can check that the points
$$
x^+ = [0,a_1,\ldots,a_{k-1},a_k,a,1^\infty]
\ens\text{and}\ens
x^- = [0,a_1,\ldots,a_{k-1},a_k-1,1,b,1^\infty]
$$
both belong to~$C_\psi$.
But if $k$ is even, $x^- < n/m < x^+$, and the order is reversed otherwise.
Therefore $[x,x']$ is contained in~$]x^-,x^+[$ (or~$]x^+,x^-[$ is $k$ is odd),
and $n/m$ is a convergent of all those points.
This implies easily that a connected component of~$U$ cannot contain more that one point
of~$\Q_\psi$. 
\smallbreak

The first inequality in~(2) follows from the second inclusion in~\excltoocl.
For the second inequality, 
consider $x^+$ and~$x^-$ as defined above
for $n/m=[0,a_1,\ldots,a_k]\in\Q_\psi\setminus\{0/1\}$, 
but this time we choose~$a$ and~$b$ maximal:
$$
\frac{1}{\psi(m)} - m < a m + m_- \le \frac{1}{\psi(m)}, \quad
\frac{1}{\psi(m)} - m < (b+1) m - m_- \le \frac{1}{\psi(m)}.
$$
By virtue of~\ineqfc,
since $m_k(x^+)=m$, $m_{k+1}(x^+)=am+m_-$, 
$m_{k+1}(x^-) = m$ and $m_{k+2}(x^-)=(b+1)m-m_-$,
$$
|x^+ - n/m|,|x^- - n/m| < \frac{\psi(m)}{m(1-m\psi(m))} \le \frac{2\psi(m)}{m}.
$$
This yields the desired inequality.
If $n/m=0/1$, one can use
$x^+ = [0,a,1^\infty]=\frac{1}{a+g}$ with $a=[\frac{1}{\psi(1)}]\ge2$ and $g=[0,1^\infty]$,
and $x^-=-1+[0,1,a-1,1^\infty]=-x^+$.
\smallbreak

(3) was already observed at the beginning of the proof.
\qed

\bigbreak\noindent {\bf c)}
We now fix $\ka\in\,]0,1[$ and~$d>0$, and study the sets~$C_{\psi,\ka}$ and~$C_{\psi,\ka,d}$
associated to the approximation function~$\psi$.
Proposition 2.2 yields a decomposition 
of~$(\R/\Z)\setminus C_\psi$ into connected components;
this will reflect in a description of the complement of~$C_{\psi,\ka}$:

\Proc{Lemma 2.4}
{For each $n/m\in \Q_\psi$, let
$$
\De_{n/m}= \{ 
x \in \C/\Z\,\mid\, 
\RE x \in\,]x_{n/m},x'_{n/m}[, \, 
|\IM x|< \ka\,\min (\RE x -  x_{n/m}, x'_{n/m} - \RE x)
\},
$$
which is an open diamond of base~$]x_{n/m},x'_{n/m}[$ and slopes~$\pm\ka$ 
with respect to the real axis.
We have
$$
C_{\psi,\ka} = (\C/\Z)\setminus 
\bigsqcup_{n/m\in {\Q}_{\psi}}\De_{n/m},
\Eqno\eqdecCpk
$$
the sets $C_{\psi,\ka,d}$ are compact subsets of~$\C/\Z$ and 
they have positive measure when $d>\ka\,\ze(4)$:
$$
\hbox{meas}\,\bigl(C_{\psi,\ka,d}\bigr) > 
2d - 8\ka\sum_{m\ge1}{\tst\bigl(\frac{\psi(m)}{m}\bigr)^2}.
$$
}

\proof
Let us rephrase the definition of~$C_{\psi,\ka}$ as
$$
C_{\psi,\ka} = (\C/\Z) \setminus \bigcap_{y \in C_\psi} \De^*_{y },
$$
with
$
\De^*_{y } = \bigl\{ x\in\C/\Z \,\mid\, |\IM\ti x| < 
               \ka\,\min(\RE\ti x - \RE\ti y , 1 + \RE y  - \RE\ti x)
             \bigr\}
$
(each~$\De^*_{y }$ is an open diamond whose trace on~$\R/\Z$ has length~1 
and coincides with the complement of~$\{y\}$).
Formula~\eqdecCpk\ is now reduced to the identity
$$
\bigsqcup_{n/m\in\Q_\psi} \De_{n/m} = \bigcap_{y \in C_\psi} \De^*_{y }.
$$
\vskip .2cm
\epsfysize=3.8cm
\centerline{\epsfbox{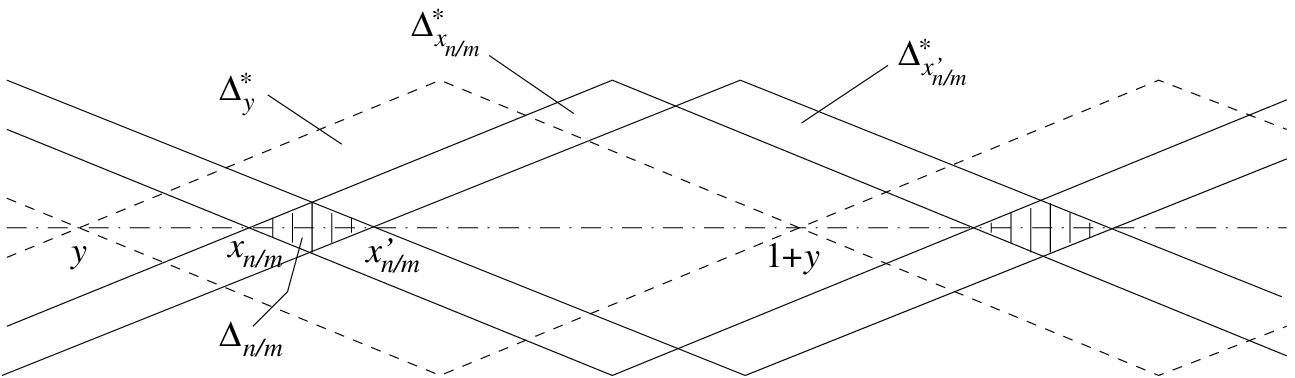}}
\vskip 1cm

If $n/m\in\Q_\psi$ and $y \in C_\psi$, the fact that $y \not\in\,]x_{n/m},x'_{n/m}[$
implies that $\De_{n/m}\subset\De^*_{y }$, hence
the union of the diamonds~$\De_{n/m}$ is contained in the intersection of the
diamonds~$\De^*_{y }$.
Let~$x$ in the intersection of the diamonds~$\De^*_{y }$.
If $x\in\R/\Z$, this means that $x\not\in C_\psi$,
thus $x\in\,]x_{n/m},x'_{n/m}[\,\subset \De_{n/m}$ for some $n/m\in\Q_\psi$.
If $x\not\in\R/\Z$,
the intersection with~$\R/\Z$ of the lines of slopes~$\pm\ka$ which pass through~$x$ 
define two points $x^-<x^+$. 
Necessarily $[x^-,x^+]\,\subset\,]x_{n/m},x'_{n/m}[$ for some $n/m\in\Q_\psi$
(because the existence of $y \in[x^-,x^+]\,\cap\, C_\psi$ would lead
to the contradiction $x\not\in\De^*_{y }$),
hence $x\in\De_{n/m}$.
Thus $x$ belongs to the union of diamonds~$\De_{n/m}$ in both cases
and this yields the reverse inclusion.

As a consequence $C_{\psi,\ka}$ is closed and 
its intersection with a strip $\bigl\{|\IM x|\le d\bigr\}$ is compact. The inequalities
$$
\forall n/m\in\Q_\psi,\quad
\hbox{meas}\,\bigl(\De_{n/m}\bigr)\, = \demi\ka(x'_{n/m}-x_{n/m})^2 
< 8 \ka {\tst\bigl(\frac{\psi(m)}{m}\bigr)^2} \le \frac{2\ka}{m^4}
$$
(which follow from Proposition~2.2~(2) and from $\psi(m)\le 1/2m$)
yield the last statement.
\qed

\remark{2.5}
{Using a suggestion by Herman ([He], Remark at p.\  81), 
one can prove that ${\cal O}(C_{\psi,\ka,d},B) = 
{\cal R}(C_{\psi,\ka,d},B)$, a result to be compared with
the general inclusion which was indicated in Section 2.1.
Notice that, since
$\dst\R \subset \bigcap_{y \in C_\psi} \ov{\De^*_{y}} =
C_\psi \cup \bigsqcup_{n/m\in {\Q}_{\psi}}\ov{\De_{n/m}}$, 
$$
\INT\bigl( C_{\psi,\ka,d}\bigr) =
\{ x\in\C/\Z\, \mid \; |\IM x|<d \}
\setminus \Bigl( 
          C_\psi \cup \bigsqcup_{n/m\in {\Q}_{\psi}}\ov{\De_{n/m}}
          \Bigr)
\subset \C\setminus \R \modZ.
$$
The idea is to apply Milnikov's theorem [Za, p.\  112] which states that,
if the inner boundary of a compact set~$K$ is a subset of an analytic curve,
${\cal O}(K,B) = {\cal R}(K,B)$. 
(The {inner boundary} of~$K$ is defined as $\partial_{I}K = \partial K\setminus
\bigsqcup \partial \De_{\ell}$, where $\bigsqcup\De_{\ell}$ is the decomposition
of~$\C\setminus K$ into disjoint connected components.
Here $\pa_I C_{\psi,\ka,d} = C_\psi \subset \R/\Z$.)
}

\remark{2.6}
{One can check that $C_{\psi,\ka,d}$ has a finite number of connected components 
and is locally connected; it is connected as soon as $d>\ka\psi(2)$. 
Also $\INT(C_{\psi,\ka,d})$ has a finite number of connected
components.
}

\vskip 1cm
\epsfysize=4cm
\centerline{\epsfbox{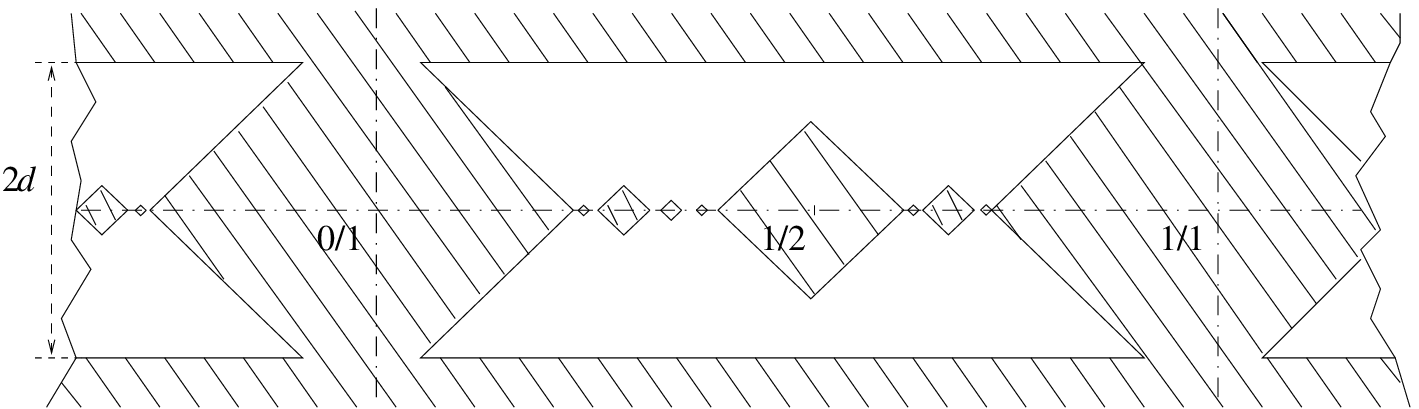}}
\vskip 1cm


\bigbreak\noindent {\bf d)}
Finally we define the sequence of compact subsets~$K_j$ of~$\C$ which will be
used in the sequel.

\Def{Definition 2.5}
{Let us fix $\ka\in\,]0,1[$, $d,\al>0$ and a decreasing sequence~$(\ga_j)_{j\ge0}$ wich tends to~0.
We assume $\ga_j < \inf ({\al e\over 2},{e^{\al}-1\over 2},1)$ for all~$j\ge0$.
We define
$$
\psi_{j}(m)= \ga_{j}e^{-\al m} \ens\text{for $m\ge1$},
\quad
K_j = \{ q=e^{2\pi i x},\; x\in C_{\psi_j,\ka,d} \},
\quad
C = \bigcup_{j\in\N} K_j.
$$
}

Observe that each $K_j$ is contained in the annulus 
$\{ e^{-2\pi d} \le |q| \le e^{2\pi d} \}$ and that
its measure tends to the measure of this annulus,
while the measure of~$K_j\cap\S^1$ tends to the measure of the circle,
as $j\to\infty$.

\remark{2.7}
{Since $C\cap {\S}^{1}=\cup_{j\in {\N}}\{\,e^{2\pi i x} \mid
x\in C_{\psi_{j}}\,\}$, by Lemma 2.3 it is a countable union of nowhere 
dense closed sets. Proposition 2.2 then shows that its complement in 
${\S}^{1}$ is a dense $G_{\de}$-set with zero 
$s$-dimensional Hausdorff measure for all $s>0$.
}

\Proc{Lemma 2.5}
{There exists a positive number~$\mu$, which depends only on~$\ka$, such that
$$
\forall j\in\N, \; \forall q\in K_j, \; \forall \La\in\cR, \quad
|q-\La| > \mu\, \frac{\psi_j(m(\La))}{m(\La)}.
$$}

\proof
Let $j\in\N$ and~$q\in K_j$.
Since $\frac{\psi_j(m)}{m}\le\frac{1}{2}$ for all~$m\ge1$, 
we may suppose that $\dist(q,\S^1)$ be less than some arbitrary constant;
thus we assume
$$
q = e^{2\pi ix}, \quad x\in C_{\psi_j,\ka}\modZ, \quad |\IM x|\le 1.
$$
We also choose $y\in C_{\psi_j}\modZ$ such that $|\IM x| \ge \ka |\RE(x-y)|$.

Let $\La\in\cR$. We choose $n/m\in\Q$ such that $\La = e^{2\pi i n/m}$ 
and $|\RE(x-\frac{n}{m})|\le\frac{1}{2}$.
According to~\excltoocl, $|y-\frac{n}{m}| > \frac{\psi_j(m)}{2m}$,
and one can check easily that
$|x-\frac{n}{m}|\ge\mu_0|y-\frac{n}{m}|$ with $\mu_0=(1+\ka^{-2})^{-1/2}$.
Thus $z=x-\frac{n}{m}$ satisfies
$$
|\RE z|\le\frac{1}{2}, \quad |\IM z|\le1, \quad |z|>\mu_0\frac{\psi_j(m)}{2m}.
$$
Hence $|q-\La| = |e^{2\pi iz}-1|$ can be bounded from below as required.
\qed

\beginsection{2.4 Monogenic regularity of the solutions}

Let~$B$ a Banach space.
We now consider the space of $B$-valued monogenic functions which corresponds to 
the sequence~$(K_j)$ of Definition~2.5.
We will see that the general solution of the cohomological equation as encoded
by the mapping~$F_{r_1,r_2}$ of Section~1.3 belongs to this space
--- recall its definition~\equadefinF\ and the notation $B_r = z H^\infty(\D_r)$;
of course $B=\cL(B_{r_1},B_{r_2})$ in that case.

More generally, we will show that the \BWD\ series of the form
$$
\fS_\cR(a):\; q \mapsto \bigl(\fS_{\cR}(a)\bigr)(q) = \sum_{\La\in\cR}
{a_{\La}\over q-\La} 
\Eqno\equadefSRa
$$
(not necessarily with the same coefficients as those of~\equabwdR\ in Remark~2.4)
are monogenic; we simply restrict ourselves to 
$$
\cS(r,B) = \bigl\{\, a=\{a_\La\}_{\La\in\cR} \;\text{sequence of~$B$ such
that}\; \exists c>0\,/\; \forall\La\in\cR, \nor a_\La \nor \le \tst
\frac{c\,r^{m(\La)}}{m(\La)}
\,\bigr\}
\Eqno\equadefSrB
$$
for some $r\in\,]0,e^{-3\al}[$.

\Proc{Theorem 2.2}
{For all $r\in\,]0,e^{-3\al}[$ the \BWD\ series of the form~$\fS_\cR(a)$, $a\in\cS(r,B)$,
belong to~$\cM((K_{j}),B)$.
In particular, this the case for the general solution~$F_{r_1,r_2}$
if $B=\cL(B_{r_1},B_{r_2})$ and $0<r_2 < r_1\,e^{-3\al}$.
}

\proof
According to Definition~2.3 we must check that $f=\fS_\cR(a)\in \cC^1_{hol}(K_j,B)$
for every $a\in\cS(r,B)$ and~$j\in\N$.
It is natural to define the function
$$
f^{(1)}(q)=-\sum_{m=1}^\infty \sum_{\Lambda\in{\cal R}^{*}_{m}}
{a_\La\over (q-\La )^{2}}
$$
whose restriction to $\INT(K_{j})$ is just the ordinary 
derivative of~$f$. 

According to Lemma~2.5,
$$
\forall q\in K_j,\;\forall\La\in\cR,\quad
|q-\La |\ge \mu\ga_j {e^{-\alpha m(\La )}\over m(\La )}. 
\Eqno\ineqDistLa
$$
Thus, for $k=0$ or~1, and for $q\in K_j$,
$$
\Vert f^{(k)}(q)\Vert 
\le 
\sum_{m=1}^\infty \sum_{\La\in\cR^*_m} \frac{c \, r^m}{|q-\La |^{k+1} m}
\le
c(\mu\ga_j)^{-k-1} \sum_{m=1}^\infty m^{k+1} (r\,e^{(k+1)\al})^m <+\infty.
$$
Note that $f$ and~$f^{(1)}$ are continuous since the convergence is uniform 
and $K_{j}$ is compact. To prove ${\cal C}^{1}$-smoothness,
we consider the remainder
$$
R(q,q') = {f(q)-f(q')\over q-q'} -f^{(1)}(q')
= -\sum_{m=1}^\infty \sum_{\Lambda\in{\cal R}^{*}_{m}}
{(q'-q)\over (q-\La )(q'-\La )^{2}} a_\La.
$$
Because of~\ineqDistLa\ and the assumption $r<e^{-3\alpha}$,
we have $\Vert R(q,q')\Vert \le c_j |q-q'|$, with
$c_j= c(\mu\ga_j)^{-3}
      \sum_{m=1}^\infty m^{3} (r\,e^{3\al})^m$.
In particular Definition~2.1 is satisfied.

The statement about~$F_{r_1,r_2}$ is a particular case of what we just proved:
choosing $a_\La = \La \cL_{m(\La)} \odot$ and $r=\frac{r_2}{r_1}$,
we use Lemma~A1.1 and see that
$\Vert a_\La \Vert_{\cL(B_{r_1},B_{r_2})} \le \Vert\cL_{m(\La)}\Vert_{B_r}
\simeq \frac{1}{m(\La)} r^{m(\La)}$.

\qed

As for the fundamental solution, notice that $f_\de \in \cM((K_j),B_r)$
as soon as $0<r<e^{-3\al}$.

\beginsection{2.5 Whitney smoothness of monogenic functions}

As already mentioned in Remark~2.4, we cannot apply Borel's Theorem 
to conclude that functions in ${\cal M}((K_{j})_{j\in\N},B)$
are ${\cal C}^\infty$-holomorphic in some subsets of the~$K_{j}$'s.
But this can be shown directly.

Let $c_{0}(\R )$ denote the classical Banach space of real sequences $s=(s_{k})_{k\ge 0}$ 
such that $s_{k}\rightarrow 0$ as $k\rightarrow +\infty$, 
endowed with the norm $\Vert s\Vert = \sup |s_{k}|$. 
A subset~$A$ of~$c_{0}(\R )$ is {\sl closed} and {\sl totally bounded} 
if and only if the following two conditions are satisfied: 
\item{(i)} $\exists C>0 \;/\; \forall s\in A,\;\Vert s\Vert \le C.$
\item{(ii)} $\forall\eps >0,\;\exists k_{0}\in \N \;/\;
             \forall s\in A,\; \forall k\ge k_{0},\; |s_{k}|\le \eps$.

\Def{Definition 2.6}
{To $\ga\in\,]0,1[$ and~$A$, totally bounded closed subset of~$c_{0}(\R )$,
we associate 
$$
W^{A}_{\ga} = \{ x\in {\R}\setminus {\Q}\modZ\, \mid\,
\exists s\in A \text{such that} \forall k\in\N,\;
m_{k+1}(x)\le \ga^{-1} e^{s_{k}m_{k}(x)} \}.
$$
If moreover $\ka\in\,]0,1[$ and $d>0$, we define
$$
W^{A}_{\ga,\ka,d} = 
  \bigcup_{y\in W^{A}_{\ga}}
  \{ x\in {\C/\Z}\, \mid \; \ka |\RE(\ti x-\ti y)| \le |\IM\ti x| \le d \},
$$
where $\ti x$ and $\ti y$ denote some lifts in~$\C$ of~$x$ and~$y$.
}

One can study the sets~$W^A_{\ga}$ and~$W^A_{\ga,\ka,d}$ with the same kind of
arguments as in Section~2.3.
For instance one can easily check that they are closed and perfect.
Notice that $W^A_{\ga}$ is non-empty as soon as~$A$ contains a 
sequence~$s$ such that $s_{k}\ge 2 G^{3-k\over 2}$ for all~$k$
(indeed $g\in W^A_{\ga}$ in that case).
%
%
Moreover, if $x\in \R\setminus \Q \modZ$ satisfies the condition
$$
\lim_{k\to\infty} \frac{\log m_{k+1}(x)}{m_{k}(x)}=0,
\Eqno\equacondarlin
$$
and if $\al>0$ is given,
there exist~$\ga\in\,]0,1[$ and~$s\in c_0(\R)$ such that 
$x\in W^{\{s\}}_\ga$ and $\Vert s\Vert\le \alpha$.

 %
 %



\Proc{Theorem 2.3}
{Let $\ga,\ka\in\,]0,1[$, $d>0$, $\psi$ an approximation function 
of the form $\psi (m)= \ga e^{-\alpha m}$ 
and $K = \{ q=e^{2\pi ix}, \; x\in C_{\psi,\ka,d} \}$.
Let $A$ a totally bounded closed subset of~$c_{0}(\R )$ such that
$\forall s\in A$, $\Vert s\Vert \le \al$,
and $K^* = \{ q=e^{2\pi ix}, \; x\in W^{A}_{8\ga,\ka,d/2} \}$.
Then $K^*\subset K$ and ${\cal C}^{1}_{hol}(K,B)\subset {\cal C}^{\infty}_{hol}(K^*,B)$
for any Banach space~$B$.
}

\proof
It is immediate to check that 
$W^{A}_{8\ga}\subset C_{\psi}=\{x\,\mid\, \forall k\in\N,\;
m_{k+1}(x)\le\ga\ii\,e^{\al m_k(x)} \}$; thus~$K^*\subset K$. 
\smallbreak

Let $f\in{\cal C}^{1}_{hol}(K,B)$. We will use Remark~2.1.
Observe that, in view of Lemma~2.4,
the connected components of~$(\C/\Z)\setminus C_{\psi,\ka,d}$ are
of the form $\De_{n/m}$ with $n/m\in\Q_\psi$,
except for one or two of them: the components of~$i\infty$ and~$-i\infty$ may
be reduced to the half-planes $\{ \pm\IM x>d \}$, or else they both coincide
with the union of these half-planes and a finite number of diamonds~$\De_{n/m}$.
>>From that we deduce the decomposition 
$\bigsqcup_{\ell\ge1} U_\ell$ of~$\C\setminus K$ into connected components
--- the index $\ell=1$ (\resp $\ell=1$ and~2) will correspond to the exceptional component
(\resp components), the next ones being numbered as 
$U_\ell = \exp(2\pi i\De_{n_\ell/m_\ell})$ with a non-decreasing sequence~$(m_\ell)$.

Moreover, for each $n/m\in\Q_\psi$, 
we recall that according to Proposition~2.2,
$$
|X-\frac{n}{m}|<r_{n/m}=\frac{2\ga}{m} \, e^{-\al m} 
\quad\text{if $X=x_{n/m}$ or~$x'_{n/m}$},
$$
hence $\pa\De_{n/m}$ has length less than~$4 r_{n/m}\sqrt{1+\ka^2}$.
The series $\sum\length(\pa U_\ell)$ is thus convergent.
\smallbreak

Let $j\in\N$. We will now check that the series
$$
\sum_{\ell\ge1} \int_{\pa U_\ell} \frac{|d\ze|}{|\ze-q|^{j+1}} 
\Eqno\equaSerLengj
$$
is uniformly convergent for $q\in K^*$. This will allow us to set
$$
f^{(j)}(q) = \frac{j!}{2\pi i}\sum_{\ell\ge1} \int_{\pa U_\ell} 
             \frac{f(\ze)}{(\ze-q)^{j+1}}
             \,d\ze.
\Eqno\equaDerivj
$$

\Proc{Lemma 2.7}
{There exists a positive number~$\mu$ (which depends only on~$\ka$) such that,
whenever $n/m\in\Q_\psi$,
$$
\forall\xi\in \ov\De_{n/m}, \quad
\dist(e^{2\pi i\xi},K^*) > \frac{2\mu\ga}{m} \, e^{-\al m}.
$$
For each $j\in\N$, there exists a positive integer~$M$ 
(which depends only on~$\ga$,~$\al$ and~$j$) such that,
whenever $n/m\in\Q_\psi$ and $m\ge M$,
$$
\forall\xi\in \ov\De_{n/m}, \quad
\dist(e^{2\pi i\xi},K^*) > \frac{2\mu\ga}{m} \, e^{-\frac{\al m}{2(j+1)}}.
\Eqno\ineqDistRac
$$
}

We end the proof of Theorem~2.3 before proving Lemma~2.7.
%
According to the first part of Lemma~2.7, each term in the series~\equaSerLengj\
is well defined when $q\in K^*$.
For~$\ell$ large enough (say $\ell\ge L$), $U_\ell = \exp(2\pi i\De_{n_\ell/m_\ell})$ with $n_\ell/m_\ell\in\Q_\psi$
and $m_\ell\ge M$, thus we can use~\ineqDistRac\ for each~$q\in K^*$:
$$
\int_{\pa U_\ell} \frac{|d\ze|}{|\ze-q|^{j+1}} 
\le 2\pi\,e^{2\pi d} \Bigl(\frac{m_\ell}{2\mu\ga}\Bigr)^{j+1} e^{\frac{\al m_\ell}{2}} 
                                                              \length(\pa\De_{n_\ell/m_\ell})
\le \frac{8\pi\,e^{2\pi d}\sqrt{1+\ka^2}}{\mu} \Bigl(\frac{m_\ell}{2\mu\ga}\Bigr)^{j} 
                                                              e^{-\frac{\al m_\ell}{2}}.
$$
The series~\equaSerLengj\ is thus convergent, and we can use~\equaDerivj\ with
$j=0$ or~1 to represent~$f$ or~$f^{(1)}$ in~$K^*$.
For $j\ge2$, we define~$f^{(j)}$ in~$K^*$ by~\equaDerivj, and
the previous computation shows the existence of~$C>0$ such that
$$
\forall\ell\ge L,\; \forall q\in K^*,\quad
\Vert \int_{\pa U_\ell} \frac{f(\ze)}{(\ze-q)^{j+1}} \,d\ze \Vert
\le C\, m_\ell^j\, e^{-\frac{\al m_\ell}{2}}
$$
(and for $\ell < L$ this expression is continuous in~$q$);
hence, by uniform convergence, $f^{(j)}$ is continuous in~$K^*$.

Let us consider the Taylor remainders
$$
R^{(j,v)}(q,q') = f^{(j)}(q') - \sum_{u=0}^v \frac{1}{u!} f^{(j+u)}(q) (q'-q)^u
$$
for $j,v\ge0$ and $q,q'\in K^*$.
Remark~2.5 applies also to~$W^{A}_{8\ga,\ka,d/2}$, and thus to~$K^*$:
these sets have a finite number of connected components and are locally
connected.
In fact, for $q,q'\in K^*$ close enough (say $|q-q'|\le\de$), one can define a
path~$\Ga(q,q')$ which joins~$q$ to~$q'$ inside~$K^*$ and which is the image by
$x\mapsto\,e^{2\pi ix}$ of the union of 1,2 or~3 segments of slopes~$\pm\ka$;
the length of~$\Ga(q,q')$ is less than $\nu|q'-q|$, where $\nu$ depends only on~$\ka$.

We now conclude the proof of Theorem~2.3 by checking that there exists 
$C>0$ such that
$$
\forall q,q'\in K^*,\quad |q-q'|\le\de \ens\Rightarrow\ens 
\Vert R^{(j,v)}(q,q')\Vert\le C \, |q'-q|^{v+1}.
\Eqno\ineqconcllem
$$
We can write 
$$
R^{(j,v)}(q,q') = \frac{j!}{2\pi i} \sum_{\ell\ge1} \int_{\pa U_\ell}
\cR^{(j,v)}(q,q',\ze) f(\ze) \,d\ze,
$$
where $\cR^{(j,v)}(q,q',\ze)$ is the Taylor remainder at order~$v$
for the function $q'\mapsto (\ze-q')^{-j-1}$, \ie
$$
\cR^{(j,v)}(q,q',\ze) = \frac{(j+v+1)!}{j!\,v!} \int_{\Ga(q,q')} 
\frac{(q'-q'')^v}{(\ze-q'')^{j+v+2}} \,dq''.
$$
>>From this identity and from Lemma~2.7 applied with~$j$ replaced by~$j+v+1$, one
can deduce the existence of a positive integer~$L$ such that,
if $\ze\in\pa U_\ell$ with $\ell\ge L$,
$$
\Vert \cR^{(j,v)}(q,q',\ze) \Vert \le \text{const} m_\ell^{j+v+2} 
\, e^{\frac{\al m_\ell}{2}} \, |q-q'|^{v+1},
$$
whereas for $\ell<L$, 
$\Vert \cR^{(j,v)}(q,q',\ze) \Vert \le \text{const}|q-q'|^{v+1}$.
Therefore, the validity of~\ineqconcllem\ follows from the inequalities
$\length(\pa U_\ell)\le\frac{\hbox{\sevenrm const}}{m_\ell} \,e^{-\al m_\ell}$.
\qed

\Pf{Proof of Lemma~2.7}
We must show that $|q-e^{2\pi i\xi}|>\text{const}\frac{\psi(m)}{m}$ for $q\in K^*$
and $\xi\in\ov\De_{n/m}$.
Notice that $\frac{\psi(m)}{m}\le\ga\,e^{-\al}\le\demi$ since $\psi$ is an
approximation function, and
$|\IM\xi|\le\ka(x'_{n/m}-x_{n/m})<\ka$.
Therefore we can assume 
$$
q=e^{2\pi ix} \quad\text{with}\ens x\in W^A_{8\ka,\ka,d}, \ens |\IM x|\le 2\ka.
$$
Moreover we can consider that $|\RE z|\le\demi$, where $z=x-\xi$, and since
$|q-e^{2\pi i\xi}|\ge \,e^{-4\pi\ka}\,|1-e^{-2\pi iz}|$,
it will be enough to bound from below~$|z|$ itself
($z$ lies indeed in a domain where $|(e^{-2\pi iz}-1)/z|$ is bounded from below).
The same reasoning holds for the proof of~\ineqDistRac\ provided that we take
$m\ge M\ge 2(j+1)$.

In fact we will prove the inequalities
$$
\forall n/m\in\Q_\psi,\; \forall y\in W^A_{8\ga},\quad
|y-\frac{n}{m}| \ge \frac{4\ga}{m}\,e^{-\al m},
\Eqno\ineqpflemun
$$
and the existence, for each $j\in\N$, of a positive integer~$M$ such that
$$
\forall n/m\in\Q_\psi,\; \forall y\in W^A_{8\ga},\quad
m\ge M \ens\Rightarrow\ens
|y-\frac{n}{m}| \ge \frac{4\ga}{m}\,e^{-\frac{\al m}{2(j+1)}}.
\Eqno\ineqpflemdeux
$$
This is enough to bound $|z|=|x-\xi|$ from below as required since
for any~$x\in W^A_{8\ga,\ka,d}$ there
exists~$y\in W^A_{8\ga}$ such that $|\IM x|\ge\ka|\RE(x-y)|$, but
then $|x-\xi|\ge (1+\ka^2)^{-1/2}\dist\bigl(y,[x_{n/m},x'_{n/m}]\bigr)$ 
for all $\xi\in\ov\De_{n/m}$,
and 
$$\tst
\dist\bigl(y,[x_{n/m},x'_{n/m}]\bigr) \ge 
|y-\frac{n}{m}| - \max\bigl(x'_{n/m}-\frac{n}{m},|x_{n/m}-\frac{n}{m}|\bigr)
< |y-\frac{n}{m}| - \frac{2\ga}{m}\,e^{-\al m}.
$$
\smallskip

Let $y\in W^A_{8\ga}$. Let $s\in A$ such that
$m_{k+1}(y)\le\frac{1}{8\ga}\,e^{s_k m_k(y)}$.
According to~\ineqfc,
$$
\forall k\ge0,\quad |y-\frac{n_k(y)}{m_k(y)}| 
> \frac{4\ga}{m_k(y)}\, e^{-s_k m_k(y)}
\ge \frac{4\ga}{m_k(y)}\, e^{-\al m_k(y)}.
$$
Let $n/m\in\Q_\psi$.
Either $\frac{n}{m}=\frac{n_k(y)}{m_k(y)}$ for some $k\ge0$ and \ineqpflemun\ is proved.
Or $\frac{n}{m}$ is not a convergent of~$y$; then
$m_{k-1}(y) \le m < m_k(y)$ for some $k\ge1$ and Proposition~A3.4 applies:
$$
m\,|y-\frac{n}{m}| > m_{k-1}(y)\,|y-\frac{n_{k-1}(y)}{m_{k-1}(y)}| 
\ge 4\ga \, e^{-s_{k-1}m_{k-1}(y)} \ge 4\ga \, e^{-\al m}.
$$
Therefore \ineqpflemun\ is true in all cases.

As for~\ineqpflemdeux, given~$j\in\N$ we first choose $k_0\ge0$ such that 
$|s_k| \le \frac{\al}{2(j+1)}$ for all~$k\ge k_0$ and~$s\in A$.
We then choose $M\ge1$ such that 
$$
\forall y\in W^A_{8\ga}, \quad m_{k_0+1}(y) < M.
$$
(For all $y\in W^A_{8\ga}$, 
$m_0(y) = 1$ thus $m_1(y)\le \frac{1}{8\ga}\,e^\al = M_1,\;
m_2(y)\le \frac{1}{8\ga}\,e^{\al M_1} = M_2, \ldots$: take $M>M_{k_0+1}$.)
According to~\ineqfc, we have now
$$
\forall y\in W^A_{8\ga},\ens \forall k\ge k_0, \quad |y-\frac{n_k(y)}{m_k(y)}| 
> \frac{4\ga}{m_k(y)}\, e^{-\frac{\al m_k(y)}{2(j+1)}}.
$$

Let $y\in W^A_{8\ga}$ and $n/m\in\Q_\psi$ with $m\ge M$.
We are faced with the same alternative as above, but we know moreover that 
if $\frac{n}{m} = \frac{n_k(y)}{m_k(y)}$
or $m_{k-1}(y) \le m < m_k(y)$, necessarily $k\ge k_0+1$.
Therefore we obtain the refined inequality~\ineqpflemdeux\ in all cases.
\qed

\Def{Definition 2.7}
{For any closed totally bounded subset~$A$ of~$c_0(\R)$ and any integer~$j$, we define
$$
K_{A,j}^* = \{ q = e^{2\pi ix}, \; x\in W^A_{8\ga_j,\ka,d/2} \}
$$
provided that $\Vert s\Vert\le\al$ for all $s\in A$,
with the same notations as in Definition~2.5.
}

According to Theorem~2.3, $K_{A,j}^* \subset K_j$ and 
${\cal C}^{1}_{hol}(K_j,B)\subset {\cal C}^{\infty}_{hol}(K_{A,j}^*,B)$.
In particular, according to Theorem~2.2, the solutions of the cohomological
equation are $\cC^\infty$-holomorphic in each~$K_{A,j}^*$.

Observe that any point of the form $\la=e^{2\pi ix}$ with~$x$ satisfying~\equacondarlin\
lies in~$K_{\{s\},j}^*$ for~$s$ well chosen and~$j$ large enough.


\vfill \eject

\beginsection{3. Carleman classes at Diophantine points}

In this section, we address the following question (directly inspired by~[He],
Question at p.\ 82):
Do the solutions of the cohomological equation belong to any quasianalytic
Carleman class?
We will treat separately some particular points of~$\S^1$ 
among those at which Theorem~2.3
yields Whitney smoothness, 
and study asymptotic expansions in disks tangent to~$\S^1$ at each of these points.

As a preliminary, in Section~3.1, we define the Carleman classes~$\CM\pm$ which
we think are the most relevant for the problem at 
hand\footnote{\noteCar}
{Carleman classes are usually defined as spaces of functions which are defined
and $C^\infty$ (in the real sense) on some --- possibly infinite --- interval~$I$
of~$\R$ and whose derivatives satisfy some uniform bounds (see~[Th]); the relationship
between such classes with $I=\R^+$ and the spaces~$\CM\pm$ defined below is
indicated in~[Ca].
}.  
We recall a well-known criterium of quasianalyticity due to Carleman, and we
also introduce spaces of functions which admit Gevrey asymptotic expansions.
Our presentation is somewhat influenced by the works of Ramis and Malgrange on
divergent series (see for instance~[Ra], [Ma]).

In Section 3.2 we prove that all functions monogenic in the compacts~$K_j$ of
Definition~2.5 admit Gevrey-$\tau$ asymptotic expansions at Diophantine points
of exponent~$\tau\ge 2$.
%
%
On the other hand, 
in the case of the fundamental solution,
we prove in Section~3.3 the sharpness of the
index~$\tau=2$ in Gevrey asymptotics for those Diophantine points which
correspond to quadratic irrationals, and conclude that no quasianalytic Carleman
class at those points contains the fundamental solution.


\beginsection{3.1 Carleman and Gevrey classes}

{\bf a)}
Let $B$ be a complex Banach space, whose norm we denote by $\nor\,.\,\nor$,
and $\la\in\S^1$.
Let us fix some sequence $\{M_n\}_{n\ge0}$ of positive numbers.

\Def{Definition~3.1}{We define the Carleman class $\CM-$
to be the vector space of all $B$-valued functions~$f$ for which there exist
an open disk $\De\subset\D$ tangent to~$\S^1$ at~$\la$, 
a formal series $\sum_{n\ge0} a_n Q^n \in B[[Q]]$
and positive numbers~$c_0$ and~$c_1$ such that
the function~$f$ is holomorphic in~$\De$ and
$$
\forall N\ge0,\; \forall q\in\De, \quad
\nor f(q) - \sum_{0\le n\le N-1} a_n (q-\la)^n \nor \le c_0\, c_1^N\, M_N\, |q-\la|^N.
$$
}

\vskip .5cm
\epsfysize=2.1cm
\centerline{\epsfbox{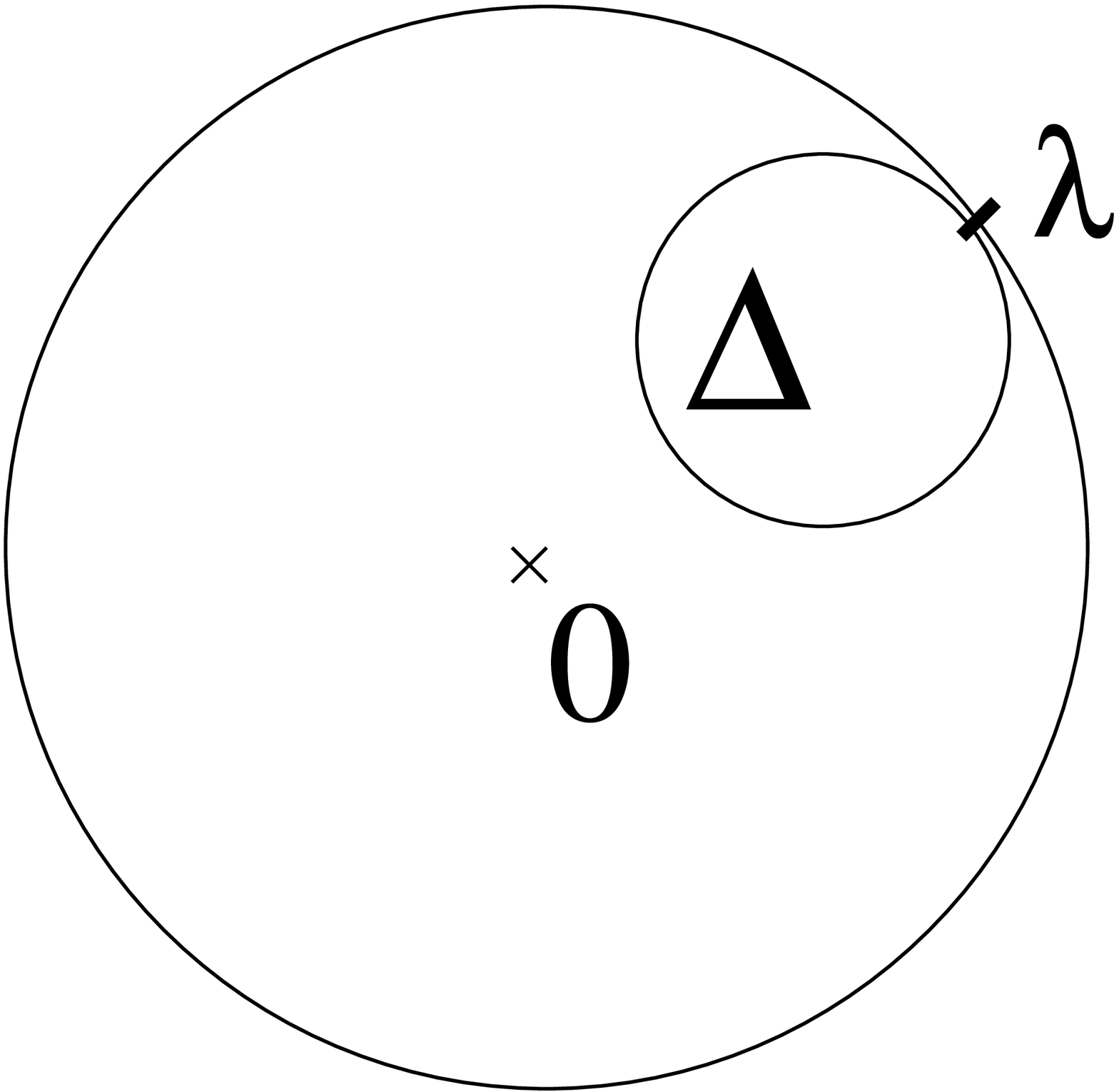}}
\vskip .3cm

The mapping
$$
J^-_\la: \; f \in \CM-  \ens\mapsto\ens \sum_{n\ge0} a_n Q^n \in B[[Q]] 
$$
which associates to a function of~$\CM-$ its asymptotic expansion at~$\la$ is well defined.
In fact the functions of~$\CM-$ are $\cC^\infty$-holomorphic 
in the sense of Definition 2.2:

\Proc{Lemma~3.1}
{If $f\in\CM-$ and $S\subset\D$ is a bounded sector of vertex~$\la$ and small
enough radius,
there exist positive numbers~$c_0$ and~$c_1$ such that
\item{--} the function $f$ is $\cC^\infty$-holomorphic in the closure~$\bar S$ of the sector,
\item{--} for all $n\ge0$ and $q\in\bar S$,
$\nor \frac{1}{n!} f^{(n)}(q) \nor \le c_0 \, c_1^n \, M_n$,
\item{--} and $J^-_\la(f)=\dst\sum_{n\ge0}\frac{1}{n!} f^{(n)}(\la)\,Q^n$.

Conversely, if a function~$f$ is $\cC^\infty$-holomorphic in a closed
disk $\bar\De\subset\overline\D$ tangent to~$\S^1$ at~$\la$
and satisfies inequalities of the form 
$\nor \frac{1}{n!} f^{(n)}(q) \nor \le c_0 \, c_1^n \, M_n$
in~$\bar\De$, then it belongs to~$\CM-$.
}

\noindent
Here, by ``bounded sector'' we mean the intersection of an open infinite sector and of some open
disk centered at its vertex, and we call ``radius'' of the sector the radius of
that disk.

We leave the proof of Lemma~3.1 to the reader (one can use the
Taylor-Lagrange formula).

As a consequence, the asymptotic expansion~$J^-_\la(f)$ of
any $f\in\CM-$ belongs to the space
$$
B[[Q]]_{\{M_n\}} = \biggl\{\, \sum_{n\ge0} a_n Q^n \in B[[Q]]\mid
\exists c_0,c_1>0 \;\text{such that}\; (\forall n\ge0)\,
\nor a_n \nor \le c_0\,c_1^n \,M_n \,\biggr\}.
$$
By definition, the space~$\CM-$ is quasianalytic at~$\la$ if and only if the mapping~$J^-_\la$ 
is injective.


Analogously we define the space~$\CM+$ by using disks~$\De$ contained in~$\E$ instead
of~$\D$, and the corresponding mapping 
$J^+_\la:\, \CM+ \rightarrow  B[[Q]]_{\{M_n\}}$.
The change of variable $q\mapsto\la^2/q$ induces an isomorphism between~$\CM-$ and~$\CM+$.

We can now state Carleman's criterium of quasianalyticity~[Ca]:

\Proc{Theorem 3.1 (Carleman's Criterium)}
{The space $\CM\pm$ is quasianalytic at~$\la$ 
(\ie $J^\pm_\la$ is injective on that space)
if and only if 
$\dst\;\sum_{n\ge1} \frac{1}{\be_n} = +\infty\;$
where $\dst\be_n=\inf_{n'\ge n} \{\, M_{n'}^{1/{n'}}\}$.
}

\remark{3.1}{The criterium is usually stated for spaces of scalar functions, but it is also
valid for spaces of $B$-valued functions (as soon as $B\neq\{0\}$ of course).
The quasianalyticity of~$\CM\pm$ is indeed equivalent to that of~$\cC^\pm(\la,\{M_n\},\C)$
because of the existence of non-trivial continuous linear functionals on any normed linear
space: if $f$ is a function in~$\CM\pm$, any continuous linear functional~$\ell$ on~$B$ 
induces a function $\ell\circ f$ which belongs to~$\cC^\pm(\la,\{M_n\},\C)$, and
$J_\la^\pm(\ell\circ f) = \ell\bigl(J_\la^\pm(f)\bigr)$ (letting $\ell$ act termwise
in~$B[[Q]]$ in order to define the right-hand side). The point is that for a function~$f$ to
be identically zero, it is necessary and sufficient that all the functions~$\ell\circ f$
vanish identically
(given any Banach space $B$ its dual separates points on $B$).
}

Let~$\CMs$ be the space of all $B$-valued functions for
which there exist disks $\De^-\subset\D$ and $\De^+\subset\E$ tangent to~$\S^1$ at~$\la$ such
that $f_{|\De^-}\in\CM-$, $f_{|\De^+}\in\CM+$ and
$J^-_\la(f_{|\De^-}) = J^+_\la(f_{|\De^+})$.
We will denote by $J_\la(f)$ simply the asymptotic expansion at~$\la$ of a function~$f$
of~$\CMs$.
As a consequence of Theorem~3.1, $\CMs$ is quasianalytic at any point
of~$\D\cup\{\la\}\cup\E$ if and only if $\;\dst\sum \frac{1}{\be_n} = +\infty$.


\bigbreak\noindent {\bf b)}
As a special case we will consider {\sl Gevrey classes}, \ie spaces of functions with
Gevrey-$\tau$ asymptotic expansion for some~$\tau\ge0$.

\Def{Definition~3.2}
{If $B$ is a Banach space, $\la\in\S^1$ and $\tau\in\R^+$, we define the Gevrey classes
$$
\cG_\tau^-(\la,B),\qquad \cG_\tau^+(\la,B),\qquad \cG_\tau(\la,B)
$$
respectively as the Carleman classes
$$
\CM-,\ens \CM+,\ens \CMs
$$ 
with the sequence $\{ M_n = \Ga(1+n\tau) \}$.
We also set
$B[[Q]]_\tau = B[[Q]]_{\{M_n\}}$ with the same sequence~$\{M_n\}$.
}

We warn the reader that not all the authors follow this convention for indexing Gevrey
classes.
For us, $\tau=0$ corresponds to the analytic class:
$B[[Q]]_0$ is the space~$B\{Q\}$ of convergent series, and $J^-_\la$ and~$J^+_\la$ are
isomorphisms in that case. 
Thus $\cG_0^\pm(\la,B)$ and~$\cG_0(\la,B)$ can all be identified to 
the space of all germs of $B$-valued holomorphic functions at~$\la$.

We retain that, by Carleman's Theorem, the space $\cG^\pm_\tau(\la,B)$ is quasianalytic
at~$\la$ if and only if~$\tau\le1$; and the same is true for $\cG_\tau(\la,B)$.
One can check that,
if $B$ is a Banach algebra,
$\cG_\tau^\pm(\la,B)$ and $\cG_\tau(\la,B)$ are in fact algebras: they are stable
by multiplication~[Ma].

\bigbreak\noindent {\bf c)}
We will now focus on the $\tau=1$ case and the relationship with the Laplace
transform. We suppose moreover that $B$ is a Banach algebra.

We denote by~$\hat\cN^\pm(B)$ the space of all $B$-valued functions~$\hat\phi$ for
which there exist some positive numbers~$\rho'<\rho$ and some real number~$\de$ such that
$\hat\phi$ is holomorphic in the open ``half-strip''
$$
H^\pm_{\rho} = \ao \xi \in\C\,/\ens \dist(\xi,\R^\pm) < \rho \af
$$
and $\xi \,\mapsto\, e^{-\de|\xi|}\,\nor \hat\phi(\xi) \nor$ is bounded in the closed
half-strip~$\bar H^\pm_{\rho'}$. The vector space~$\hat\cN^\pm(B)$ is stable by convolution,
the convolution of two holomorphic functions~$\hat\phi_1$ and~$\hat\phi_2$ being defined as
$\hat\phi_1*\hat\phi_2(\xi) = \int_0^\xi \hat\phi_1(\xi_1)\hat\phi_2(\xi-\xi_1) \, d\xi_1$.

We also introduce a symbol~$\de_0$ which one may think of as the Dirac distribution at
the origin:
identifying any pair $(a_0,\hat\phi)\in B\times\hat\cN^\pm(B)$ with the symbolic sum
$a_0\de_0+\hat\phi\in B\de_0\oplus\hat\cN^\pm(B)$
and extending the convolution to the space~$B\de_0\oplus\hat\cN^\pm(B)$ by treating $\de_0$
as a unit, we get an algebra.
The following theorem is due to Nevanlinna [Ma]:

\Proc{Theorem 3.2 (Nevanlinna's Theorem)}
{The Laplace transform 
$$
\L_\la^\pm\,:\;
a_0\de_0+\hat\phi \;\mapsto\; f^\pm \;\text{such that}\; 
f^\pm(\la(1+t)) = a_0+\int_0^{\pm\infty} \hat\phi(\xi)\,e^{-\xi/t}\,d\xi
$$
defines an isomorphism between the algebras $B\de_0\oplus\hat\cN^\pm(B)$
and~$\cG^\pm_1(\la,B)$.
}

\remark{3.2}
{Again we mention that the replacement of scalar functions by $B$-valued functions, with
respect to the usual statement, is innocuous.
Notice that Nevanlinna's Theorem implies that~$\cG_1^\pm(\la,B)$ is a differentiable
algebra: it is stable by derivation
(see Appendix~A.5 for a description of the counterpart
in the {convolutive model}~$\hat\cN^\pm(B)$ of such elementary operations as
differentiation). 
Also, with respect to the notations of Definition~3.1, 
we have incorporated in our statement the change of infinitesimal variable
$Q=q-\la \,\mapsto\, t=\la\ii Q$ 
in order to deal with Laplace integrals on~$\R^\pm$ only
(the counterpart in~$\hat\cN^\pm$ of such homotheties
and of more general changes of variable is described in Appendix~A.5).
}

Theorem~3.2 shows that the quasianalyticity of~$\cG^\pm_1(\la,B)$ is in some sense
constructive, the reciprocal operator of~$J^\pm_\la$ being described in terms of {\sl
Borel-Laplace resummation}:

\Def{Definition~3.3}
{If $\ti f=\sum a_n Q^n\in B[[Q]]_1$, we define a formal series 
$\ti\phi(t) = \sum \phi_n t^n\in B[[t]]_1$ 
by
$\ti\phi(t)= \ti f(\la\,t)$, and its {\sl formal Borel transform} by
$\phi_0\de_0+\hat\phi$ where
$$
\hat\phi(\xi) = \sum_{n\ge0} \phi_{n+1} \frac{\xi^n}{n!}  
              = \sum_{n\ge0} \la^{n+1} a_{n+1} \frac{\xi^n}{n!}\in B\{\xi\} ;
$$
the series $\ti f$ belongs to~$J^\pm_\la\bigl(\cG^\pm_1(\la,B)\bigr)$ if and only if
$\hat\phi$ can be analytically continued to an element of~$\hat\cN^\pm(B)$, and its preimage 
is then equal to~$\L_\la^\pm(a_0\de_0+\hat\phi)$: it is called the {\sl Borel-Laplace sum}
of~$\ti f$ (in the direction of~$\R^\pm$). }

The reader is referred to Appendix~A.5 for more details on the Borel-Laplace summation
process.

\Def{Definition~3.4}
{Let $\hat\cN(B) = \hat\cN^-(B)\cap\hat\cN^+(B)$.
We define $\L_\la$ in~$B\de_0\oplus\hat\cN(B)$ by gluing 
$\L_\la^-$ and~$\L_\la^+$:
we obtain an isomorphism between~$B\de_0\oplus\hat\cN(B)$ and~$\cG_1(\la,B)$.
}


\beginsection{3.2 Gevrey asymptotics at Diophantine points for monogenic functions}

Let $B$ a Banach space.
According to Theorem~2.3, monogenic functions of~$\cM((K_j),B)$ 
are ${\cal C}^\infty$-holomorphic in the compacts~$K^*_{A,j}$, 
with the notations of Definitions~2.5 and~2.7.
In particular, such a function admits as asymptotic expansion its Taylor series 
at any point of~$K^*_{A,j}\cap {\S}^1$.
Among those points, some of them have further arithmetic properties which
will yield Gevrey asymptotic expansions.

\Def{Definition 3.5}
{Let $\ga>0$, $\tau\ge 2$. We define $\DC(\ga,\tau)$ to be the set of all irrational
numbers~$y$ which satisfy Diophantine inequalities of constant~$\ga$ and
exponent~$\tau$, \ie
$$
\forall n/m\in{\Q}, \quad 
|y - n/m|\ge \ga \, m^{-\tau}. 
$$
We also set $\DC_\tau = \bigcup_{\ga >0} \DC(\ga ,\tau )$ 
and $\oDC_\tau = \{ \la = e^{2\pi iy},\; y\in\DC_\tau \}$.
}

It is well-known that $\DC_\tau$ has full measure as soon as~$\tau>2$
and that $\DC_2$ (which has measure zero) coincides with the set of
constant-type irrationals (irrationals with bounded quotients).
If~$y\in\DC(\ga,\tau)$, the property 
$$
\forall k\ge0,\quad
m_{k+1}(y) < \ga\ii \,m_k(y)^{\tau-1}
$$ 
allows one to find~$A$ and~$\ga'$ such that $y\in W^A_{\ga'}$
(\eg $A=\{s\}$ with $s_k = m_k(y)^{-\de}$ for some $\de\in]0,1[$, and
$\ga' = \ga\,\min_{_{k\ge0}}\{m_k(y)^{1-\tau}\exp(m_k(y)^{1-\de})\}$).
In particular, each point of~$\oDC_\tau$ is contained in some~$K^*_{A,j}$.

\Proc{Theorem 3.3}{Let $\tau\ge2$. 
If $\la\in\oDC_\tau$, monogenic functions of~$\cM((K_j),B)$ admit Gevrey-$\tau$
asymptotic expansions at~$\la$:
$$
{\cal M}((K_{j}),B)\subset {\cal G}_\tau (\la ,B ).
$$
}

In particular, according to Theorem~2.2,
the solution~$F_{r_1,r_2}$ of the cohomological equation belongs 
to~$\cG_\tau(\la,\cL(B_{r_1},B_{r_2}))$ as soon as $0<r_2<r_1$
(using the fact that the positive number~$\al$ which enters into the definition of the
sequence~$(K_j)$ can be chosen arbitrarily small).
Similarly $f_\de \in \cG_\tau(\la,B_r)$ if $0<r<1$.

The proof of Theorem~3.3 is somewhat analogous to that of Theorem~2.3. We first state a lemma
about the relation between Diophantine points and the geometry of the compacts~$K_j$,
which parallels Lemma~2.7.

\Proc{Lemma 3.2}{Let $\tau\ge 2$ and $\la\in\oDC_\tau$.
There exist $\mu>0$, $j\ge1$ and two open disks~$\De^-\subset{\D}$ and~$\De^+\subset{\E}$ 
tangent to~${\S}^1$ at~$\la$ such that
the set $\De^- \cup \{\la\} \cup \De^+$ is contained in~$K_j$ and,
for every~$n/m\in\Q_{\psi_j}$ and~$\xi \in \ov\De_{n/m}$,
the point~$\ze=e^{2\pi i\xi}$ satisfies
$$
\dist(\ze,\ov{\De^- \cup \De^+}) \ge \mu |\ze-\la|^2
\quad\text{and}\quad
|\ze-\la| \ge \mu\,m^{-\tau}.
\Eqno\ineqlem
$$
}

\proof
Let $\ga>0$ and $y\in\DC(\ga,\tau)$ such that $\la=e^{2\pi iy}$.
We choose~$j$ large enough to ensure 
$\ga_{j} \le \frac{1}{4}\ga\,\min_{_{m\ge1}}\{m^{1-\tau}\,e^{\al m}\}$.
According to Definition~3.5, 
$$
\forall n/m\in\Q,\quad
|y-n/m| \ge\tst \ga_{j}\,\frac{e^{-\al m}}{m} = \frac{\psi_{j}(m)}{m},
$$
hence $y\in C_{\psi_{j}}\modZ$ and $\la\in K_{j}$ by~\excltoocl.

Let us define the function $f(X) = 2\ka\ga_{j}\,\exp(-c\,|X|^{-1/\tau})$, 
with $c=\al(\frac{\ga}{2})^{1/\tau}$.
We can use Lemma~2.4 and Proposition~2.2 to show that
$$
\cK_f=\{\xi\in\C \,\mid\; |\IM\xi|\ge f(\RE(\xi-y))\} \subset \ti C_{\psi_{j},\ka},
$$
where $\ti C_{\psi_{j},\ka}$ denotes the lift of~$C_{\psi_{j},\ka}$ in~$\C$.
Indeed, if $\xi\in\C\setminus\ti C_{\psi_{j},\ka}$, there exists $n/m\in\Q_\psi$
such that $\xi\in\De_{n/m}\modZ$;
according to Proposition~2.2~(2), 
$$
|\RE(\xi-n/m)| < 2\ga_{j}\, \tst \frac{e^{-\al m}}{m}
\quad\text{and}\quad
|\IM\xi| \le \demi\ka(x'_{n/m}-x_{n/m}) < \dst 2\ka\ga_{j}\,e^{-\al m};
$$
but $X=\RE(\xi-y)$ satisfies $|X| \ge |y-n/m| - |\RE(\xi-n/m)| \ge \demi\ga\,m^{-\tau}$,
hence $|\IM\xi|<f(X)$.

Since~$\cK_f$ has a contact of infinite order with~$\R$ at~$y$, 
we obtain $\De^-\cup\De^+\subset\exp(2\pi i \, C_{\psi_{j},\ka})$
by taking the radius of these disks small enough.
Reducing this radius if necessary, we make them contained in the 
annulus $\{ e^{-2\pi d} \le |q| \le e^{2\pi d} \}$ and thus in~$K_{j}$.

Finally, by compactness, it is sufficient to prove~\ineqlem\ 
for~$\ze=e^{2\pi i\xi}$ close to~$\la$.
On the one hand, the estimate
$$
\dist(\ze,\ov{\De^-\cup\De^+}) 
\mathrel{\mathop\sim_{\ze\to\la}}
\text{const} |\ze-\la|^2
$$
follows from the fact that, for all $\ze\in\C\setminus \exp(2\pi i \,\cK_f)$,
$|\ze| = 1 + \ti f(\ze-\la)$,
where the function~$\ti f(X)$ is exponentially small for small~$|X|$.
On the other hand,
$|\ze-\la|\ge\text{const}|\xi-y|\ge\text{const}m^{-\tau}$ 
if~$\xi\in\ov\De_{n/m}$,
according to the previous computation.
%
\qed

\remark{3.3}
{The exponent~``2'' in the right-hand side of the first inequality of~\ineqlem\ 
corresponds to the order of contact of the disks~$\De^\pm$ 
in which we ask for asymptotic expansions.
But the proof of Theorem~3.3 which follows would be valid with any other
exponent as well. This means that a monogenic function of~$\cM((K_j),B)$ admits
a Gevrey-$\tau$ asymptotic expansion at~$\la$ in compacts with arbitrarily high
order of contact at~$\la$, not only disks.}

\Pf{Proof of Theorem 3.3}
Let $\tau\ge 2$, $\la\in\oDC_\tau$ and $f\in\cM((K_j),B)$.
Let~$\mu$, $j$, $\De^\pm$ as in Lemma~3.2.
We proceed as in the proof of Theorem~2.3:
the connected components of~$\C\setminus K_j$ are of the form
$U_\ell = \exp(2\pi i\,\De_{n_\ell/m_\ell})$ with $n_\ell/m_\ell\in\Q_{\psi_j}$,
except for~$\ell=1$ or~$\ell=1,2$.
We recall that $\length(\De_{n/m})\le\text{const}\frac{e^{-\al m}}{m}$.
Formula~\equaDerivj\ leads us to define the coefficients
$$
a_k = {f^{(k)}(\la)\over k!} = 
\frac{1}{2\pi i}\sum_{\ell\ge1}
\int_{\pa U_\ell}\frac{f(\ze)}{(\ze-\la)^{k+1}}\,d\ze,
\qquad
k\ge0.
$$
Cauchy's formula (extended to monogenic functions) applies for $q\in\ov{\De^-\cup\De^+}$:
$$
f(q) = \frac{1}{2\pi i}\sum_{\ell\ge1}
\int_{\pa U_\ell}\frac{f(\ze)}{\ze-q}\,d\ze.
$$
Using the identity 
$$
\frac{1}{\ze-q} = 
                  \sum_{k=0}^{N-1}\frac{(q-\la)^k}{(\ze-\la)^{k+1}} 
                  + \frac{(q-\la)^N}{(\ze-\la)^N(\ze-q)},
$$
we find that 
$$
f(q)-\sum_{k=0}^{N-1}a_k(q-\la )^k = 
\frac{1}{2\pi i} \sum_{\ell\ge1}
\int_{\pa U_\ell} \frac{f(\ze)}{(\ze-\la)^{N}}\frac{(q-\la)^N}{\ze-q}\,d\ze.
$$
We now use~\ineqlem\ to bound the contributions of the 
components~$U_\ell = \exp(2\pi i\,\De_{n_\ell/m_\ell})$,
noticing that, by compactness, such inequalities hold for the exceptional components as well 
provided that~$\mu$ is small enough:
if $\ze\in\pa U_\ell$, $|\ze-q|\ge\mu|\ze-\la|^2$ and
$|\ze-\la|\ge\mu\,m_\ell^{-\tau}$
(extending the definition of~$m_\ell$ by the value~1 for the exceptional
components),
hence
$$
\Vert f(q)-\sum_{k=0}^{N-1}a_k(q-\la )^k\Vert \le 
\frac{\text{const}}{\mu^{N+3}} \sum_{\ell\ge1} 
                               \frac{e^{-\al m_\ell}}{m_\ell}m_\ell^{(N+2)\tau}
\le \frac{\text{const}}{\mu^{N+3}} \Phi((N+2)\tau),
\Eqno\ineqtruncated
$$
with $\Phi(X) = \sum_{m\ge1} m^X\,e^{-\al m}$.
Comparing the sum~$\Phi(X)$ and the integral
$\int_0^{+\infty} m^X\,e^{-\al m}\,dm = \al^{-X-1}\Ga(X+1)$,
we obtain
$\Phi(X) \le \al^{-X-1}(\Ga(X+1) + 2 \al \,X^X\,e^{-X})$ and the Stirling formula
yields the result.
\qed


\beginsection{3.3 Borel transform at quadratic irrationals for the fundamental solution} 

We fix in this section $r\in\,]0,1[$.
If $\la\in\oDC_\tau$ for some $\tau\ge2$, according to Theorem~3.3
the solutions of the cohomological equation are contained in the corresponding Gevrey
class, which is not quasianalytic at~$\la$. 
But would it be possible for them to be contained in some smaller, quasianalytic Carleman class?
We now show that the answer is {\em negative} if $\tau=2$ and $\la$ belongs to a
subset~$\IQ$ of~$\oDC_2$. 

\Def{Definition~3.6}
{For any point $\la=e^{2\pi i \al}$ in~$\oDC_2$ (say with $\al\in\,]0,1[$),
we define the {\sl Lagrange spectral constants}~$\nu_\pm(\la)>0$ by
$$
\frac{1}{\nu_-(\la)} = - \liminf_{(D,N)\in\N^*\times\Z} \ao D^{-2}(\frac{N}{D} - \al)\ii
\af,
\quad
\frac{1}{\nu_+(\la)} =   \limsup_{(D,N)\in\N^*\times\Z} \ao D^{-2}(\frac{N}{D} - \al)\ii \af.
$$
We will use the notation $\ka_\pm(\la) = (\nu_\pm(\la))^{1/2}$ too.
}

\Def{Definition~3.7}
{We define~$\IQ$ to be the subset of~$\oDC_2$ consisting of all $\la=e^{2\pi i \al}$ with
$\al$ quadratic irrational, \ie $\al\in\R\setminus\Q$ algebraic of degree~2.}

The {\sl Lagrange spectrum} can be defined as the set $\ao \nu(\la) = \min\{
\nu_-(\la),\nu_+(\la) \}, \; \la\in\oDC_2 \af$
(\ie it is the set of the numbers $\nu(\la) = \liminf \{ D^2|\frac{N}{D}-\al| \}$ 
for $\la\in\oDC_2$), 
but here we need an asymmetric version of it
because we will separate the rational approximations of~$\al$ by the left from its rational
approximations by the right.
We will need the restriction $\la\in\IQ$ because
of the following lemma, which is an arithmetical result about the way the quantities 
$D^2 |\frac{N}{D} - \al|$ approach~$\nu_\pm(\la)$, and for which we do not know of any
analogue when $\la\in\oDC_2\setminus\IQ$.

\Proc{Lemma~3.3}
{Let $\al\in\,]0,1[$ be irrational and algebraic of degree~2. Let $\la=e^{2\pi i\al}$ and
$$
\cE^-=\ao (D,N)\in\N^*\times\Z \,|\ens N/D < \al \af
\quad\text{and}\quad
\cE^+=\ao (D,N)\in\N^*\times\Z \,|\ens N/D > \al \af.
$$
For each of these sets, there exist a partition 
$$
\cE^\pm = \cF^\pm \cup \cE^\pm_* \cup \cA^\pm 
$$
and a number $\ka_\pm'>\ka_\pm(\la)$ such that:
\pppar
-- the set $\cF^\pm$ is finite;
\pppar
-- for all $(D,N)\in\cE^\pm_*$, $D^2|\frac{N}{D}-\al| \ge (\ka_\pm')^2$;
\pppar
-- the set $\cA^\pm$ can be written 
$$
\cA^\pm = \ao (D_p^\pm,N_p^\pm),\; p\ge0 \af
$$
with $\{D_p^\pm\}$ increasing sequence of $\N^*$,
$\; \sum (D_p^\pm)^{-1/2} < \infty$
and
$
(D_p^\pm)^2|\frac{N_p^\pm}{D_p^\pm} - \al| = \nu_\pm(\la) + o(\frac{1}{D_p^\pm}).
$
\pppar \noindent
Moreover, if $\eps\in\{+,-\}$ satisfies $\ka_\eps(\la)\le\ka_{-\eps}(\la)$,
the sequence $\{D_{p+1}^\eps/D_p^\eps\}$ is bounded.
}

Notice that
$$
\nu_\pm(\la) = \ka_\pm(\la)^2 = \liminf_{(D,N)\in\cE^\pm} \ao D^2 |\frac{N}{D} - \al|\af
             = \lim_{p\rightarrow\infty}  (D_p^\pm)^2 |\frac{N_p^\pm}{D_p^\pm} - \al|.
$$
In the case of the golden mean $\al=\frac{1+\sqrt{5}}{2}$, one may check that
$\nu_+(\la)=\nu_-(\la)=\frac{1}{\sqrt{5}}$. But for $\al=\sqrt{3}$, one finds
$\nu_-(\la)=\frac{1}{\sqrt{3}} > \nu_+(\la)=\frac{1}{2\sqrt{3}}$.
In both examples one can take the even convergents for the
sequence~$\{\frac{N_p^+}{D_p^+}\}$ and the odd convergents for the
sequence~$\{\frac{N_p^-}{D_p^-}\}$.

The proof of Lemma~3.3 is given in Appendix~A.4, since it is purely
arithmetical.
It is the only place where we use the hypothesis $\la\in\IQ$.

\Proc{Theorem~3.4 (Non-quasianalyticity and sharpness of Gevrey-2 asymptotics for quadratic
irrationals)}
{Let $\la\in\IQ$. We know by Theorem~3.3 that $f_\de\in\cG_2(\la,B_r)$, thus we may consider
its asymptotic expansion at~$\la$:
$$
\ti f = J_\la(f_\de) = \sum_{n\ge0} F_n Q^n \in B_r[[Q]]_2,
$$
and the formal Borel transform of~$Q^{1/2} \ti f(Q)$ with respect to~$Q^{1/2}$:
$$
\hat F(\xi,z) = \sum_{n\ge0} F_n(z) \frac{\xi^{2n}}{(2n)!} \in \C\{\xi,z\}.
$$
(a)
The holomorphic germ~$\hat F$ extends analytically to the set
$\ao (\xi,z)\in\C\times\D_r\,/\ens \xi\in \REC(z) \af$ where, for each $z\in\D_r$,
the rectangle $\REC(z)$ is defined as the set of the complex numbers~$\xi$ such
that
$$
|\RE\bigl((2\pi i\la)^{-1/2}\xi\bigr)| < \ka_+(\la) \log\frac{1}{|z|}
\ens\text{and}\ens
|\IM\bigl((2\pi i\la)^{-1/2}\xi\bigr)| < \ka_-(\la) \log\frac{1}{|z|}.
$$
(b)
For each $z\in\D_r$, $\pa(\REC(z))$ is a natural boundary for the analytic function~$\xi
\mapsto
\hat F(\xi,z)$. 
\ppar \noindent
(c) 
Suppose that $\{M_n\}$ is a non-decreasing sequence of positive numbers such
that 
$$
f_\de\in\CMr- \ens\text{or}\ens f_\de\in\CMr+.
$$
Necessarily $\CMr\pm$ is not quasianalytic at~$\la$.
}


\vskip .15cm
\epsfysize=6cm
\centerline{\epsfbox{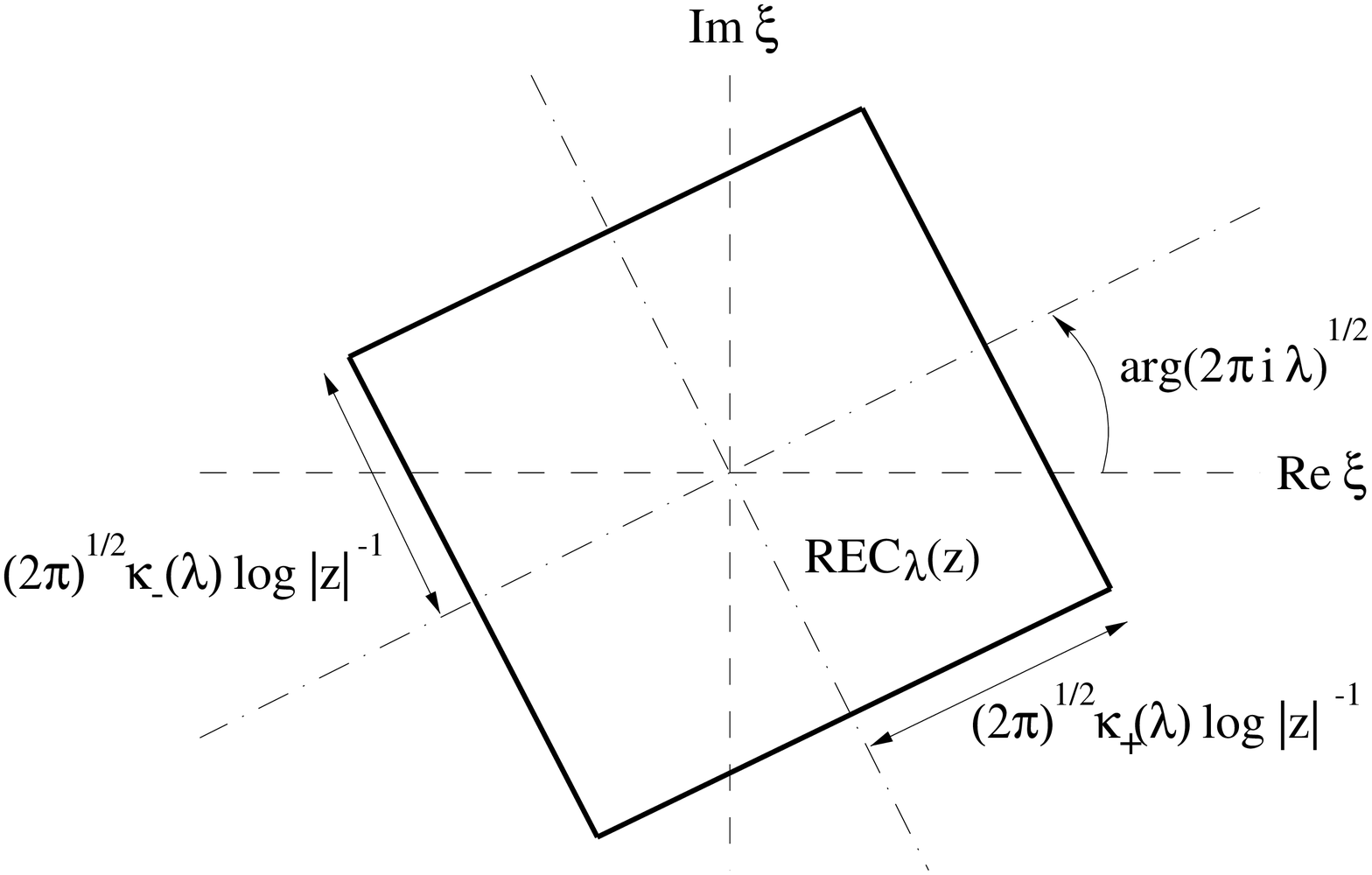}}
\vskip .3cm

The assumption that $\{M_n\}$ be non-decreasing seems only technical, but we were
not able to get rid of it.
With that restriction the spaces of solutions $\{f_g,\; g\in B_{r'}\}$ with $r'>r$ 
and {\sl a fortiori} the spaces of monogenic functions $\cM(\{C_j\},B_{r'})$ 
are contained in none of the quasianalytic Carleman classes 
at~$\la$ that we have defined in Section~3.1.

Note that this theorem holds for the fundamental solution of the cohomological equation,
because of its very specific features, 
but we claim no such result for a general \BWD\ series with poles in~$\cR$ nor for any class
of monogenic functions.

We will obtain that theorem itself as a consequence of a more precise result.
In the statement of this result, we will make use of the variables
$h=\frac{1}{2\pi i}\log\frac{q}{\la}$ and $s=\log z$ rather
than~$Q=q-\la$ and~$z$. 
Since we are dealing with functions of~$B_r$, the variable~$s$
will move in the half-plane $\ao\RE s<\log r<0\af$ and these functions decrease 
at least like~$e^{\RE s}$ when $\RE s$ tends to~$-\infty$.

\Proc{Theorem~3.5 (Borel transform of order~2 at quadratic irrationals)}
{Let $\la\in\IQ$. One can give a decomposition of the fundamental solution
$$
f_\de(\la \,e^{2\pi i h},e^s) = f_\de(\la,e^s) + \frac{1}{2\pi i}\bigl(
\chi^+(h,s) + \chi^-(h,s)
\bigr)
$$
satisfying the following properties:
\ppar\noindent
(a)
the function $\chi^\pm$ is analytic for $\RE s<\log r$ and $h\in\C\setminus\R^\pm$, with
$$
\chi^+(h,s) = h^{1/2} \int_0^{+i\infty} \hat\psi^+(\ze,s) \,e^{-\ze h^{-1/2}} \, d\ze,
\quad
\chi^-(h,s) = h^{1/2} \int_0^{+\infty} \hat\psi^-(\ze,s) \,e^{-\ze h^{-1/2}} \, d\ze,
$$
the Borel transform~$\hat\psi^\pm$ being analytic in
$$
\ao (\ze,s) \,/\ens \RE s<\log r \;\text{and}\; |\RE\ze| <
\ka_+(\la) (-\RE s) \af \ens\text{for $\hat\psi^+$},
$$
$$
\ao (\ze,s) \,/\ens \RE s<\log r \;\text{and}\; |\IM\ze| <
\ka_-(\la) (-\RE s) \af \ens\text{for $\hat\psi^-$},
$$ 
and, for each~$s$,
even with respect to~$\ze$ and bounded in any 
substrip
$$
\left.\eqalign
{%
\ao |\RE\ze| \le \text{const} \!\af \ens&\text{for $\hat\psi^+$}\;\cr
\ao |\IM\ze| \le \text{const} \!\af \ens&\text{for $\hat\psi^-$}\;
}\right|
\quad \text{const} < \ka_\pm(\la) (-\RE s) ;
$$
\ppar\noindent
(b)
for each~$s$, the Borel transform~$\ze\mapsto\hat\psi^\pm(\ze,s)$ has a dense set of
singular points on the boundary of its strip of definition;
more precisely, if one defines the points
$$
\ze^+_{k,l}(s) = \ka_+(\la) \bigl( -s + 2\pi i (k\al+l) \bigr),
\quad
\ze^-_{k,l}(s) = i \ka_-(\la) \bigl( -s + 2\pi i (k\al+l) \bigr),
\qquad
k,l\in\Z,
$$
the real part of the function~$\hat\psi^\pm(\ze,s)$
tends to~$-\infty$ when $\ze$ tends to one of the points~$\ze^\pm_{k,l}(s)$,
horizontally from the left for~$\hat\psi^+$,
vertically from below for~$\hat\psi^-$
(\ie $\ze=\ze^\pm_{k,l}(s) + \xi$, $\xi\rightarrow0$, with
$\xi\in\R^-$ for~$\hat\psi^+$ and $\xi\in i\,\R^-$ for~$\hat\psi^-$).
\ppar\noindent
(c)
for each real $s<\log r$ there exists a positive integer~$j_0$ and
a non-decreasing sequence of positive
numbers $\{\de_j\}_{j\ge j_0}$ such that 
$$
\sum_{j\ge j_0} (\de_j)^{-3/4} < +\infty \qquad\text{and}\qquad
\forall j\ge j_0, \quad
|\chi_{2j-1}(s)| \ge (\de_j)^{2j-1},
$$
with the following notation for the Taylor series of~$\hat \psi=\hat\psi^++\hat\psi^-$:
$$
\hat\psi(\ze,s) = \sum_{n\ge0} \chi_{n+1}(s) \frac{\ze^{2n}}{(2n)!}.
\Eqno{\defichin}
$$
}

\vskip .15cm
\epsfysize=6cm
\centerline{\epsfbox{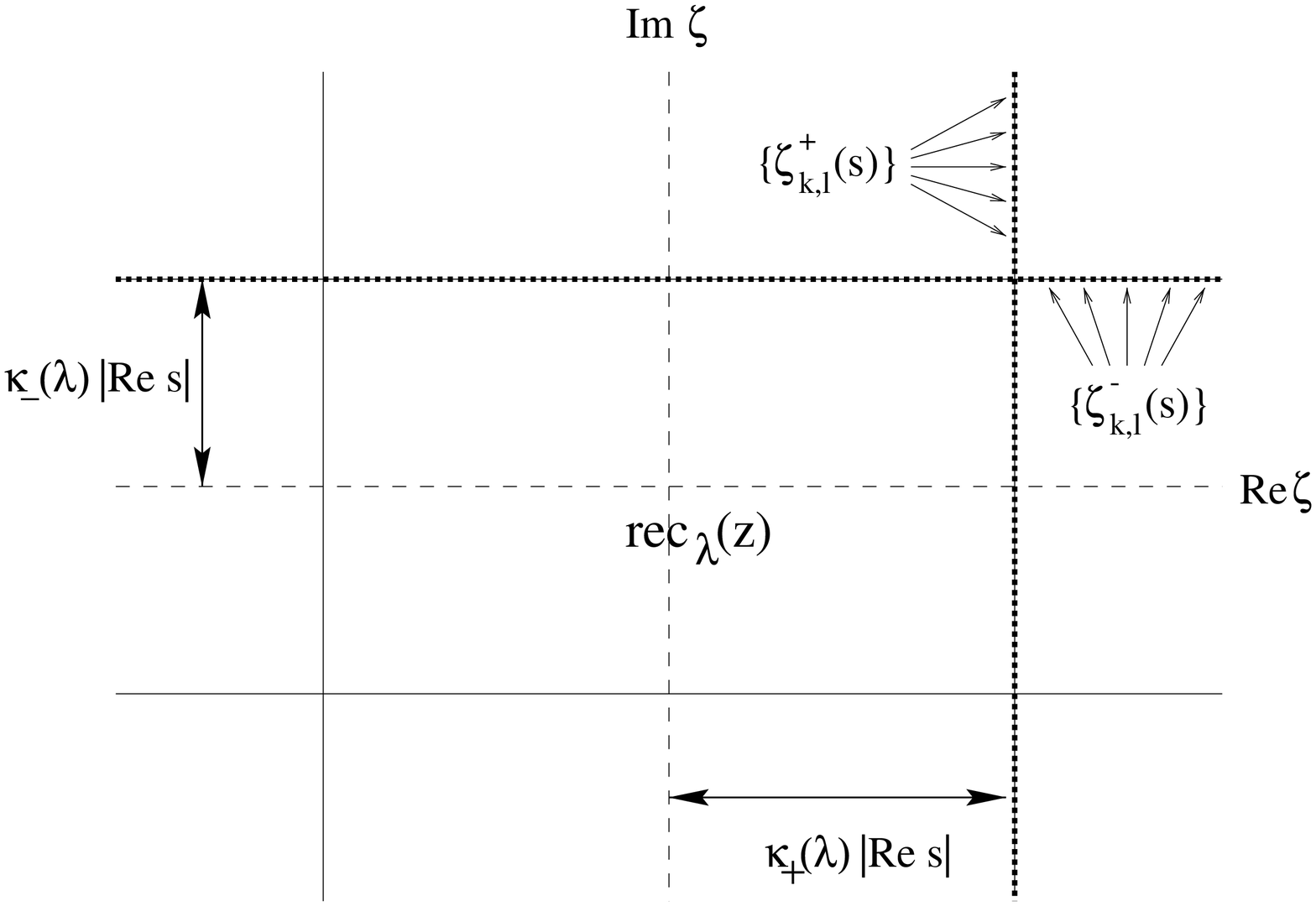}}
\vskip .3cm

\remark{3.4}
{The function~$\ze\mapsto\hat\psi^\pm(\ze,s)$ is in fact the {\sl integral} Borel transform
of~$h\mapsto h^{-1/2}\chi^\pm(h,s)$ with respect to~$h^{1/2}$, whereas its Taylor series
at~$\ze=0$ is the {\sl formal} Borel transform with respect to~$h^{1/2}$ of the asymptotic
expansion of~$h^{-1/2}\chi^\pm$ at~$h=0$. In the formulas of Part~(a) which indicate how to
recover~$\chi^\pm$ from~$\hat\psi^\pm$ by Laplace transform, 
there is an implicit choice of determination of~$h^{1/2}$: for~$\chi^+$ one chooses the
determination which is holomorphic in~$\C\setminus\R^+$ and has always positive imaginary
part, while for~$\chi^-$ one chooses the
determination which is holomorphic in~$\C\setminus\R^-$ and has always positive real
part (in order to ensure the decrease of~$|e^{-\ze h^{-1/2}}|$).
}

\vskip .3cm
\epsfysize=7.5cm
\centerline{\epsfbox{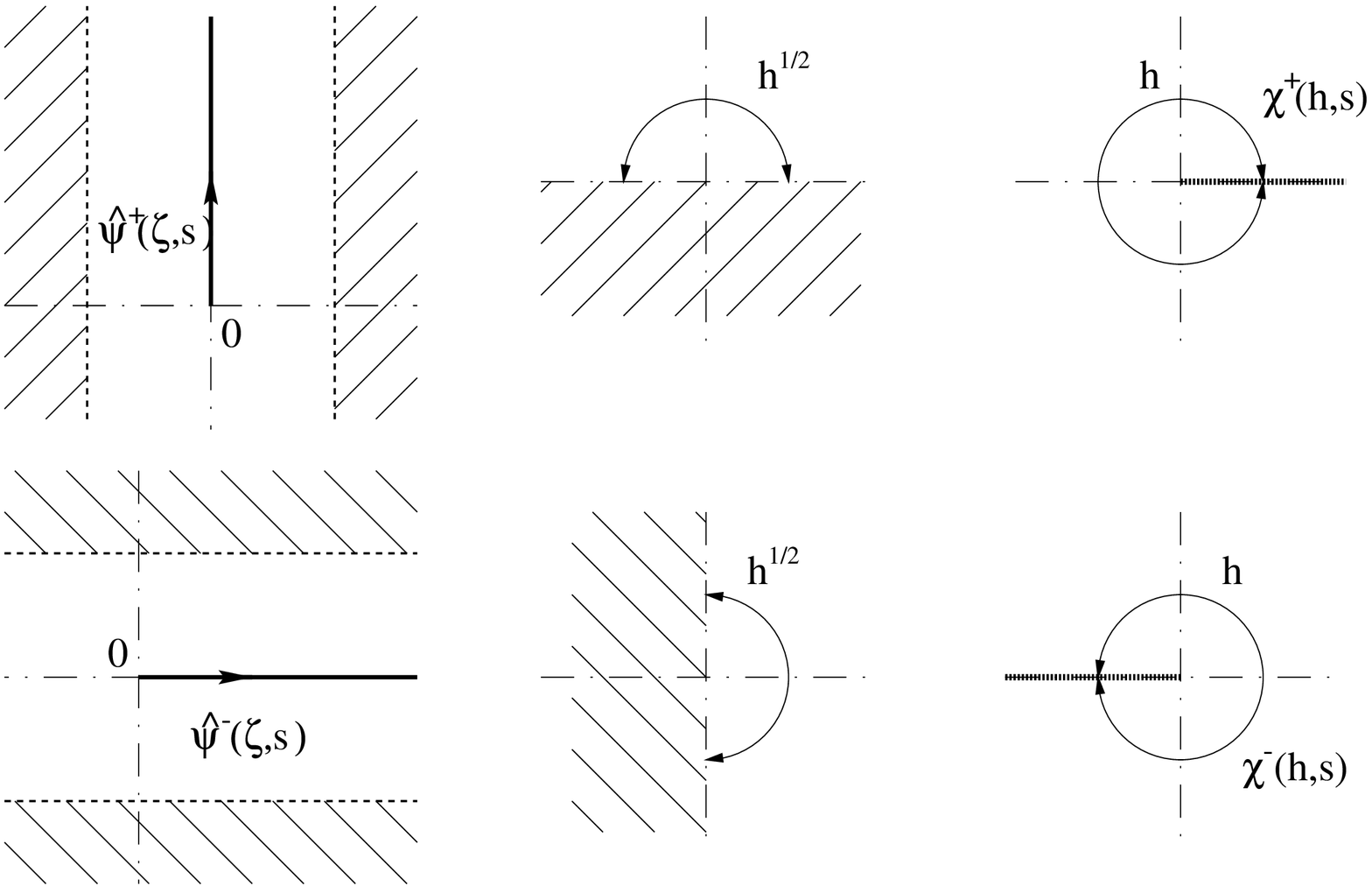}}
\vskip .3cm

\noindent
One could make the opposite choice as well: 
since $\hat\psi^\pm$ is even with respect
to~$\ze$, one would simply have to compute the Laplace 
integral along the opposite ray.
Besides, Parts~(a) of Theorem~3.4 or~3.5 do not require $\la\in\IQ$ but only~$\la\in\oDC_2$.


\beginsection{3.4 Deduction of Theorem~3.4 from Theorem~3.5}

Parts~(a) and~(b) are an exercise of application of the general 
theory of which Appendix~A.5 gives a brief account. 
We will relate $\hat F(\xi)$ and
$$
\hat\psi(\ze) = \hat\psi^+(\ze) +
\hat\psi^-(\ze)
\Eqno{\defipsi}
$$ 
(from now on the variable~$z=e^s$ will be understood).
Part~(a) of Theorem~3.5 implies that $\hat\psi$ is analytic in the 
rectangle~$\rec(z)$
defined by
$$
|\RE\ze| < \ka_+(\la) (-\RE s) 
\ens\text{and}\ens
|\IM\ze| < \ka_-(\la) (-\RE s).
$$

Let $Q_1=Q^{1/2}$ and $\ti F(Q_1) = Q_1 \ti f(Q_1^2)$: 
$\hat F(\xi)$ is the formal Borel transform of~$\ti F$ with 
respect to~$Q_1$, which we will
indicate by the notation
$$
\ti F = \LF{\xi}{Q_1} \hat F
$$
in order to be able to deal with changes of variables in the formal model.
By definition of~$\ti f$, we have the asymptotic expansion
$f_\de(\la+Q_1^2) \sim Q_1\ii \ti F(Q_1)$, thus
$$
f_\de(\la+Q_1^2) \sim \hat F(0) + \LF{\xi}{Q_1}(\pa_\xi\hat F).
$$

On the other hand we can introduce $h_1=h^{1/2}$. 
According to Part~(a) of Theorem~3.5, 
$f_\de(\la \,e^{2\pi i h_1^2}) 
\sim f_\de(\la) + \frac{1}{2\pi i} h_1 \LF\ze{h_1} \hat\psi 
   = f_\de(\la) + \LF\ze{h_1} (\frac{1}{2\pi i} * \hat\psi)$.
We deduce that $\hat F(0)=f_\de(\la)$, and
setting
$$
\hat G_1 = \frac{1}{2\pi i} * \hat\psi,
\Eqno{\defiGun}
$$
we have the identity
$
\LF{\xi}{Q_1}(\pa_\xi\hat F) = \LF\ze{h_1} \hat G_1
$
under the change of variable
$$
h_1 =\bigl[ \frac{1}{2\pi i}\log(1 + \la\ii Q_1^2) \bigr]^{1/2} 
    = (2\pi i \la)^{-1/2} Q_1 (1+O(Q_1^2)).
$$

This change of variable is the composition of the dilatation 
$h_1 \mapsto Q_2 = (2\pi i \la)^{1/2} h_1$ 
and of the transformation 
$Q_1 \mapsto Q_2 = \bigl[ \la \log(1+\la\ii Q_1^2)\bigr]^{1/2}$.
The dilatation is responsible for the passage from $\hat G_1$ analytic for~$\ze\in\rec(z)$
to a function
$$
\hat G_2(\xi_2) = (2\pi i \la)^{-1/2} \hat G_1((2\pi i \la)^{-1/2} \xi_2)
\quad\text{analytic for~$\xi_2\in\REC(z)=(2\pi i \la)^{1/2}\rec(z)$,}
\Eqno{\defiGdeux}
$$
such that
$
\LF\ze{h_1} \hat G_1 = \LF{\xi_2}{Q_2} \hat G_2
$.
According to Part~(b) of Theorem~3.4, $\pa\REC(z)$ is a natural boundary for~$\hat G_2$.

Finally 
$$
\hat F = \hat F(0) + 1*\hat G,
$$ 
where the function~$\hat G(\xi)$ is determined
from~$\hat G_2$ 
by {\sl composition-convolution}: indeed
$$
\LF\xi{Q_1} \hat G = \LF{\xi_2}{Q_2} \hat G_2
$$
under a change of variable
$$
Q_2\ii = Q_1\ii + L_{12}(Q_1) \ens \Leftrightarrow\ens
Q_1\ii = Q_2\ii + L_{21}(Q_2), \qquad
L_{12}(X),L_{21}(X) \in X\C\{X\},
$$
hence
$$
\hat G = \hat G_2 + \sum_{r\ge1} \frac{1}{r!} (\hat L_{12})^{*r} * \hat\pa^r \hat G_2,
\quad
\hat G_2 = \hat G + \sum_{r\ge1} \frac{1}{r!} (\hat L_{21})^{*r} * \hat\pa^r \hat G,
\Eqno{\defiG}
$$
where $\hat\pa$ denotes  the 
multiplication by~$-\xi$ or~$-\xi_2$, Borel
counterpart of differentiation with respect to~$X_1=Q_1\ii$ or~$X_2=Q_2\ii$.
Here $\hat L_{12}$ and~$\hat L_{21}$ are entire functions and $\REC(z)$ is star-shaped with
respect to the origin, hence the above series are uniformly convergent in any compact
subset of~$\REC(z)$.
Therefore $\hat G$ is holomorphic in~$\REC(z)$, and if $\pa\REC(z)$ were not a natural
boundary for~$\hat G$, neither would it be for~$\hat G_2$.
This proves the statements of Parts~(a) and~(b) of Theorem~3.4.

\medbreak

As for Part~(c), we now suppose that
$f_\de\in\CMr\pm$ for some non-decreasing sequence of positive numbers~$\{M_n\}$.
In particular $\ti f = J_\la^\pm(f_\de) \in B_r[[Q]]_{\{M_n\}}$.
Let us fix a real number $s_0<\log r$ at which all the subsequent $s$-dependent
functions will be evaluated.
For instance $F_n$ will denote the value at~$s_0$ of the function $s\mapsto F_n(e^s)$, and we
have
$$
\forall n\ge1,\quad |F_n| \le c_0 c_1^n M_n
$$
for some $c_0,c_1>0$.

Part~(c) of Theorem~3.5 yields a sequence~$\{\de_j\}$ which allows to bound from 
below half of the coefficients of
$$
\ti \chi(h) = \sum_{n\ge0} \chi_n h^n,
$$
where we use the convention of~\defichin\ for denoting the Taylor coefficients
of~$\hat\psi$ and
we set $\chi_0=2\pi i f_\de(\la)$ for conveniency.
We have
$$
f_\de(\la+Q) \sim \ti f(Q) = \sum_{n\ge0} F_n Q^n
\quad\text{and}\quad
f_\de(\la \,e^{2\pi i h}) \sim \frac{1}{2\pi i} \ti\chi(h),
$$
therefore $\ti\chi(h) = 2\pi i \ti f(\la(e^{2\pi ih}-1))$
and in particular, for all $n\ge1$,
$$
\chi_n = (2\pi i)^{n+1} \sum_{r=1}^n \la^r b_{r,n} F_r
\qquad
\text{with}
\quad
b_{r,n} = \sum_{n_1+\cdots+n_r=n,\,n_i\ge1} \frac{1}{n_1!\cdots n_r!}.
$$
Since $b_{r,n} < r^n/n!$ and $\{M_n\}$ is non-decreasing, we deduce that
$$
|\chi_n| \le c_0 (2\pi)^{n+1} M_n \frac{n^{n+1}}{n!} (\max\{1,c_1\})^n ,
$$
and for $n$ large enough, $M_n \ge c_2^n |\chi_n|$ with some $c_2>0$.

We are now in a position to apply Theorem~3.1: let
$\dst\be_n=\inf_{n'\ge n} \{\, M_{n'}^{1/{n'}}\}$.
For $j$ large enough,
$$
M_{2j} \ge M_{2j-1} \ge c_2^{2j-1} (\de_j)^{2j-1},
$$
thus, for $j$ large enough,
$$
M_{2j}^{\frac{1}{2j}} \ge c_2^{1-\frac{1}{2j}} (\de_j)^{1-\frac{1}{2j}}
\ge \text{const} (\de_j)^{3/4}
\quad\text{and}\quad
M_{2j-1}^{\frac{1}{2j-1}} \ge c_2 \de_j \ge \text{const} (\de_j)^{3/4}.
$$
Since the sequence~$\{\de_n\}$ is non-decreasing, 
$\be_{2j} \ge \text{const} (\de_j)^{3/4}$,
$\be_{2j-1} \ge \text{const} (\de_j)^{3/4}$,
thus $\sum \frac{1}{\be_n} < +\infty$.
\qed


Let $\la=e^{2\pi i\al}\in\IQ$ with $\al\in\,]0,1[$.
Before proceeding to the proof of Theorem~3.5,
we give a ``decomposition into simple elements'' of~$f_\de$ with respect to
the variable $h = \frac{1}{2\pi i}\log \frac{q}{\la}$.

\Proc{Lemma~3.4}
{$$
f_\de(\la\,e^{2\pi ih},z) = f_\de(\la,z) + \frac{1}{2\pi i}
\sum_{D\in\N^*,N\in\Z} Z\ii \cdot \frac{h}{h-Z} \cdot \frac{z^D}{D},
\quad \text{with}\ens Z = \frac{N}{D}-\al.
$$}

\Pf{Proof of Lemma~3.4}
We start with the decomposition which is relative to the variable~$q$
and which can be written
$$
f_\de(q,z) = \sum_{D\ge1} \frac{z^D}{D} \sum_{\La\in\cR_D} (\frac{q}{\La} - 1)\ii.
$$
Thus
$$
f_\de(\la\,e^{2\pi ih},z) = \sum_{D\ge1} \sum_{0\le N\le D-1} 
          \frac{z^D}{D} (e^{2\pi i(h+\al-\frac{N}{D})} - 1)\ii.
$$
We now use the identity
$$
\frac{d}{dx} (e^{2\pi ix}-1)\ii =
                        \frac{1}{2\pi i}\sum_{M\in\Z} \frac{d}{dx} (x-M)\ii,
$$
which yields
$$
\eqalign
{\frac{d}{dh} [f(\la\,e^{2\pi ih},z)] 
     &= \sum_{D\ge1} \sum_{0\le N\le D-1} 
          \frac{1}{2\pi i} \frac{z^D}{D} 
          \sum_{M\in\Z} \frac{d}{dh} (h+\al-\frac{N+MD}{D})\ii \cr
     &= \frac{1}{2\pi i} \sum_{D\ge1} \sum_{N\in\Z} \frac{z^D}{D} 
          \frac{d}{dh} (h+\al-\frac{N}{D})\ii 
      = \frac{1}{2\pi i} \sum_{D\ge1} \sum_{N\in\Z} \frac{z^D}{D} 
          \frac{d}{dh}[ Z\ii\cdot\frac{h}{h-Z} ],
}
$$
hence the result by integration.
\qed


\beginsection{3.5 Proof of Theorem~3.5}

-- Using the notations of Lemma~3.3 and the change of variable $z=e^s$, we
introduce the functions~$\chi^+$ and~$\chi^-$:
$$
f_\de(\la\,e^{2\pi ih},z) = f_\de(\la,z) + \frac{1}{2\pi i}\bigl(
\chi^+(h,s) + \chi^-(h,s)
\bigr),
\quad
\chi^\pm(h,s) = \sum_{(D,N)\in\cE^\pm} Z\ii \cdot \frac{h}{h-Z} \cdot \frac{e^{Ds}}{D},
$$
still with $Z = \frac{N}{D} - \al$.

Each term
$$
\chi_{(D,N)}(h,s) = Z\ii \cdot \frac{h}{h-Z} \cdot \frac{e^{Ds}}{D},
$$
being analytic at the origin with respect to~$h^{1/2}$, may be written as the Laplace
integral in any direction of its Borel transform; we find it convenient to let a 
factor~$h^{1/2}$ outside the integral:
$$
\chi_{(D,N)}(h,s) = h^{1/2} \int_0^\infty \hat \psi_{(D,N)}(\ze,s) \,
e^{-\ze h^{-1/2}} \,d\ze.
$$
One computes easily
$$
\hat \psi_{(D,N)}(\ze,s) = -Z^{-2}\cdot \frac{e^{Ds}}{D} \sum_{n\ge0} Z^{-n}
\cdot \frac{\ze^{2n}}{(2n)!}
$$
which is entire and of exponential type in any direction:
according to the sign of~$Z$ we obtain a hyperbolic or a trigonometric cosine.
Part~(a) of Theorem~3.5 will thus derive from the study of the 
convergence of the series
$$
\hat\psi^+(\ze,s) = - \sum_{(D,N)\in\cE^+} Z^{-2} \cosh(Z^{-1/2}\ze)
\frac{e^{Ds}}{D}
\Eqno{\defipsip}
$$
and
$$
\hat\psi^-(\ze,s) = - \sum_{(D,N)\in\cE^-} Z^{-2} \cos(|Z|^{-1/2}\ze)
\frac{e^{Ds}}{D}.
\Eqno{\defipsim}
$$

Let us consider~$\hat\psi^+$ for instance. It is of course the even part (with
respect to~$\ze$) of
$$
\Psi^+(\ze,s) = - \sum_{(D,N)\in\cE^+}
D\ii Z^{-2} e^{Z^{-1/2}\ze + Ds}
\Eqno{\defiPsip}
$$
Let $\de>0$ and $\ka_0<\ka_+(\la)$:
we obtain the uniform convergence of this series for
$$
\RE\ze + \ka_0\RE s \le -\de
$$
by observing that $\ka_+(\la)=\liminf_{(D,N)\in\cE^+} \{D Z^{1/2} \}$.
Indeed, for almost all $(D,N)\in\cE^+$ (\ie all of them except a finite number),
$D Z^{1/2}\ge\ka_0$, therefore
$\RE (Z^{-1/2}\ze+Ds) \le - \de \ka_0\ii D$;
and for each $D\ge1$,
$$\eqalign
{\sum_{N\in\Z, (D,N)\in\cE^+} Z^{-2} &= D^2 \sum_{N>\al D} (N-\al D)^{-2} \cr
         &\le D^2 \Bigl( \dist(\al D,\Z)^{-2} + \ze(2) \Bigr) \le \text{const} D^4,
}$$
hence 
$$
|\Psi^+(\ze,s)| \le \text{const} \sum_{D\ge1} D^4 \, e^{- \de \ka_0\ii D}
$$
in that domain.
In particular, $\Psi^+$ is analytic and bounded in
$$
\ao (\ze,s) \,/\ens \RE s<\log r \;\text{and}\; \RE\ze <
\ka_0 (-\RE s) -\de \af
$$
for all $\de>0$ and $\ka_0<\ka_+(\la)$,
and this is enough to establish the analyticity of
$\hat \psi^+(\ze,s)=\frac{1}{2}(\Psi^+(\ze,s)+\Psi^+(-\ze,s))$
in
$$
\ao (\ze,s) \,/\ens \RE s<\log r \;\text{and}\; |\RE\ze| <
\ka_+(\la) (-\RE s) \af.
$$

The same analysis can be performed on~$\hat \psi^-$ which is the even part of
$$
\Psi^-(\ze,s) = - \sum_{(D,N)\in\cE^-}
D\ii Z^{-2} e^{-|Z|^{-1/2}i\ze + Ds},
\Eqno{\defiPsim}
$$
but the factor~$-i$ in front of~$\ze$ in the exponentials is responsible for a
rotation of~$\pi/2$ of the whole picture.

\bigbreak \noindent -- Lemma~3.3 will be useful in the proof of Part~(b) of Theorem~3.5.
The partition of~$\cE^+$ yields a decomposition of~$\Psi^+$:
$$
\Psi^+ = \Psi^+_{\cF^+} + \Psi^+_{\cE^+_*} + \Psi^+_{\cA^+},
\quad \text{with}\ens
\Psi^+_\cB(\ze,s) = - \sum_{(D,N)\in\cB} D\ii Z^{-2} e^{Z^{-1/2}\ze + Ds}.
$$
Because of the properties of~$\cF^+$ and~$\cE^+_*$, the function 
$\Psi^+_{\cF^+} + \Psi^+_{\cE^+_*}$ is analytic in a domain
$$
\ao (\ze,s) \,/\ens \RE s<\log r \;\text{and}\; \RE\ze <
\ka_+' (-\RE s) \af 
$$
which is larger than the domain of analyticity that we just obtained
for~$\Psi^+$, 
as one can see by the same arguments as above.

Let us fix $s\in\C$ with $\RE s<\log r$ and let us consider a
point~$\ze^+_{k,l}(s)$.
When $\ze$ tends to~$\ze^+_{k,l}(s)$, the function 
$\Psi^+_{\cF^+} + \Psi^+_{\cE^+_*}$ tends to its value
at~$(\ze^+_{k,l}(s),s)$,
thus its real part remains finite 
and we now focus on the third term,~$\Psi^+_{\cA^+}$.
According to Lemma~3.3, we can write
$$
\Psi^+_{\cA^+}(\ze,s) = - \sum_{p\ge0} c_p (D_p^+)^3
\,e^{(Z^+_p)^{-1/2}\ze+D^+_p s},
$$
with $Z^+_p = \frac{N^+_p}{D^+_p}-\al$ and $c_p= (Z^+_p)^{-2} (D^+_p)^{-4}$.
Moreover we can study the asymptotic behaviour with respect to~$p$ of these
quantities:
$$
Z_p^+ = (\nu_+(\la) + \rho_p^+) (D_p^+)^{-2}
$$
and $\lim\rho_p^+ = 0$, thus $\lim c_p = \nu_+(\la)^{-2}$.
Let us introduce 
$$
\sig_p = (Z_p^+)^{-1/2} - \ka_+(\la)\ii D_p^+ 
       \sim -\frac{1}{2}\ka_+(\la)^{-3} D_p^+ \rho_p^+.
$$
We know that $\lim \sig_p = 0$, and we can define a function
$$
\Phi^+(X,\ze) = \sum_{p\ge0} c_p (D_p^+)^3 X^{D_p^+} \, e^{\sig_p\ze}
$$
such that
$$
\Psi^+_{\cA^+}(\ze,s) = -\Phi^+(e^{\ka_+(\la)\ii\ze+s},\ze).
$$
According to the definition of~$\ze^+_{k,l}(s)$,
when $\ze$ tends to~$\ze^+_{k,l}(s)$ horizontally by the left, the new variable
$X = e^{\ka_+(\la)\ii\ze+s}$ tends to $e^{2\pi i(k\al+l)} = \la^k$ along the
ray~$]0,\la^k[$. 
Moreover, since $\dist(k\al D_p^+,\Z)$ tends to~0 as $p$ tends to infinity, we
have $\lim \la^{k D_p^+} = 1$.
We are in a position to apply the following elementary result:

\Proc{Lemma~3.5}
{Let $\Phi(X,\ze) = \sum_{p\ge0} c_p d_p^3 \, e^{\sig_p\ze}X^{d_p}$.
Assume that the $\sig_p$ and~$c_p$ are real numbers, with $\lim_{p\to\infty} \sig_p=0$ and
$c_p$ bounded from above and from below by some positive constants,
and that $\{d_p\}$ is an increasing sequence of integers such that $\lim_{p\to\infty}
\la^{k d_p}=1$.
Let $\cK$ be a compact subset of~$\C$.
\item{--} The series which defines~$\Phi$ converges uniformly in~$\cK_0\times\cK$ for
any compact subset~$\cK_0$ of~$\D$.
\item{--} The function~$\RE \Phi(X,\ze)$ tends to~$+\infty$ as $X$ tends to~$\la^k$ along the
ray~$]0,\la^k[$, uniformly with respect 
to $\ze\in\cK$.
}

\Pf{Proof of Lemma~3.5}
The convergence of the series is obvious.
Let $M>0$. The quantity $e^{\sig_p\ze}\la^{k d_p}$ tends to~1 as $p$ tends
to infinity uniformly with respect to~$\ze$, thus we can chose $p_0$ large
enough so that, for all $\ze\in\cK$,
$$
\RE \sum_{p=0}^{p_0} c_p d_p^3 \, e^{\sig_p\ze} \la^{k d_p} \ge 2M
\quad\text{and}\quad
\forall p\ge p_0,\ens
\RE (e^{\sig_p\ze} \la^{k d_p}) \ge \frac{1}{2}.
$$
Let $\de>0$, small enough so that, for all $\ze\in\cK$ and $X\in\C$,
$$
|X-\la^k| \le \de  
\Rightarrow 
|\sum_{p=0}^{p_0} c_p d_p^3 \, e^{\sig_p\ze} (X^{d_p} - \la^{k d_p})| \le M.
$$

We see that, if $\ze\in\cK$ and $X\in\,]0,\la^k[$ with $|X-\la^k|\le\de$,
$$\eqalign
{\RE\Phi(X,\ze) &=
\RE \sum_{p=0}^{p_0} c_p d_p^3 \, e^{\sig_p\ze} (X^{d_p} - \la^{k d_p})
+ \RE \sum_{p=0}^{p_0} c_p d_p^3 \, e^{\sig_p\ze} \la^{k d_p}
+ \RE \sum_{p> p_0} c_p d_p^3 \, e^{\sig_p\ze}X^{d_p} \cr
               &< -M +2M
}$$
(the third term is positive since $X=t \la^k$ with $t\in\,]0,1[$ and
$\RE(e^{\sig_p\ze} \la^{k d_p})>0$).
\qed

\Pf{Continuation of the proof of Theorem~3.5}
We have obtained that $\RE \Psi^+_{\cA^+}$ and thus $\RE \Psi^+$ tend
to~$-\infty$ as $\ze$ tends to~$\ze^+_{k,l}(s)$ horizontally by the left. 
This allows to reach the desired conclusion for~$\hat\psi^+$.
The previous work is easily adapted to the case of~$\hat\psi^-$, with the
introduction of
$$
\Psi^-_{\cA^-}(\ze,s) = - \sum_{p\ge0} c_p (D_p^-)^3
\,e^{-(Z^-_p)^{-1/2}i\ze+D^-_p s},
$$
(with real numbers $c_p$ and~$\sig_p$ associated to~$\cA^-$) and
$$
\Phi^-(X,\xi) = \sum_{p\ge0} c_p (D_p^-)^3 X^{D_p^-} \, e^{\sig_p\xi},
$$
but this time the correspondence is 
$\Psi^-_{\cA^-}(\ze,s) = \Phi^-(e^{-\ka_-(\la)\ii i\ze + s},-i\ze)$.
This ends the proof of Part~(b) of Theorem~3.5.

\bigbreak\noindent -- We now come to Part~(c). 
Let us fix $s<\log r$ and $z=e^s$ (thus $z\in\,]0,1[$).
We recall the notation
$$
\hat\psi = \hat\psi^+ + \hat\psi^- = \sum_{n\ge0} \chi_{n+1}(s) \frac{\ze^{2n}}{(2n)!}.
$$
Our aim is to bound from below half of the coefficients of that series.
According to the formulas~\defipsip--\defipsim, 
$$
\forall n\ge0, \quad
-\chi_{n+1}(s) = \sum_{(D,N)\in\cE^+} z^D D\ii Z^{-n-2}
         + (-1)^n\sum_{(D,N)\in\cE^-} z^D D\ii |Z|^{-n-2},
$$
with the usual notation $Z=\frac{N}{D}-\al$.
Let us choose $\eps\in\{+,-\}$ so that $\ka_\eps(\la) \le \ka_{-\eps}(\la)$.
When we restrict ourselves to even~$n$, only positive quantities appear in the
right-hand side of the above equation, thus we obtain a lower bound for the
left-hand side by retaining only the terms which correspond to~$(D,N)\in\cE^\eps$:
$$
\forall j\ge1, \quad
-\chi_{2j-1}(s) > \sum_{p\ge0} z^{D_p^\eps} (D_p^\eps)\ii |Z_p^\eps|^{-2j}.
$$

According to Lemma~3.3, $|Z_p^\eps| = (\nu_\eps(\la)+\rho_p^\eps)(D_p^\eps)^{-2}$
and $\rho_p^\eps$ tends to~0 as $p$ tends to infinity, thus we can fix $p_0$
large enough and
$c= \frac{3}{2} \nu_\eps(\la)$
so that 
$$
\forall p\ge p_0,\quad  |Z_p^\eps| \le c (D_p^\eps)^{-2}.
$$
For $j\ge D_{p_0}^\eps/4$, we define
$$
E_j = \max_{D_p^\eps \le 4j} \{D_p^\eps\}.
$$
Thus $E_j\le 4j$, and since $E_j\in\ao D_p^\eps, \, p\ge p_0\af$,
we can choose to retain only the corresponding contribution in the previous sum:
$$
-\chi_{2j-1}(s) > z^{E_j} c^{-2j} E_j^{4j-1} > (c\ii z^2)^{2j} E_j^{4j-2},
$$
and for $j$ large enough,
$$
|\chi_{2j-1}(s)|^{\frac{1}{2j-1}} > \de_j := \frac{1}{2} c\ii z^2 E_j^2.
$$

The sequence $\{\de_j\}$ that we just defined is obviously non-decreasing and
there remains only to check that $\sum \de_j^{-3/4} < +\infty$, 
\ie that $\sum E_j^{-3/2} < +\infty$.
We observe that $E_j=F_{4j}$ with
$$
\forall m\ge D_{p_0}^\eps,\quad
F_m = \max_{D_p^\eps \le m} \{D_p^\eps\},
$$
\ie
$$
\eqalign
{F_{D_{p_0}^\eps} &= F_{D_{p_0}^\eps+1} = \ldots = F_{D_{p_0+1}^\eps-1} = D_{p_0}^\eps, \cr
F_{D_{p_0+1}^\eps} &= F_{D_{p_0+1}^\eps+1} = \ldots = F_{D_{p_0+2}^\eps-1} = D_{p_0+1}^\eps,
}
$$
and so on.
Hence, for $P>p_0$,
$$
\sum_{m=D_{p_0}^\eps}^{D_{P+1}^\eps-1} F_m^{-3/2} = 
\sum_{p=p_0}^P (D_p^\eps)^{-3/2} (D_{p+1}^\eps-D_p^\eps)
\le (-1+\sup\{\frac{D_{p+1}^\eps}{D_p^\eps}\}) \sum_{p=p_0}^P (D_p^\eps)^{-1/2},
$$
and the series $\sum F_m^{-3/2}$ and $\sum E_j^{-3/2}$ converge.
\qed


\vfill \eject

\beginsection{4. Resummation at resonances and constant-type points}

For a class of monogenic functions (to which the solutions of the cohomological equation belong),
we have obtained asymptotic expansions at Diophantine points of the unit circle.
Now, restricting ourselves to the subspace of \BWD\ series with poles at resonances,
we will study asymptotic behaviour at resonances.

Then, we will address the question: Is it possible to recover any solution in a constructive way
from its asymptotic expansion at a particular point of~$\S^1$?
We will provide refined results on Gevrey-1 asymptotics at resonances
and Gevrey-2 asymptotics at constant-type points
which show that the answer is positive for each of these points. 
In the latter case there is no contradiction with the non-quasianalyticity of~$\cG_2(\la,B)$
nor with Part~(c) of Theorem~3.4,
since the question amounts to working in a smaller quasianalytic subspace without
demanding it to be a Carleman class.

At resonances a rigid structure appears, which is an elementary case of {\sl
resurgence} [E1] in the case of the fundamental solution~$f_\de$. 
The Borel transform of a given solution
$f=f_\de\odot g$ at a resonance~$\La_0$ can be completely described, the appropriate
Laplace transform then yields the function inside or outside the unit disk, and one can even
recover all the other residues
$\La\cL_{m(\La)} \odot g$ from the singularities of the Borel transform at~$\La_0$
by computing the {\sl Stokes phenomenon}.
In some sense,
this means passing from local information (one particular singular point~$\La_0$) 
to global information (the whole set of ``poles'').

For constant-type points, although it is likely that no quasianalytic Carleman class
contains the solutions (as is the case for quadratic irrationals), one can still
define a quasianalytic space which contains them and in which an adaptation of Borel-Laplace
summation process provides constructive quasianalytic continuation, like for
resonances.


\beginsection{4.1 Asymptotic expansions at resonances}

{\bf a)}
Recall the formulas~\equadefSRa\ and~\equadefSrB\ which, by Theorem~2.2, define
$\fS_\cR:\; \cS(r,B) \to \cM((K_j),B)$.

\Proc{Theorem~4.1}
{Let $r\in\,]0,1[$, $B$ a Banach space and $\La_0\in\cR$.
If $a\in\cS(r,B)$,
the function~$q \mapsto (q-\La_0)\bigl(\fS_\cR(a)\bigr)(q)$ belongs to~$\cG_1(\La_0,B)$
and the constant term in its asymptotic expansion $J_{\La_0}\bigl((q-\La_0)\fS_\cR(a)\bigr)$
is equal to~$a_{\La_0}$.
In particular, if $0<r_1<r_2$, the solution~$F_{r_1,r_2}$ belongs to
$(q-\La_0)\ii \cG_1(\La_0,\cL(B_{r_1},B_{r_2}))$
and the constant term in~$J_{\La_0}\bigl((q-\La_0)F_{r_1,r_2}\bigr)$ 
is~$\La_0 \cL_{m(\La_0)} \odot$.
}

Therefore the \BWD\ series of~$\fS_\cR\bigl(\cS(r,B)\bigr)$ or the solutions of
the cohomological equation are contained in quasianalytic spaces
$(q-\La_0)\ii \cG_1(\La_0,B)$.
Moreover Nevanlinna's Theorem ensures the possibility of following the quasianalytic
continuation of any such \BWD\ series~$f$ across~$\S^1$ ``through~$\La_0$'':
the Borel transform of~$J_{\La_0} \bigl((q-\La_0)f\bigr)$ necessarily belongs
to~$\hat\cN(B)$, and the appropriate Laplace transform restores the function on one
side or the other of~$\S^1$. But much more can be said about the Borel transform 
in the case of the solutions, as will be shown in Section~4.2.

Unfortunately nothing indicates that such a quasianalytic property could be shared by all
the monogenic functions of~$\cM((K_j),B)$ or~$\CH((K^*_{A,j}),B)$.

\Pf{Proof of Theorem~4.1}
Let $a\in\cS(r,B)$ and 
$$
F(q) = \sum_{\La\in\cR,\,\La\neq\La_0} \frac{a_\La}{q-\La} 
     = \bigl(\fS_\cR(a)\bigr)(q) - \frac{a_{\La_0}}{q-\La_0}.
$$
It is sufficient to prove that $F\in\cG_1(\La_0,B)$.

We have $\La_0=e^{2\pi i \al}$ with,
$\forall n/m \in\Q\setminus\{\al\}$, $|\al-n/m|\ge 1/(m(\La_0)|m|)$,
and one checks easily the existence of a positive constant $\ga_1$ such that
$$
\forall \La\in\cR\setminus\{\La_0\},\quad |\La_0-\La| \ge \frac{\ga_1}{m(\La)}.
\Eqno{\ineqdiophS}
$$
Therefore the series 
$$
A_n = (-1)^n\sum_{\La\in\cR,\,\La\neq\La_0} \frac{a_\La}{(\La_0-\La)^{n+1}},
\qquad n\in\N,
$$
are absolutely convergent in~$B$. 
In fact there exists $c>0$ such that, $\forall n\ge0$, 
$\nor A_n \nor \le c \, \ph(n+1)$,
where the function~$\ph$ is defined by 
$$
\forall n\ge0,\quad \ph(n) = \sum_{\La\in\cR,\,\La\neq\La_0} 
\frac{r^{m(\La)}}{m(\La)}\,|\La_0-\La|^{-n}.
$$

\Proc{Lemma~4.1}
{Let $\cK$ be a compact subset of~$\C$ which intersects $\S^1$
at~$\La_0$ only, with finite order of contact~$\be>0$
(\ie $\exists c>0$ such that $\forall q\in\cK$, $\forall q'\in\S^1$, 
$|q-q'|\ge c|\La_0-q'|^\be$).
There exists $c_0>0$ such that
$$
\forall N\ge0,\; \forall q\in\cK,\quad
\nor F(q) - \sum_{0\le n\le N-1} A_n (q-\La_0)^n \nor \le c_0 \,|q-\La_0|^N\, \ph(N+\be).
$$
}

\Pf{Proof of Lemma~4.1}
One computes easily the identity
$$
F(q) - \sum_{0\le n\le N-1} A_n (q-\La_0)^n = (-1)^N (q-\La_0)^N
\sum_{\La\in\cR,\,\La\neq\La_0} a_\La \, \frac{(\La_0-\La)^{-N}}{q-\La}.
$$
But for $q\in\cK$ and $\La\in\cR$,
$|q-\La|\ge c|\La_0-\La|^\be$,
whereas $\nor a_\La \nor \le \text{const} \frac{r^{m(\La)}}{m(\La)}$.
\qed

\Pf{End of the proof of Theorem~4.1}
Let us check the existence of~$c_1>0$ such that
$$
\forall n\ge0,\quad \ph(n) \le c_1^{n+1} n!.
$$
Using the inequality~\ineqdiophS\ we obtain
$$
\forall n\ge0,\quad
\ph(n) \le \ga_1^{-n} \sum_{m\ge1} m^n\, r^m.
$$
If we set $r=e^{-s}$ with $s>0$ and compare the sum $\sum_{m\ge1} m^n\, e^{-ms}$ and the integral 
$\int_0^{+\infty} m^n\,e^{-m s}\,dm = s^{-n-1}n!$,
we obtain
$\sum m^n\, e^{-ms} \le s^{-n-1}(\Ga(n+1) + 2 s \,n^n\,e^{-n})$; the Stirling formula
yields the desired inequality.

We now choose for $\cK$ a closed disk $\bar\De^\pm$ contained in~$\overline\D$ or~$\overline\E$,
then
$\be=2$ and ${F}_{|\bar\De^\pm}\in\cG^\pm_1(\La_0,B)$ with 
$J^\pm_{\La_0}({F}_{|\bar\De^\pm}) = \sum A_n Q^n$.
\qed

Notice that, according to the proof of Theorem~3.5, if $a\in\cS(r,B)$
the \BWD\ series~$\fS_\cR(a)$ admits a
Gevrey-$1$ asymptotic expansion at~$\La_0$ in compact subsets~$\cK$ with arbitrarily high
order of contact at~$\La_0$.



\vfill\eject

\beginsection{4.2 Resurgence of the fundamental solution at resonances}

We fix in this section a resonant
point~$\La_0\in\cR$ with $m_0=m(\La_0)$.
We denote by $n_0$ the integer such that 
$$
\La_0 = e^{2\pi i n_0/m_0},\quad 0\le n_0\le m_0-1,\quad (n_0|m_0)=1.
$$
We know by Theorem~4.1 that the function~$(q-\La_0) f_\de$ belongs to~$\cG_1(\La_0,B_r)$ for
all $r\in\,]0,1[$,  with an asymptotic expansion
$$
J_{\La_0}\bigl((q-\La_0)f_\de\bigr) = \sum_{n\ge0} a_n Q^n, \quad
a_0 = \La_0\cL_{m_0}, \quad
(\forall n\ge0)\; a_n\in B_r.
$$
According to Theorem~3.2 and Definitions~3.3 and~3.4, the Borel transform
$$
\hat\Phi^\de(\xi) = \sum_{n\ge0} \La_0^{n+1} a_{n+1} \frac{\xi^n}{n!}
$$
belongs to~$\hat\cN(B_r)$ for all $r\in\,]0,1[$, and $f_\de$ can be recovered
from~$\hat\Phi^\de$ by the formula
$(q-\La_0)f_\de = \La_0\cL_{m_0} + \L_{\La_0} \hat\Phi^\de$,
which can be rephrased as
$$
f_\de\bigl(\La_0(1+t)\bigr) = t\ii \cL_{m_0} + \La_0\ii t\ii \Phi^\de(t),
\qquad
\Phi^\de(t) = \int_0^{\pm\infty} \hat\Phi^\de(\xi)\, e^{-\xi/t} \, d\xi.
$$ 
We may consider $\hat\Phi^\de$ as a holomorphic function of two variables as
well, by setting 
$\hat\Phi^\de(\xi,z)=\hat\Phi^\de(\xi)(z)$. Our goal is now to study the analytic
continuation with respect to~$\xi$ of this Borel transform.

\Def{Definition~4.1}
{For $a\in \Z^*$ and $b\in \Z$, we define the {\sl moving singular point}
$$
z\in\D^* \ens\mapsto\ens
\xi_{a,b}(z) = \frac{2\pi a}{m_0}(-i\log z + \frac{2\pi b}{m_0}) \in\C,
$$
where $\D^*=\D\setminus\{0\}$ and we have chosen some determination of the logarithm once
for all. We also attach to it a complex number:
$$
C_{a,b} = -\frac{1}{m_0} \, e^{2\pi i a b n'_0/m_0},
$$
where $n'_0+m_0\Z$ is the multiplicative inverse of~$n_0+m_0\Z$ in the ring~$\Z/m_0\Z$.
}

\vskip .3cm
\epsfysize=6cm
\centerline{\epsfbox{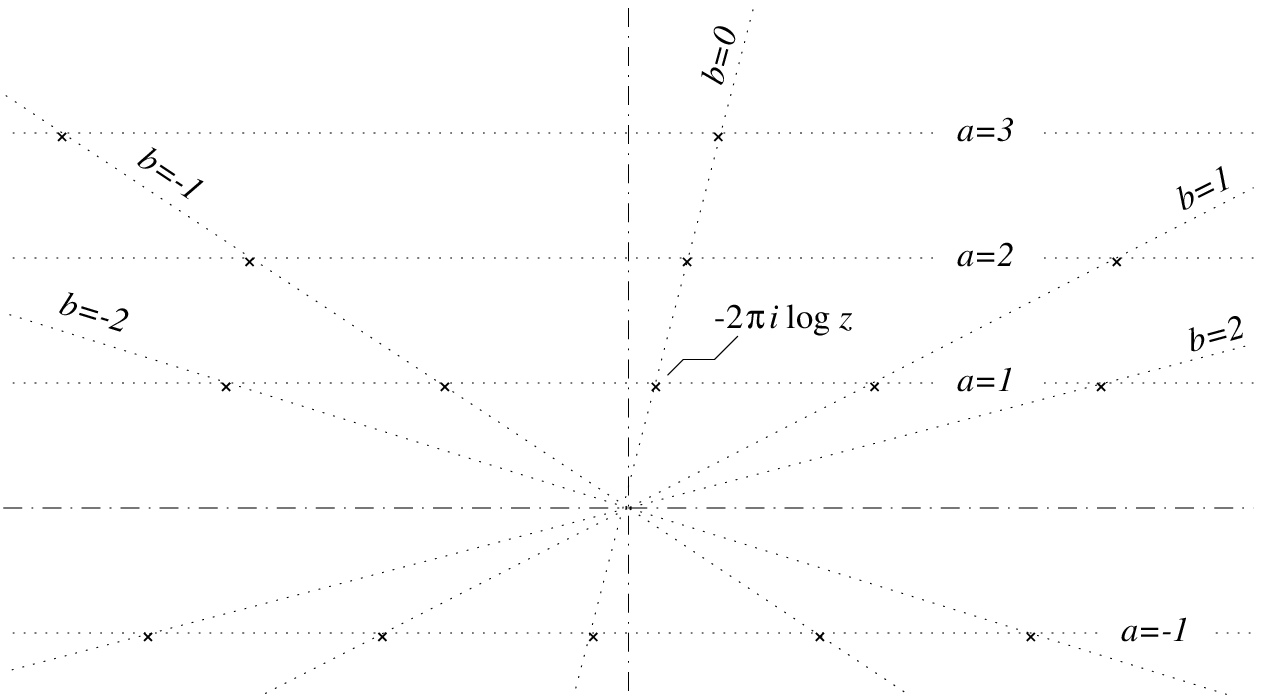}}
\vskip .3cm

\centerline{\eightrm The points 
$\sst \xi_{a,b}(z)$
lie at the intersection}
\centerline{\eightrm of two family of lines parametrized by $\sst a\in\Z^*$ or
$\sst b\in\Z$.}
\vskip .3cm

\Proc{Theorem~4.2 (Resurgence at resonances)}
{For each $z\in\D^*$, the function $\xi\mapsto\hat\Phi^\de(\xi,z)$
extends analytically to the universal covering\footnote{\noteRiemsurf}
{\rm This simply means that for $\hat\Phi^\de(\xi,z)$ viewed as an analytic germ in~$\xi$ at the
origin, analytic continuation can be followed along any path issuing from the origin and lying
in~$\C\setminus\{\xi_{a,b}(z)\}$. We obtain a Riemann surface by considering homotopy classes
of such pathes; its {\sl main sheet} corresponds to rectilinear paths and can be identified to
the holomorphic star of our germ. }  
of $\C\setminus\{\xi_{a,b}(z),\;  a\in \Z^*,\; b\in \Z\}$;
near a moving singular point~$\om=\xi_{a,b}(z)$ on the main sheet of this Riemann
surface, one can write
$$
\hat\Phi^\de(\xi,z) = 
\La_0 C_{a,b} \bigl( {e^{-\om/2}\over\xi-\om} +
\hat L_\om(\xi-\om)\log(\xi-\om) \bigr) + \hbox{regular function,}
$$
where $\hat L_\om$ is an entire function.
Moreover, for any $z\in\D^*$ and for any line~$\De$ of~$\C$ passing through the origin and
avoiding the singular points~$\xi_{a,b}(z)$, the function~$\hat\Phi^\de(\xi,z)$ has at
most exponential growth for~$\xi\in\De$.
}

It is even possible to compute the entire functions~$\hat L_\om$:
they are the Borel transforms of the convergent series
$L_\om(t) = -e^{-\om/2} + \bigl( 1+tL(t)\bigr) e^{-\om L(t)}=O(t)$,
where 
$$
L(t)= \bigl( \log(1+t) \bigr)\ii - t\ii = \demi + O(t).
$$
This theorem will appear as a consequence of Theorem~4.3 below.

In the terminology of resurgence, $(q-\La_0)f_\de(q)$ would be called a {\sl simple resurgent
function} (see Appendix~A.5).
Theorem~4.2 shows that the index~1 in the Gevrey asymptotics provided by
Theorem~4.1 is optimal, since the Borel transform~$\hat\Phi^\de$ has finite radius of
convergence with respect to~$\xi$ for each nonzero~$z$.

\remark{4.1}
{There is some analogy between the first line of moving singular 
points~$\xi_{1,b}(z)$ and
the points~$\ze_{k,l}^\pm(s)$ of Theorem~3.5~(b).
Both cases deal with the Borel transform of some Gevrey-$\tau$ asymptotic
expansion at a point of~$\oDC_\tau$, at $\la=e^{2\pi i\al}\in\IQ$ in Section~3.3
($\tau=2$) and at $\La_0=e^{2\pi i n_0/m_0}\in\cR$ here 
($\tau=1$; indeed \ineqdiophS\ leads us to set $\oDC_1=\cR$). 
We have $s=\log z$,
but in Section~3.3 we were expanding with respect to~$h$ defined by
$q=\la\,e^{2\pi i h}$ (and then computing a Borel transform with respect
to~$h^{1/\tau}$) instead of $t=\frac{q-\la}{\la}$,
and this is responsible for a scaling by a factor~$2\pi i$ between the
variables~$\ze$ and~$\xi$ for the Borel transforms.
The special singular points~$\ze_{k,l}^+(s)$ can be defined by
$$
s + \ka\ii \ze_{k,l}^+(s) = 2\pi i(k\al+l), \qquad k,l\in\Z,
$$
where $\ka=\ka_+(\la)$ is the largest number such that $|\frac{N}{D}-\al|\ge
(\frac{\ka}{D})^\tau$ for all $\frac{N}{D}>\al$ except a finite number of them
(recall that $\tau=2$ in that case).
In the resonant case we can set $\ka=\frac{1}{m_0}$: this is the largest number such that
$|\frac{N}{D}-\frac{n_0}{m_0}|\ge \frac{\ka}{D}$ for all $\frac{N}{D}\neq\al$
($\tau=1$ in this case and we need not distinguish left and right rational
approximations of~$n_0/m_0$). The first line of moving singular points appears
to be defined by
$$
s + \frac{1}{2\pi i}\ka\ii \xi_{1,b}(z) =
-2\pi i \frac{b}{m_0},
\qquad b\in\Z,
$$
but the group $\ao -\frac{b}{m_0};\; b\in\Z \af = 
\ao k\frac{n_0}{m_0}+l;\; k,l\in\Z \af$ is discrete, thus the singular points
are isolated (hence the resurgence property),
whereas $\ao k\al+l;\; k,l\in\Z \af$ was dense in~$\R$, hence the natural
boundary for~$\hat\psi^+(\ze,s)$.
}

If a function $g\in B_{r_1}$ is given, for some $r_1>0$, one can deduce results for the
corresponding solution $f=f_\de\odot g$:
we know by Theorem~4.1 that $(q-\La_0) f\in \cG_1(\La_0,B_{r_2})$ 
and the function
$\hat\Phi^g = \cB \circ J_{\La_0}\bigl(-\La_0\cL_{m_0}\odot g+(q-\La_0)f\bigr)$
belongs to~$\hat\cN(B_{r_2})$ for all $r_2\in\,]0,r_1[$.
In fact, for each $\xi\in\C$, $\hat\Phi^g = \hat\Phi^\de\odot g$
and the singularities with respect to~$\xi$ of~$\hat\Phi^g$ depend on the
singularities (with respect to~$z$) of~$g$.
More precisely, the location of the moving singular points of~$\hat\Phi^\de$ shows that
$\hat\Phi^\de$ is holomorphic in $\ao (\xi,z)\,|\; z\in\D,\;
|\IM\xi|<2\pi\ln\frac{1}{|z|} \af = \ao |z|<\exp(-\frac{|\IM\xi|}{2\pi}) \af$;
thus $\hat\Phi^g$ is holomorphic in $\ao |z|<r_1\exp(-\frac{|\IM\xi|}{2\pi})
\af$, which means that for each $z\in\D_{r_1}$, $\hat\Phi^g$ is holomorphic with 
respect to~$\xi$ in a horizontal strip of width 
$4\pi r_1\ln\frac{1}{|z|}$.
(But $\hat\Phi^g$ may have a natural boundary with respect to~$\xi$ if this is the 
case for $g$ with respect to~$z$.)

\medbreak
So far we were dealing with Borel transforms with respect to~$t=\frac{q-\La_0}{\La_0}$, 
but in fact the variable $\eta = \log(q/\La_0)$ is more convenient.
Thus we consider the function 
$$
\eta \mapsto \Psi^g(\eta) = \eta f(\La_0\,e^\eta)
$$ 
still for a general solution $f=f_\de\odot g$: it admits a Gevrey-1 asymptotic expansion
$$
\ti\Psi^g(\eta) = \sum_{p\ge0} \Psi^g_p \eta^p
$$
for $\eta$ tending to zero by the left or by the right, whose constant term
is~$\Psi^g_0=\cL_{m_0}\odot g$, and we are interested in the Borel transforms
$$
\hat\Psi^g   = \sum_{n\ge0}  \Psi^g_{n+1}  \frac{\xi^n}{n!},
\qquad
\hat\Psi^\de = \sum_{n\ge0} \Psi^\de_{n+1} \frac{\xi^n}{n!}.
$$

\Proc{Theorem~4.3 (Borel transform at resonances)}
{When viewed as a holomorphic function of two variables, $\hat\Psi^g$ can be written
$$
\hat\Psi^g(\xi,z) = \sum_{k=0}^{m_0-1} (\frac{k}{m_0} - \demi) g(\La_0^k z)
                  - \sum_{k=0}^{m_0-1} \sum_{a\in\Z^*}
                                   \frac{e^{2\pi i\frac{ka}{m_0}}}{2\pi i a} \bigl[
                                   g(\La_0^k z\, e^{\frac{m_0\xi}{2\pi i a}}) - g(\La_0^k z)
                                   \bigr]
$$
for $|\IM\xi|$ and $|z|$ small enough.
In particular, for each $z\in\D^*$, the function $\xi\mapsto\hat\Psi^\de(\xi,z)$ is
meromorphic with simple poles only, located at the points~$\xi_{a,b}(z)$, with
$C_{a,b}$ as corresponding residues.
Moreover, for any $z\in\D^*$ and for any line~$\De$ of~$\C$ passing through the origin and
avoiding the poles~$\xi_{a,b}(z)$, the function $(1+|\xi|)\ii\hat\Psi^\de(\xi,z)$ is
bounded for~$\xi\in\De$. 
}

\longremark{4.2}
{The knowledge of the residues of~$\hat\Psi^\de$ with respect to~$\xi$ allows to
compute the ``residues'' of~$f_\de$ with respect to~$q$, \ie to determine the
sequence~$(a_\La)$ such that $f_\de=\fS_\cR((a_\La))$. In other words the complete
asymptotic expansion of~$f_\de$ at one resonance contains the information on
the leading term in the asymptotics at all other resonances.

Indeed let us fix $\La=e^{2\pi i n/m}\in\cR$, with $\frac{n}{m}>\frac{n_0}{m_0}$ 
for conveniency (and as always $m\in\N^*$, $m\in\Z$, $(n|m)=1$), and $z\in\D$,
$s=\log z$ (the dependence on~$z$ of the various functions below will be usually omitted).
We will check directly from Theorem~4.3 that
$f_\de(q,z) \sim \frac{a_\La}{q-\La} = \frac{\La}{q-\La}\cL_m(z)$ for $q$ tending non-tangentially 
w.r.t.~$\S^1$
to~$\La$, which is obviously equivalent to
$$
f_\de(e^{2\pi i h},z) \sim \frac{\cL_m(z)}{2\pi i(h-\frac{n}{m})}
$$
for $h$ tending non-tangentially w.r.t.~$\R$ to~$\frac{n}{m}$.

Let us choose a direction~$\th$ in~$]0,\pi[$ such that
$\arg(\xi_{1,0}(z)) < \th < \arg(\xi_{1,-1}(z))$.
By Cauchy theorem, we can compare the two Laplace transforms
$$
\Psi^\de(\eta) = \eta f_\de(\La_0 e^\eta) = \cL_{m_0} + 
\int_0^{+\infty} \hat\Psi^\de(\xi) \, e^{-\xi/\eta} \, d\xi
\quad \text{for} \ens \RE\eta>0
$$
and
$$
\Psi_\th^\de(\eta) =  \int_0^{e^{i\th}\infty} \hat\Psi^\de(\xi)
\, e^{-\xi/\eta} \, d\xi
\quad \text{for} \ens \RE(\eta\,e^{-i\th})>0.
$$

\vskip .3cm
\epsfysize=3.5cm
\centerline{\epsfbox{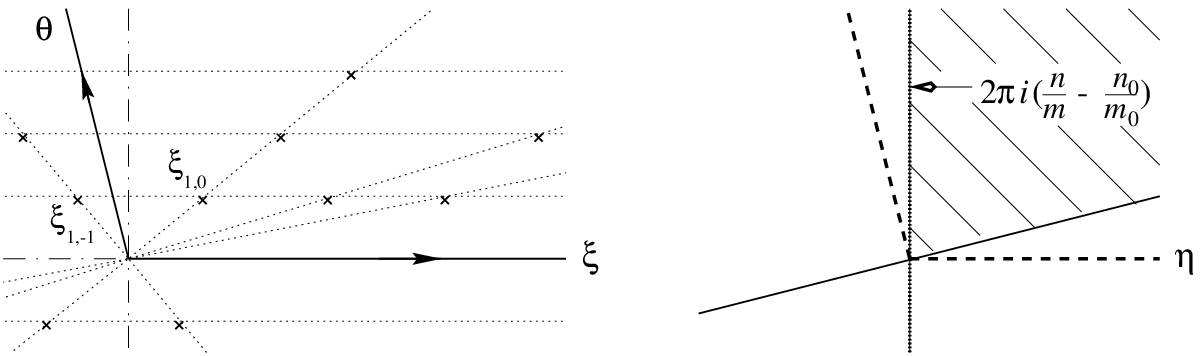}}
\vskip .3cm

If $\RE\eta>0$ and $\RE(\eta\,e^{-i\th})>0$,
\ie $\eta$ belongs to the intersection of the two half-planes, 
$$
\eta f_\de(\La_0 e^\eta) = \cL_{m_0} + \Psi_\th^\de(\eta) 
   + 2\pi i\sum_{a\ge1,b\ge0} C_{a,b} \, e^{-\xi_{a,b}/\eta}.
$$
We are interested in $\eta$ tending to~$2\pi i(\frac{n}{m}-\frac{n_0}{m_0})$
from the right. The term $\Psi_\th^\de(\eta)$ is regular there and will yield no 
contribution in the singular behaviour that we want to analyze.
On the contrary,
for each $a\ge1$, the sum of the geometric series
$$
2\pi i\sum_{b\ge0} C_{a,b} \, e^{-\xi_{a,b}/\eta} =
-\frac{2\pi i}{m_0} \cdot  \frac{e^{\frac{2\pi i a s}{m_0\eta}}}{
		1 - e^{\frac{2\pi i a}{m_0}(
			n'_0 + \frac{2\pi i}{m_0\eta}
			)}
		}
$$
defines a function which is meromorphic w.r.t.~$1/\eta$.
Translating this in the variable~$h=\frac{n_0}{m_0}+\frac{\eta}{2\pi i}$
($h$ tends to~$n/m$ with $\IM h<0$),
we obtain
$$
(h-\frac{n_0}{m_0}) f_\de(e^{2\pi i h}) = -\frac{1}{m_0} \,
\sum_{a\ge1} \frac{e^{\frac{a s}{m_0 h - n_0}}}{
		1 - e^{2\pi i a \frac{n'_0 h+m'_0}{m_0 h- n_0}}
		}
+ \text{regular function,}
$$
where we have introduced $m'_0\in\Z$ defined by $m_0 m'_0 + n_0 n'_0 =1$.

The image of~$\frac{n}{m}$ by the linear fractional map $h\mapsto \frac{n'_0 h+m'_0}{m_0
h-n_0}$ is~$\frac{N}{M}$, where
$N = n'_0 n + m'_0 m$,
$M = m_0 n  - n_0 m$
and $(N|M)=1$.
The only terms contributing to the singularity at $h=\frac{n}{m}$ correspond thus to
$a=jM$, $j\ge1$, and an easy computation allows to conclude that 
$$
(h-\frac{n_0}{m_0}) f_\de(e^{2\pi i h}) \sim 
\frac{1}{2\pi i} \cdot \frac{M}{m_0 m} \cdot \frac{1}{h-\frac{n}{m}} \cdot
\sum_{j\ge1} \frac{e^{jms}}{jm},
$$
hence
$\lim (h-\frac{n}{m}) f_\de(e^{2\pi i h}) = \frac{1}{2\pi i} \cL_m(z)$.
}


\vfill\eject

\beginsection{4.3 Proof of Theorems 4.2 and 4.3}

{\sl Theorem~4.3 implies Theorem~4.2:}
This is an exercise of application of the general theory of which Appendix~A.5
gives a brief account. We will relate $\hat\Phi^\de(\xi)$ and~$\hat\Psi^\de(\xi)$
(from now on we will omit the dependence on the variable~$z$), and first prove that
$\hat\Phi^\de$ extends analytically to the universal covering~$\cC$ of
$\C\setminus\{\xi_{a,b}\}$ with at most exponential growth at infinity
just because $\hat\Psi^\de(\xi)$ has that property.

In the vicinity of the resonant point~$\La_0$, we have two local variables $t$
and~$\eta$:
$$
q = \La_0(1+t) = \La_0 \,e^\eta,
$$
and correspondingly two representations of~$f_\de$ as a Laplace transform:
$$
t f_\de = \cL_{m_0} + \La_0\ii \L_{(\xi\rightarrow t)} \hat\Phi^\de,
\quad
\eta f_\de = \cL_{m_0} + \L_{(\xi\rightarrow \eta)} \hat\Psi^\de. 
$$
We retain that, under the change of variable $t = e^\eta -1 \Leftrightarrow
\eta = \log(1+t)$, 
$$
\L_{(\xi\rightarrow t)} \hat\Phi^\de = \La_0 tL(t) \cL_{m_0} + \La_0 (1+tL(t))
\L_{(\xi\rightarrow \eta)} \hat\Psi^\de,
\quad
L(t)= \bigl( \log(1+t) \bigr)\ii - t\ii = \demi + O(t).
\Eqno{\PhiPsi}
$$

Now we can write $\L_{(\xi\rightarrow \eta)} \hat\Psi^\de =
\L_{(\xi\rightarrow t)} \hat\chi$, \ie we can interpret the change of variable
in the Borel plane, by defining~$\hat\chi(\xi)$ as follows:
$$
\hat\chi(\xi) = e^{-\xi/2}\hat\Psi^\de(\xi) + \sum_{r\ge1} \hat \ell^{*r} *
{\hat\pa^r(e^{-\xi/2}\hat\Psi^\de)\over r!},
$$
where $\hat\ell(\xi)$ is the Borel transform of $\ell(t) = -\frac{1}{2} + L(t)$
and is thus an entire function of exponential type.
(This is because $\eta\ii = t\ii + \frac{1}{2} + \ell(t)$: the translation
by~$1/2$ is responsible for the multiplication by~$e^{-\xi/2}$, and we are then
left with {\em composition-convolution} as described in Appendix~A.5. 
The
notation~$\hat\pa$ simply means multiplication by~$-\xi$, the Borel counterpart
of differentiation w.r.t.~$t\ii$.)

We observe that $\hat\chi$ extends analytically to~$\cC$ with at most
exponential growth at infinity, thus this is also the case for
$$
\hat\Phi^\de = \La_0 \bigl( \hat M \cL_{m_0} + \hat\chi + \hat M * \hat\chi \bigr),
$$
where the entire function~$\hat M$ is simply the Borel transform of~$tL(t)$
(thus $\hat M = \frac{1}{2} + 1*\hat\ell$).

We must now compute the singularity of~$\hat\Phi^\de$ at a point
$\om=\xi_{a,b}$.
For that purpose we can use \'Ecalle's formalism of {\em alien calculus}:
in our particular case, the result to be checked is equivalent to the formula
$$
\De_{(\om\rightarrow t)} \Phi^\de = 2\pi i\La_0 C_{a,b} (e^{-\om/2} + L_\om(t)),
$$
whereas the indications of Theorem~4.3 on the poles of~$\hat\Psi^\de$ amount to
$$
\De_{(\om\rightarrow\eta)} \Psi^\de = 2\pi i C_{a,b}.
$$ 
The operator~$\De_{(\om\rightarrow t)}$ is the {\em alien derivation} of
index~$\om$ relative to the variable~$t$; 
it is defined so to measure the singular behaviour at~$\om$ of the
Borel transform w.r.t.~$t$ of the function on which it is evaluated. For
instance it vanishes on~$tL(t)$ since the corresponding Borel transform is
entire.  
The result to be checked is a consequence of the relation~\PhiPsi\ and
of the fact that $\De_{(\om\rightarrow t)}$ is a derivation and 
$e^{-\om t\ii}\De_{(\om\rightarrow t)} = e^{-\om \eta\ii}\De_{(\om\rightarrow\eta)}$
under the change of variable $\eta\ii = t\ii+L(t)$.  
Indeed, when applied to~\PhiPsi, these rules imply that 
$$ e^{-\om t\ii} \De_{(\om\rightarrow t)} \Phi^\de = \La_0 (1+tL(t)) e^{-\om(t\ii+L(t))}
2\pi i C_{a,b},
$$
while precisely $(1+tL(t)) e^{-\om L(t)} = e^{-\om/2} + L_\om(t)$.
\qed


\Pf{Proof of Theorem~4.3}
Since $\hat\Psi^g = \hat\Psi^\de \odot g$ and
$g(\la z) = \de(\la z) \odot g(z)$ for all $\la\in\C$, it is sufficient to
consider the case where $g=\de$.
>>From now on we will omit the superscript~$\de$.
We also replace the variable~$z$ by $s=\log z$ (and still keep the same names for 
some of our functions), so that 
$$
\Psi(\eta,s) = \eta f_\de(\La_0 e^\eta,e^s) \sim \ti\Psi(\eta,s)
\qquad\text{as $\eta\rightarrow 0$,}
$$
and our goal is to study the Borel transform~$\hat\Psi(\xi,s)$ of that
asymptotic series~$\ti\Psi(\eta,s)$.

>>From the cohomological equation that $f_\de$ satisfies, we deduce an equation
which admits $\Psi$ as solution (and thus $\ti\Psi$ as formal solution):
$$
\Psi(\eta,s+\Om+\eta) - \Psi(\eta,s) = \eta\, \varphi(s)\,,
\Eqno{\eqFE}
$$
where 
$$
\Om = 2\pi i \frac{n_0}{m_0},
\qquad
\varphi(s) = \frac{e^s}{1-e^s}.
$$
In fact, at this level, one can retain this sole equation and forget everything else.

\Proc{Lemma~4.2}
{The equation~\eqFE\  admits a unique formal solution
$$
\ti\Psi(\eta,s) = \sum_{p\ge0} \eta^p \Psi_{p}(s)
$$
with coefficients analytic in $z=e^s$ and vanishing for~$z=0$.
This solution is explicitly given by formulas~\eqSF, \eqphimr\ 
and~\eqGa\  below;
in particular, $\Psi_{0}(s) = \cL_{m_0}(e^s) =
-{1\over {m_0}}\log(1-e^{{m_0}s})$.
}

\proof
Keeping in mind that the solution is required to be $2\pi i$-periodic 
in~$s$, we introduce the following linear combinations of the 
$\Om$-translations of~$\Psi$:
$$\matrix
{\hfill\sig_{r}(\eta,s)\hskip -.7em &=
     & \hskip -.7em\dst\sum_{k=0}^{{m_0}-1}{\La^{-kr}\over {m_0}}\,\Psi^{[k]}(\eta,s)\hfill
     & \hbox{for}\ens r=0,1,\ldots,{m_0}-1, \cr\noalign{\smallskip}
\hfill\Psi^{[k]}(\eta,s)\hskip -.7em &=
     & \hskip -.7em \Psi(\eta,s+k\Om) \hfill
     & \hbox{for}\ens k=0,1,\ldots,{m_0}-1.
}
$$
The identities
$$
\sum_{r=0}^{{m_0}-1} {\La^{-kr}\over {m_0}} = \cases{
1 & if $k=0$ \cr\noalign{\smallskip}
0 & if $k=1,\ldots,{m_0}-1$
}
$$
yield the inverse formulas
$$
\Psi^{[k]}(\eta,s) = \sum_{r=0}^{{m_0}-1} \La^{kr} \sig_{r}(\eta,s) \quad
\hbox{for} \ens k=0,1,\ldots,{m_0}-1.
$$

By combining the $\Om$-translations of equation~\eqFE, we obtain the 
system of equations
$$
\La^r \sig_{r}(\eta,s+\eta) - \sig_{r}(\eta,s) = \eta \,
\varphi_{{m_0},r}(s)
\Eqno{(*)_{r}}
$$
where
$$
\varphi_{{m_0},r}(s) 
      = \sum_{k=0}^{{m_0}-1} {\La^{-kr}\over {m_0}}\,\varphi(s+k\Om)
      = \sum_{\la\in\cR_{m_0}} {\la^{-r}\over {m_0}}\,\varphi(s+\log\la)
 \,,\Eqno{\eqphimr}
$$
for $r=0$,~1, \dots,~${m_0}-1$. 
The left-hand side of equation~$(*)_{r}$ may be viewed as a 
``differential operator of infinite order'' $(\La^r e^{\eta\pa_{s}}-\id)$
acting on~$\sig_{r}$.
Let us introduce some elementary functions which are analytic at the 
origin:
$$
\Ga_{a}(X) = {X\over a e^{X}-1} = \sum_{p\ge0} \ga_{p}(a) X^p \quad
\hbox{for}\ens a\in\C^{*}.
\Eqno{\eqGa}
$$
Note that $\ga_{0}(a) =0$ if $a\neq1$, but $\ga_{0}(1) =1$ and in fact
$$
\Ga_{1}(X) = 1 -{X\over2} -\sum_{l\ge1} (-1)^l B_{l} {X^{2l}\over (2l)!}
$$
where the coefficients~$B_{l}$ are the Bernoulli numbers.

The functions~$\Ga_a$ allow us to solve explicitly the system:
$$
(*)_{r} \quad \Leftrightarrow \quad 
         \sig_r = \pa_s\ii\Ga_{\La^r}(\eta\pa_s)\varphi_{{m_0},r} =
                  \ga_0(\La^r)\pa_s\ii\varphi_{{m_0},r} + \sum_{p\ge1}
                  \eta^p \ga_p(\La^r) \pa_s^{p-1} \varphi_{{m_0},r} 
$$
for $r=0,1,\ldots,{m_0}-1$,
with the notation~$\pa_{s}\ii$ for the unique primitive with respect 
to~$s$ which vanishes when $z=e^s$ vanishes.

Thus, we obtain only one possible formal solution of~\eqFE:
$$
\ti\Psi = \pa_s\ii\sum_{r=0}^{{m_0}-1} \Ga_{\La^r}(\eta\pa_{s}) 
           \varphi_{{m_0},r} 
        = \pa_{s}\ii\varphi_{{m_0},0} + \sum_{p\ge1} \eta^p \pa_{s}^{p-1}
            \sum_{r=0}^{{m_0}-1} \ga_{p}(\La^r) \varphi_{{m_0},r}
\Eqno{\eqSF}
$$
Since
$$
\pa_{s}\ii \varphi (s) = - \log(1-e^s)\,,
$$
we recognize the function~$\cL_{m_0}$ in the constant term:
$$\eqalign
{\Psi_{0}(s) = \pa_{s}\ii \varphi_{{m_0},0}(s) 
           &= -{1\over {m_0}} \log \prod_{k=0}^{{m_0}-1} (1-e^{s+k\Om}) \cr
           &= -{1\over {m_0}} \log (1-e^{{m_0}s}) \,.
}$$

The formal series $\ti\Psi$ that we just defined is indeed a solution 
of equation~\eqFE: for any $k=1$, \dots,~${m_0}-1$, the formal series
$$
\ti\Psi^{[k]}(\eta,s) = \sum_{r=0}^{{m_0}-1} \La^{kr} \sig_{r}(\eta,s)
$$
is actually equal to the translation $\ti\Psi(\eta,s+k\Om)$ of~$\ti\Psi$, since
for each~$r$ the series $\sig_r$ is obtained from $\varphi_{{m_0},r}$ by applying
an operator which commutes with the translations, and
$$
\La^{kr} \varphi_{{m_0},r}(s) = \varphi_{{m_0},r}(s+k\Om) \,.
$$
This remark ends the proof of the lemma.
\qed

\remark{4.3}
{The formula that we obtained is reminiscent of the Euler-MacLaurin
formula, one of the early sources of divergent asymptotic series. We will
analyze it by using the formal Borel transform.\footnote{\noteCDD}
{In [CCD] too Borel transform is used in relation with the Euler-MacLaurin
formula, but not with respect to the same variable; our problem pertains rather
to {\sl parametric resurgence} according to \'Ecalle's terminology.
}}

The above work will now allow us to compute the Borel transform w.r.t.~$\eta$ of
$\ti\Psi-\Psi_0$.
The starting point  is the following 
decomposition of the functions $\Ga_{\La^r}$ which appear
in the formula~\eqSF:
$$
\Ga_{\La^r}(X) = -{X\over 2} + 
\Eisen{\nu\in2i\pi\Z}
{X\over X+r\Om-\nu} \,,
$$
where the symbol~$\sum^e$ denotes {Eisenstein summation} [We]: 
terms corresponding to opposite indices are grouped in order to ensure 
convergence, \ie
$$
\Eisen{l\in\Z} = \lim_{L\rightarrow +\infty} \sum_{l=-L}^{+L} \,.
$$
This decomposition results from the identity
$$
\Ga_{\La^r}(X) = {X\over2}(\coth {X+r\Om\over2} -1)
$$
and from the classical decomposition 
$$
\coth X = \Eisen{l\in\Z} {1\over X-il\pi} \,.
$$

It implies that, for $r=0$,~1, \dots,~${m_0}-1$,
$$
\sum_{p\ge0} \ga_{p+1}(\La^r) X^{p+1} = -{\txt{1\over2}}X - 
    \BigEisen{\nu\in2i\pi\Z}{\nu\neq0{\rm\ if\ } r=0}
        \sum_{p\ge0} (\nu-r\Om)^{-p-1} X^{p+1},
$$
so
$$
\sum_{p\ge0} \ga_{p+1}(\La^r) {(\xi\pa_s)^p\over p!} = -{\txt{1\over2}}\,\id - 
    \BigEisen{\nu\in2i\pi\Z}{\nu\neq0{\rm\ if\ } r=0} 
        (\nu-r\Om)\ii e^{(\nu-r\Om)\ii\xi\pa_s}.
$$
According to the formula~\eqSF\  and because of the Taylor formula, the Borel
transform of~$\ti \Psi - \Psi_0$ can thus be written
$$\eqalign
{\hat \Psi(\xi,s) &= - \sum_{r=0}^{{m_0}-1} \biggl(
                 {\txt{1\over2}}\varphi_{{m_0},r}(s) + 
                 \BigEisen{\nu\in2i\pi\Z}{\nu\neq0{\rm\ if\ } r=0}
                     (\nu-r\Om)\ii \varphi_{{m_0},r}(s+(\nu-r\Om)\ii\xi)
                 \biggr) \cr\noalign{\smallskip}
                 &= -{\txt{1\over2}}\varphi(s) - 
                 \BigEisen{l\in\Z,0\le r\le 
                 {m_0}-1}{(l,r)\neq0}\,\sum_{k=0}^{{m_0}-1}\,
                 {\La^{-kr}\over {m_0}}
                 (2\pi il-r\Om)\ii \varphi(s+k\Om+(2\pi il-r\Om)\ii\xi)\,.
}$$

This is an equality between formal series of powers of~$\xi$, the 
right-hand side being considered as a formal Taylor expansion 
(Eisenstein summation ensures that each of its coefficients is well 
defined). 
But we can now identify the right-hand side with a series of 
meromorphic functions, which is easily seen to be convergent since 
$\varphi$ and~$\varphi'$ are bounded in any domain of~$\C$ obtained 
by removing small disks around their poles.
So we can conclude that~$\hat\Psi$ converges at the origin and 
extends to a meromorphic function.
The convergence can be made more obvious and the expression of~$\hat\Psi$
more convenient; we will give these details now.

The value of~$\hat\Psi$ at~$\xi=0$ is already known from~\eqSF:
$$
\Psi_{1} = \sum_{r=0}^{{m_0}-1} \ga_{1}(\La^r) \varphi_{{m_0},r} =
-{\txt{1\over2}} \varphi_{{m_0},0} - \sum_{r=1}^{{m_0}-1} {1\over 1-\La^r} 
\varphi_{{m_0},r}\,,
$$
so we have now two expressions for it:
$$\eqalignno{\hat\Psi(0,s)=
\Psi_{1}(s) &= -{1\over {m_0}}\sum_{k=0}^{{m_0}-1} ({1\over 2} +
              \sum_{r=1}^{{m_0}-1} {\La^{-kr}\over 1-\La^r}) 
              \varphi(s+k\Om) \EspNo\defiPsiun \cr\noalign{\medskip}
            &= -{\txt{1\over2}}\varphi(s) -\BigEisen{l\in\Z,0\le r\le 
               {m_0}-1}{(l,r)\neq0}\,\sum_{k=0}^{{m_0}-1}\,
               {\La^{-kr}\over {m_0}}(2\pi il-r\Om)\ii \varphi(s+k\Om)\,.
}$$
Substracting it from~$\hat\Psi$, we can write a uniformly 
convergent sum
$$
\displaylines{\quad 
   \hat\Psi(\xi,s) - \Psi_{1}(s) = 
\hfill\cr 
\hfill
           \bsum{l\in\Z,0\le r\le {m_0}-1}{(l,r)\neq0}\,\sum_{k=0}^{{m_0}-1}\,
           {\La^{-kr}\over {m_0}} (2\pi il-r\Om)\ii 
           \left(\varphi(s+k\Om+(2\pi il-r\Om)\ii\xi) 
           -\varphi(s+k\Om)\right) 
\quad\cr}
$$
without using Eisenstein summation. 

We have $\Om = 2 \pi i {n_0}/{m_0}$ with ${m_0}m_0'+{n_0}{n_0'}=1$ for some integers ${m_0'},{n_0'}$. 
The application
$$
\left\{ \matrix
{\Z\times\{0,\ldots,{m_0}-1\} &\longrightarrow &    \Z      \cr\noalign{\smallskip}
            (l,r)     &\longmapsto     & a = l{m_0}- r{n_0} \cr
} \right.
$$
is a bijection, the inverse of which is given by
$$
l = a {m_0'} - c {n_0}\,,\quad r = - a {n_0'} - c {m_0}
$$
where $c$ is the integer part of~$-a{n_0'}/{m_0}$.
Thus, we can use it as a change of indices:
$2\pi il-r\Om = 2\pi i a/{m_0}$, and $\La^{-kr} = e^{2\pi ika/{m_0}}$ because 
$r{n_0}\equiv -a\!\!\!\pmod {m_0}$, so we end up with the formula
$$
\hat\Psi(\xi,s) = \Psi_{1}(s) - \sum_{a\in\Z^*} \sum_{k=0}^{{m_0}-1}
    {1\over 2\pi ia}e^{2\pi ika/{m_0}} 
    \bigl( \varphi(s+k\Om + {{m_0}\xi\over 2\pi ia})-\varphi(s+k\Om)
    \bigr) .
\Eqno{\eqBT}
$$
Here is an argument for proving that this series converges uniformly and
defines a function which is meromorphic with respect to~$\xi$ for~\hbox{$\RE s<0$}:
it is sufficient to check, for any positive constant~$\rho$, 
the uniform convergence in the set 
$$
E_{\rho} = \bigl\{\, (\xi,s)\in\C\times\C\,| \ens \RE s\le -\rho \ens\hbox{and}\ens
\forall a\in\Z^{*},\, \forall b\in\Z,\, |\xi-\xi_{a,b}(s)|\ge |a|\rho \,\bigr\}
$$
(working in~$E_\rho$ means removing a small disk around each singularity in
the $\xi$-plane).
Let us fix~$\rho$ and define the set
$$
\cD_\rho = \{\, s\in\C \,|\ens \forall l\in\Z,\, 
|s - 2\pi il| \ge {m_0}\rho/2\pi \,\}
$$
so to have the following relation between $E_{\rho}$ and $\cD_\rho$:
$$
(\xi,s)\in E_{\rho} \quad \Leftrightarrow \quad \cases{
\RE s\le -\rho,                                 &  \cr\noalign{\smallskip}
s + k\Om + {{m_0}\xi\over 2\pi ia} \in \cD_\rho \ens& for $0\le k\le {m_0}-1$ and
$a\in\Z^*$; }
$$ 
note that $\RE s\le -\rho$ implies that the points $s+k\Om$
belong to~$\cD_\rho$ too.  The function~$\varphi$ is $2\pi i$-periodic and its
derivative is bounded in~$\cD_\rho$; there exists $c_\rho>0$ such that any two
points $s$ and~$s'$ in~$\cD_\rho$ can be joined inside~$\cD_\rho$ by a path of
length less than~$c_\rho|s-s'|$ followed by an integer number of 
$2\pi i$-translations, hence 
$$
\forall s,s'\in\cD_\rho , \quad
|\varphi(s')-\varphi(s)| \le M_\rho |s-s'|
\quad \hbox{with} \quad M_\rho = c_\rho \sup\{\,|\varphi'(s)|, \ens s\in\cD_\rho\,\}.
$$
This implies the uniform convergence of our series, with an explicit bound
$$\forall (\xi,s)\in E_\rho\,,\quad
|\hat\Psi(\xi,s)| \le |\Psi_1(s)| + 
    |\xi| \sum_{a\in\Z^*} {{m_0}^2 M_\rho\over 4\pi^2|a|^2}
$$
which shows the slow growth of~$\hat\Psi$ with respect to~$\xi$.
Note that the function~$\Psi_1$ is bounded in~$\cD_\rho$ (since $\varphi$ is
bounded in~$\cD_\rho$).

The function~$\varphi$ is meromorphic with only simple poles, 
located at the points~$2\pi il$ for $l\in\Z$; 
the corresponding residue is~$-1$.
Thus, for fixed~$s$, the function $\varphi(s+k\Om+{{m_0}\xi\over 2\pi ia})$ is meromorphic
with respect to~$\xi$, with only simple poles located at the points
$$
{2\pi ia\over {m_0}}(-s +{2\pi i\over {m_0}}(l{m_0} - k{n_0})) = \xi_{a,b}(s)
$$
with $b= - l{m_0} + k{n_0}$; the corresponding residue is $2\pi ia/{m_0}$ and
$k\equiv b{n_0'}\!\!\!\pmod {m_0}$, hence the value of the residue 
of~$\hat\Psi$ at~$\xi_{a,b}(s)$.

We let the reader check that
$$
\Psi_{1}(s) = \sum_{k=0}^{{m_0}-1} ({k\over {m_0}}-{1\over2}) \,\varphi(s+k\Om)\,.
\Eqno{\eqTB}
$$
from the identity~\defiPsiun.
\qed

\remark{4.4}
{Using again a decomposition formula, but this time for~$\varphi$:
$$
\varphi(s) = 1 + {1\over e^s-1} = {1\over2} + 
       \Eisen{\nu\in2i\pi\Z}{1\over s-\nu}\,,
$$
one finds the formula
$$
\hat\Psi(\xi,s) = \Psi_{1}(s) + \sum_{a\in\Z^{*},\,b\in\Z} 
C_{a,b} \left( {1\over \xi - \xi_{a,b}(s)} + {1\over \xi_{a,b}(s)} \right)
$$
with uniform convergence in any compact subset of $E_{\rho}$. 
}

\remark{4.5}
{One can write a different proof of Theorem~4.3 by starting from the
decomposition of~$f_\de$ as a sum of simple poles. We prefered to use a method
which relies only on the equation~\eqFE\ because it can be adapted in
some nonlinear problems (see Section~5.3).
}


\beginsection{4.4 A property of quasianalyticity at constant-type points}

Let us fix $\la=e^{2\pi i \al}\in \S^1$ with $\al\in[0,1[$ and a Banach space~$B$.
We now introduce spaces of functions which admit Gevrey asymptotics inside cardioids 
with cusp at~$\la$.

\Def{Definition 4.2}
{For any $\tau>0$, we define $\cGg{\tau}{\la}{B}$ to be the space of all $B$-valued
functions~$f$ such that
$u\,\mapsto\, f(\la(1-(u-1)^2))$ defines a function of~$\cG_\tau^+(1,B)$.
Analogously, we define $\cGd{\tau}{\la}{B}$ to be the space of all $B$-valued
functions~$f$ such that
$u\,\mapsto\, f(\la(1+(u-1)^2))$ defines a function of~$\cG_\tau^+(1,B)$.
}

\vskip .3cm
\epsfysize=3.3cm
\centerline{\epsfbox{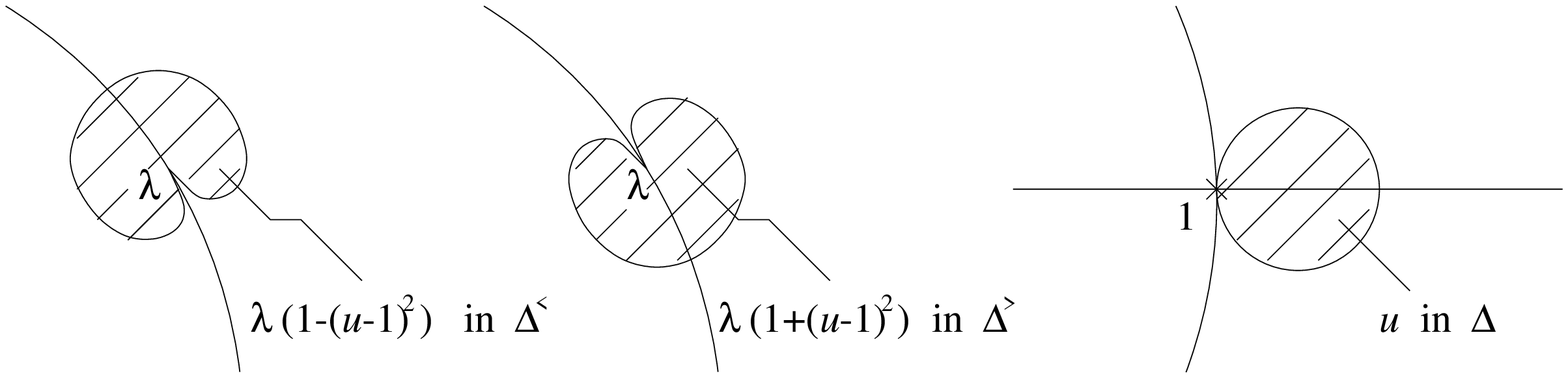}}
\vskip .3cm

Equivalently, $\cGg{\tau}{\la}{B}$ or~$\cGd{\tau}{\la}{B}$ 
is the set of the functions~$f$ which are analytic in some open set whose
boundary is a cardioid~$\Deg$ or~$\Ded$ with its cusp at~$\la$ and its axis tangent
to~$\S^1$ at~$\la$, oriented according to the picture above (such a cardioid~$\Degd$ is
nothing but the image of some disk~$\De$ by $u\,\mapsto\, \la(1\pm(u-1)^2)$),
and for which
there exist a formal series $\sum_{n\ge0} a_n Q^n \in B[[Q]]$ and positive
numbers $c_0,c_1$ such that
$$
\forall N\ge0,\; \forall q\in\Degd, \quad
\nor f(q) - \sum_{0\le n\le N-1} a_n (q-\la)^n \nor \le c_0\, c_1^N\,
\Ga(1+2\tau N) \, |q-\la|^N.
$$
Thus such a function admits Gevrey-$2\tau$ asymptotics inside the cardioid.
In particular, for $\tau=1$, we observe that $\cGg{1}{\la}{B}$
and~$\cGd{1}{\la}{B}$ are quasianalytic spaces whose members admit
Gevrey-2 asymptotics at~$\la$.

\Def{Definition 4.3}
{We define two mappings 
$\fSg,\;\fSd\,:\; \ell^1(\cR,B) \;\rightarrow\; \cO(\D\cup\E,B)$ 
by the formulas
$$
\fSgd(a)(q) = \sum_{\displaystyle\La\in\cR\cap\Sgd} \,\frac{a_\La}{q-\La}
\qquad\text{if}\quad
a=(a_\La)_{\La\in\cR}\in\ell^1(\cR,B),
$$
where
$\Sg = \ao e^{2\pi i x},\; x\in\,]\al-1/2,\al[ \af$
and
$\Sd = \ao e^{2\pi i x},\; x\in\,]\al,\al+1/2[ \af$.
}

This way, we obtain a decomposition of any \BWD\ series with
poles in~$\cR$: 
if $\la\notin\cR$, $\fS_\cR = \fSg\, + \,\fSd$
(if $\la\in\cR$, one should add the contributions of~$\la$ and~$-\la$).
This is quite reminiscent of the decomposition of the fundamental solution at
the beginning of the proof of Theorem~3.5, except that the starting point there was
Lemma~3.4 which decomposes the function according to its poles with respect to
$h=\frac{1}{2\pi i}\log\frac{q}{\la}$ rather than with respect to~$q$.

\Proc{Lemma~4.3}
{Let $r\in\,]0,1[$.
If $\la\in\oDC_\tau$ with $\tau=1$ or $\tau\ge2$,
the inclusions
$\fSg \bigl(\cS(r,B)\bigr) \subset \cGg{\tau/2}{\la}{B}$ and 
$\fSd \bigl(\cS(r,B)\bigr) \subset \cGd{\tau/2}{\la}{B}$ hold.
}

\proof
Follow the lines of the proof of Theorem~4.1, in particular adapt 
Lemma~4.1 and
choose for~$\cK$ a compact set bounded by a cardioid ($\be=3/2$).
\qed


For our purpose the previous lemma will not be of any particular
interest for~$\tau=1$, \ie for resonant points, 
whereas for $\tau=2$ it has the advantage of letting appear the quasianalytic
spaces~$\cGgd{1}{\la}{B}$ in connection with constant-type points.
But of course, for a given \BWD\ series $f=\fS_\cR(a)$, instead of dealing
with~$f$ itself that result only tells that two series~$\fSg(a)$ and~$\fSd(a)$, 
whose sum is~$f$, belong to~$\cGg{1}{\la}{B}$ or~$\cGd{1}{\la}{B}$,
and adding functions belonging to different quasianalytic classes is known to be 
a delicate matter (\cf Mandelbrojt's theorem quoted in~[Th] or~[E3], but also~[P2]).
In fact, in our situation, the relevant question is to know whether we can recover the
series~$\fSg(a)$ and~$\fSd(a)$ directly from~$f$.
A first answer is provided by the following



\Proc{Lemma~4.4}
{Assume $\la\in\oDC_\tau$ with $\tau\ge2$.
Let $r\in\,]0,1[$, $a\in\cS(r,B)$ and $q\in\D\cup\E$.
One can write
$$
\fSg(a)(q) = \frac{1}{2\pi i} \int_{\Gag(q)} \frac{\fS_\cR(a)(q_1)}{q_1-q} \,dq_1,
\quad
\fSd(a)(q) = \frac{1}{2\pi i} \int_{\Gad(q)} \frac{\fS_\cR(a)(q_1)}{q_1-q} \,dq_1,
$$
if $\Gag(q)$ (\resp $\Gad(q)$) is a simple loop with anticlockwise orientation, intersecting
$\S^1$ at~$\la$ and~$-\la$ only, transversally, and enclosing the point~$q$ and the set~$\Sd$
(\resp the set~$\Sg$).}


\vskip .3cm
\epsfysize=6.5cm
\centerline{\epsfbox{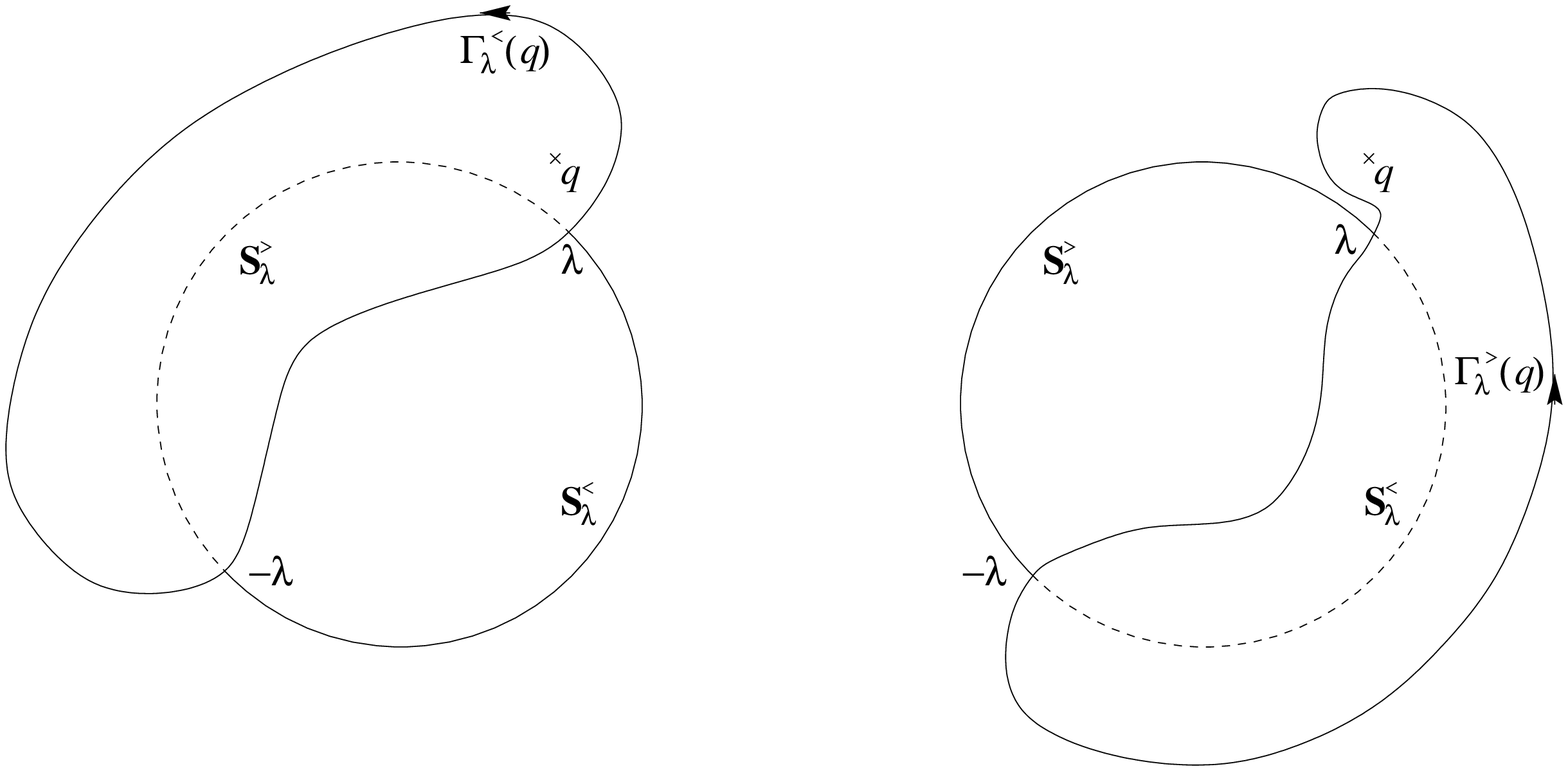}}
\vskip .3cm

The proof of Lemma~4.4 is left to the reader.

But the formulas above are ``global'' with respect to~$q$, in the sense that 
$(\fSg(a)(q))(z)$
and~$(\fSg(a)(q))(z)$ depend there on the numbers~$(\fS_\cR(a)(q_1))(z)$.
It would be more interesting to have formulas which are local in~$q$ and global in~$z$.
This turns out to be possible {\em when restricting to solutions of the cohomological equation.}

\Proc{Lemma~4.5}
{Let $\la\in\S^1\setminus\cR$ and define the coefficients
$$
\degnl = \frac{1}{n} \sum_{\La\in\cR_n\cap\Sg} \La^{-\ell}, \quad
\dednl = \frac{1}{n} \sum_{\La\in\cR_n\cap\Sd} \La^{-\ell}, \qquad
n\ge1, \; n-1\ge \ell\ge0
$$
(we recall that $\cR_n = \ao \La\in\C\,/\ens \La^n=1 \af$)
and let $r\in\,]0,1[$.
For each $q\in\D_{1/r}$, the formulas
$$
\degd(q)\,: \; z \,\mapsto\,
\sum_{\ell\ge0,n\ge \ell+1} \degdnl q^\ell \, z^n
$$
define two members $\deg(q)$ and~$\ded(q)$ of~$B_r=z H^\infty(\D_r)$.
The functions~$\deg$ and~$\ded$ are $B_r$-valued holomorphic functions in~$\D_{1/r}$
which satisfy
$$
\forall q\in\D_{1/r}, \quad
\deg(q) + \ded(q) = \de.
$$
}

\Proc{Lemma~4.6}
{Let us suppose that $\la\in\oDC_\tau$ with $\tau\ge2$,
$0<r_2<r_1$ and $r\in[r_2/r_1,1[$.
Let $g\in B_{r_1}$ and consider the corresponding solution
$f=f_g$, written as $f=\fS_\cR(a)$ where $a\in\cS(r,B_{r_2})$: for all
$q\in\D_{1/r}\setminus\S^1$,
$$
\fSg(a)(q) = \deg(q)\odot\fS_\cR(a)(q),
\quad
\fSg(a)(q) = \deg(q)\odot\fS_\cR(a)(q).
$$
}

\Pf{Proof of Lemma~4.5}
Let $n\ge1$ and $0\le \ell\le n-1$. We have obviously $|\degdnl|\le1$ and
$$
\degnl + \dednl = 
\cases{
        1 & if $\ell=0$, \cr
        \noalign{\vskip1pt}
        0 & if $1\le \ell\le n-1$. \cr
      }
$$
The Taylor series $\sum_{\ell\ge0,n\ge \ell+1} \degdnl q^\ell \, z^n$ can be written
$z \Egd(qz,z)$ with a series
$$
\Egd(x,z) = \sum_{\ell\ge0,r\ge0} \degdlrl x^\ell \, z^n
$$
which is convergent for $(x,z)\in\D\times\D$.
Thus we get functions which are holomorphic for $(q,z)\in\D_{1/r}\times\D_r$,
and for each $q\in\D_{1/r}$ we get functions $\deg(q)$ and~$\ded(q)$ which
belong to~$B_r$ and whose sum is constant and equal to~$\de$.
\qed

\Pf{Proof of Lemma~4.6}
It is sufficient to the consider the case of the fundamental solution, \ie to
prove those identities for $a=\de$.
In that case, $\fS_\cR(a) = f_\de$ and
$$
\fSgd(a) = \sum_{\dst \La\in\cR\cap\Sgd} \, \frac{\La}{q-\La} \cL_{m(\La)}(z)
         = \sum_{n\ge1} \Agdn(q) z^n,
$$
with Taylor coefficients which can written
$$
\Agdn(q) = \frac{1}{n} \sum_{\dst \La\in\cR_n\cap\Sgd} \, \frac{\La}{q-\La}
$$
(because $\dst \cL_m(z) = \sum_{n\ge m \text{\tiny s.t.} m|n} \frac{z^n}{n}$).
The identities to be proved amount to 
$$
\forall n\ge1, \qquad
\Agdn(q) = \frac{1}{q^n-1} \,\sum_{\ell=0}^{n-1} \degdnl q^\ell,
$$
which is easy to check.
\qed


\Def{Definition 4.4}
{For $\tau>0$ and $r>0$, we define $\cG^\odot_{\tau}(\la,B_r)$ to be the subspace
of~$\cG_{2\tau}(\la,B_r)$ consisting of all the functions~$f$ such that
$f\odot\deg$ extends to a function of~$\cGg{\tau}{\la}{B_r}$ and
$f\odot\ded$ extends to a function of~$\cGd{\tau}{\la}{B_r}$.
}

Putting things together we obtain

\Proc{Theorem~4.4 (Quasianalyticity at constant-type points)}
{Let $\la\in\oDC_\tau$ with $\tau\ge2$. 
For each $r\in\,]0,1[$, the fundamental solution~$f_\de$ belongs 
to the space~$\cG^\odot_{\tau/2}(\la,B_r)$, which is quasianalytic at~$\la$ if~$\tau=2$.
Thus, if $0<r_2<r_1$ and $g\in B_{r_1}$, the corresponding solution~$f_g$ belongs
to the space~$\cG^\odot_{\tau/2}(\la,B_{r_2})$, which is quasianalytic at~$\la$ if~$\tau=2$.
}



This means in particular that a solution~$f$ can be recovered from its asymptotic
expansion~$\ti f$ at a constant-type point~$\la$ by computing and ``resumming''
independently the series $\ti f\odot\deg$ and~$\ti f\odot\ded$.

\vfill \eject

\beginsection{5. Conclusions and applications}

In this final Chapter we first describe an unexpected connection
of our work with a conjecture by Gammel. Then we apply the results of 
Section 3.2 to the problem of linearization of analytic 
diffeomorphisms of the circle and we briefly sketch how the results 
of Section 4.2 can be generalized to a nonlinear small divisor problem. 


\beginsection{5.1 Gammel's series}

In a paper [Gam] published in  1974 Gammel 
studied the convergence of Pad\'e approximants to 
quasianalytic functions beyond natural boundaries
(see also [GN]). In particular he considered the 
\BWD\ series 
$$
G(q) = \sum_{m=2}^\infty\sum_{\La\in\cR^{*}_{m}}
{e^{-m}\over q-\La}. \Eqno\eqGaser
$$
As we have seen in our discussion in Section 2.2 this 
defines two complex-valued holomorphic functions, one in $\D$ and the other 
in $\E$, which have the unit circle as a natural boundary of 
analyticity. 
Gammel asked whether the function defined in $\D$ could be continued 
to the one defined in $\E$ through the natural boundary, as his 
numerical results suggested.\footnote{\noteGam}
{More precisely, Gammel asked whether the series~\eqGaser\ belongs to some
quasianalytic space of \BWD\ series, and he showed numerically that the Pad\'e
approximants~$[N/N+1]$ of~$G$ at~$q=0$ compute the value of~$G$ at $q=2$ within numerical
accuracy. Since the Pad\'e approximants depend only on the Taylor series of~$G$ 
at $q=0$, this suggested that one could continue quasianalytically~$G$ byond its 
natural boundary~$\S^1$.}

Here we want to show how our results give an affirmative 
answer to this question, but we leave untouched the quetion of the connection between
convergence of Pad\'e approximants and quasianalyticity.\footnote{\noteGN}
{Gammel's numerical results showing convergence of Pad\'e 
approximants of~$G$ beyond its circle of convergence 
could probably be justified by adapting [GN]
(which deals with the classical quasianalytic class of \BWD\ series of the form 
$\sum_{\nu=1}^\infty {A_{\nu}\over 1-q\alpha_{\nu}}$, 
with $\alpha_{\nu}$ dense on the unit circle but 
$|A_{\nu}| \le C e^{-\nu^{1+\eps}}$ for some 
$\eps >0$, which is not true for $G(q)$ which 
has $|A_{\nu}|\approx \exp (-\sqrt{\nu})$).
}

\Proc{Theorem~\thmConjGam}
{There exist $r>1$ and $g\in B_r$ such that the function 
$\chi(q,z) = q\ii(f_g(q,z)-f_g(0,1))$ satisfies
$\chi(q,1)=G(q)$ for all $q\in\D\cup\E$. As a consequence,
\item{i)} 
for all $\La_0\in\cR$, Gammel's series~$G$ belongs to the space
$(q-\La_0)\ii\cG_1(\La_0,\C)$, which is quasianalytic at~$\La_0$;
\item{ii)}
for all $\la\in\oDC_2$ and $r'\in\,]1,r[$, the function~$\chi$ belongs to the space
$\cG^\odot_1(\la,B_{r'})$, which is quasianalytic at~$\la$.
}

All the results on the Whitney smoothness and 
monogenic dependence with respect to~$q$ proved in the previous sections
apply to the function~$\chi$, thus to Gammel's series~$G$.

As for quasianalyticity, Part~{\sl i)} shows that the function~$G$ in~$\E$ can be
recovered from the knowledge of~$G_{|\D}$: one can choose any resonance~$\La_0$ and
use Borel-Laplace summation of the asymptotic expansion at~$\La_0$.
Part~{\sl ii)} yields another possibility, using the asymptotic expansion at any
constant-type point~$\la$, but for~$\chi(q,z)$ rather than for~$G$ itself:
the dependence on~$z$ is essential for that kind of quasianalyticity.

\proof
Let $A_1=0$, $A_m=e^{-m}$ for $m\ge2$.
Denoting by~$\ph$ Euler's totient function, $\ph(m)=\text{card}\cR^*_m$, we have
$$
\sum_{\La\in\cR_{m}^{*}}\frac{1}{q-\La} = 
q\ii\bigl(\ph(m)+\sum_{\La\in\cR_{m}^{*}}\frac{\La}{q-\La} 
\bigr),
$$
thus
$$
G(q) = q\ii(F(q)-F(0)), \quad
\text{with}\ens
F(q) = \sum_{m\ge1}\sum_{\La\in\cR_{m}^{*}}\frac{\La A_m}{q-\La}.
$$
In view of Proposition~A2.1, we only need to find
$g(z)=\sum_{n\ge1} g_{n}z^n$
such that 
$$
A_m = (g\odot\cL_m)_{|z=1}=\sum_{j\ge1} {g_{mj}\over mj},
\qquad
m\ge1. \Eqno\eqidexp
$$
Since $\sum m|A_m|<\infty$, we can use M\"obius inversion formula
([HW], Theorem~270, p.\ 237): we set
$$
g_n = n \sum_{j\ge1} \mu(j) A_{nj},\quad n\ge1,
$$
where the M\"obius function $\mu(j)$ is defined by~1 if $j=1$, 
$(-1)^r$ if $j$ is the product of $r$ distinct primes,
and~0 if $j$ has a squared factor.
This yields a solution of~\eqidexp,
because of the relation 
$\sum_{d|n} \mu(d) = 1$ if $n=1$ and~0 if $n\ge2$.
We observe that the radius of convergence of~$g(z)$ is~$>1$.
We have $F(q)=f_g(q,1)$,
thus we set $\chi=q\ii(f_g(q,z)-F(0))$
and we can apply Theorems~4.1 and~4.4.
\qed

In the previous example, one can check moreover that $g(z)$ has a radius of convergence
equal to~$e$ and that it defines a meromorphic function: 
$$
g(z) = -e^{-1}\,z + z\sum_{j\ge1} \mu(j)\frac{e^{-j}}{(1-z\,e^{-j})^2}.
$$
The constant~$\chi(0,1)$ involved in the description of~$G(q)$ is
$$
\sum_{m=2}^\infty e^{-m}\varphi(m) = 0.311413131378555402046127705506\ldots
$$

As is easily seen from the above proof, the statement of Theorem~\thmConjGam\ holds
for any series 
$$
G(q) = \sum_{m=2}^\infty\sum_{\La\in\cR^{*}_{m}}
{A_m\over q-\La}
$$
with $\limsup |A_m|^{1/m} < 1$.
But Gammel studies also in his paper the example corresponding to
$A_m=e^{-\sqrt{m}}$, for which quasianalyticity seems to fail as well as the
convergence of Pad\'e approximants. 
Indeed, in that case, or more generally if
$\sum m|A_m|<\infty$ but
$\limsup |A_m|^{1/m} = 1$, our arguments do not apply any longer:
there is a series $g(z)$ such that $G(q)=q\ii(f_g(q,1)-f_g(0,1))$,
but it has a radius of convergence equal to~1, which
prevents us to take~$r>1$ and thus to conclude anything for those series.

\beginsection{5.2 An application to the problem of linearization of analytic 
diffeomorphisms of the circle}

As already mentioned in the introduction, the problem of the local 
conjugacy of analytic diffeomorphisms of the circle leads to the 
linearized equation~\lineqcircle. Here we show how one can use the results of 
Section 3.2 on the existence at Diophantine 
points of Gevrey asymptotic expansions of monogenic functions in order to 
make a recent result of E. Risler [Ris] more precise.

Let $\Delta >0$, $\varepsilon >0$, $\alpha\in\C$, $\mu\in\C$. 
Following [Ris] we define:
$$
\eqalign{
\Bd{\Delta} &= \{z\in\C\,\mid\, |\IM z|<\Delta\}\; ,\cr
\Bda{\Delta}{\alpha} &= \{z\in\C\,\mid\, -\Delta  <\IM z < \Delta 
+\IM\alpha\text{if}\IM\alpha\ge 0\; ,
\cr
&\phantom{= \{z\in\C\,\mid\,  } -\Delta+\IM\alpha <\IM z < \Delta 
\text{if}\IM\alpha\le 0\}\; ,\cr
\cD(\Delta ) &= \{G\, :\Bd{\Delta} \rightarrow\C\text{analytic and 
commuting with integer translations}\}\; , \cr
\cD(\Delta ,\alpha ) &= \{G\, :\Bda{\Delta}{\alpha} 
\rightarrow\C\text{analytic and 
commuting with integer translations}\}\; , \cr
\cD_{\mu}(\Delta) &= \{G\in\cD(\Delta ) \,\mid\,
\int_{0}^{1}(G(z)-z)dz=\mu\}\; ,\cr
\cD_{\mu}^{\varepsilon}(\Delta) &= \{G\in\cD_{\mu}(\Delta ) \,\mid\,
\sup_{z\in\Bd{\Delta}}|G(z)-z-\mu |<\varepsilon\}\; ,\cr
\cD^{\varepsilon}(\Delta ) 
&=\bigcup_{\mu\in\C}\cD_{\mu}^{\varepsilon}(\Delta) \; ,\cr
\cD_{\mu}(\Delta ,\alpha) &= \{G\in\cD(\Delta ,\alpha ) \,\mid\,
\int_{0}^{1}(G(z)-z)dz=\mu\}\; ,\cr
\cD_{\mu}^{\varepsilon}(\Delta ,\alpha) &= \{G\in \cD_{\mu}(\Delta ,\alpha)
\,\mid\,
\sup_{z\in\Bda{\Delta}{\alpha}}|G(z)-z-\mu |<\varepsilon\}\; ,\cr
\cD^{\varepsilon}(\Delta ,\alpha ) &=\bigcup_{\mu\in\C}
\cD_{\mu}^{\varepsilon}(\Delta ,\alpha) \; .\cr}
$$
We will denote with $\cD_\mu^{\varepsilon ,inj}(\Delta ,\alpha )$
the set of maps in $\cD_{\mu}^{\varepsilon}(\Delta ,\alpha )$ which are 
injective on $\Bda{\Delta}{\alpha}$.

Let $\ga >0$, $\ka >0$, $d>0$ and $\beta >0$. We consider the approximation function
$$
\psi (m) = \ga \exp \Bigl( -{m\over (\log m)^{1+\beta}}\Bigr), \Eqno\afbrjuno 
$$
and the associated domain $C_{\psi ,\ka ,d}$ as in Definition~2.4.
We retain from Theorem~4, p.\ 12 of~[Ris], the following slightly weaker result:

\Proc{Theorem~\thmris\ (Local conjugacy of analytic diffeomorphisms of the 
circle with real or complex rotation numbers)}
{For all $\Delta>\delta>0$ there exist $\varepsilon >0$ and a continuous map
$$
(\alpha ,F)\in C_{\psi ,\ka ,d}\times\cD^{\varepsilon}(\Delta )
\mapsto (\ell (\alpha ,F),h_{\alpha ,F})\in\C\times
\cD_\mu^{\delta ,inj}(\Delta -\delta ,\alpha )
\Eqno\Rieslermap
$$
such that for all $(\alpha ,F)\in C_{\psi ,\ka ,d}
\times\cD^{\varepsilon}(\Delta )$ and for all $z\in \Bd{\Delta 
-\delta}$ one has
$$
\ell (\alpha , F)+F(h_{\alpha ,F}(z)) = h_{\alpha ,F}(z+\alpha )\; .
\Eqno\Rieslermapz
$$
Moreover the map~\Rieslermap\ is analytic on 
$\INT(C_{\psi ,\ka ,d})\times\cD^{\varepsilon}(\Delta )$ 
and, for all $F\in \cD^{\varepsilon}(\Delta )$, the function
$$
\ell_{F}\, :\; \al\in C_{\psi ,\ka ,d} \;\mapsto\; \ell(\alpha ,F)\in\C
\Eqno\lF
$$
is ${\cal C}^{\infty}$-holomorphic.
}

Theorem~\thmris\ is indeed a generalization of Yoccoz's theorem [Y1,Y2,Y3] on the 
linearization of analytic diffeomorphisms of the circle close to 
rotations (inasmuch as rotation numbers are allowed 
to be complex) and of Herman's [He] theorem (see also Arnol'd [Ar]) 
since the required arithmetical condition is weaker
(in [He] the real rotation numbers are assumed to be Diophantine of exponent 
$\tau\in [0,1]$). The statement in [Ris] is slightly more general 
than Theorem \thmris\ since, instead of using an 
approximation function, the real rotation numbers belong  
to any fixed relatively compact subset of the set of Brjuno 
numbers (w.r.t. a topology, finer than 
the topology induced by the usual one of $\R$, 
induced by the embedding of the Brjuno 
numbers into the space ${\ell}^{1}$ of summable sequences: 
see [Ris, pp. 6--9] for details).

The choice of the two positive constants 
$\ga$ and $\beta$ in the definition of the 
approximation function~\afbrjuno\ is arbitrary. 
Let $\psi_{j}$ denote the approximation function obtained choosing $\ga 
=\ga_{j}$, $\beta =\beta_{j}$ where $(\ga_{j})_{j\in \N}$
and $(\beta_{j})_{j\in \N}$
are two positive  decreasing sequences which tend to $0$.
>>From the previous theorem it follows that
$$
\ell_F \in {\cal M}((K_{j})_{j\in\N},\C ), \qquad K_{j} = C_{\psi_{j}, \ka ,d}.
$$
We define the Gevrey classes~$\ti\cG_\tau(y,\C)$ for~$\tau>0$ and~$y\in\R$
simply by substituting the unit circle~$\S^1$ with the real line in
Definitions~3.1 and~3.2.

\Proc{Theorem \thmrisgev}{Let $y\in\DC_{\tau}$.
The function~$\ell_F$ belongs to~$\ti\cG_{\tau'} (y ,\C )$
for all~$\tau'>\tau$.
}

\proof
This is a minor adaptation of Theorem~\thmmongev.
Following Section 3.2 very closely,
it is immediate to adapt the first part of the proof of Lemma~\lemmongev\ 
in order to see that $y\in C_{\psi_{j}, \ka ,d}$;
in fact the whole statement of Lemma~\lemmongev\ holds because again the 
points $\zeta_{n/m}$ lie between two curves with an infinite order of 
tangency to the real axis.


We then follow the proof of Theorem~\thmmongev\ and obtain inequalities which are analogous 
to~\ineqtruncated\ but involve~$\psi_j(m_\ell)$ instead of~$\text{const}e^{-\al m_\ell}$.
In order to conclude we only need to show that,
for all $j$ large enough and for all $\tau' >\tau$, there exist
two positive constants $c_{0}, c_{1}$ such that 
$$
\forall N\ge1,\quad
\sum_{m=1}^{\infty}m^{\tau (N+2)}\psi_{j}(m)\le c_{0}c_{1}^{N}
\Ga (\tau' (N+2))
$$
But this is an easy consequence of the fact that for 
all $\varepsilon >0$ one has $\lim_{m\rightarrow +\infty}
\exp (m^{1-\varepsilon})\psi_{j}(m)=0$ and one can therefore bound the 
above series using the integral 
$$
\int_{1}^{+\infty}x^{\tau (N+2)}\exp (-x^{1-\varepsilon})dx
\le {1\over 1-\varepsilon}
\Ga\left({\tau (N+2)+\varepsilon\over 1-\varepsilon}
\right).
$$
\qed


\beginsection{5.3 An application to a nonlinear small divisor problem
(semi-standard map)}

In Section~4.2 we have studied the behaviour of the solution~$f(q,z)$ 
of the
linear equation
$$
f(q,qz) - f(q,z) = g(z)
$$
for $q$ close to a resonance $\La_0=e^{2\pi i n_0/m_0}$. 
For $q$ inside or outside the unit circle~$\S^1$, the solution could be 
recovered
from its asymptotic series via Borel-Laplace summation:
$$
t f(\La_0(1+t),z) = g\odot\cL_{m_0} + \int_0^{\pm\infty} 
\hat\Phi(\xi,z)\,e^{-\xi/t}\,d\xi,
$$
and the analytic continuation of the Borel transform~$\hat\Phi$ 
w.r.t.~$\xi$ was
carefully investigated. 

We now indicate briefly that the same techniques can be adapted to a 
particular nonlinear 
equation. The reader is referred to a forthcoming paper for the proof 
of what
follows.
As for the motivation, the reader is referred to~[BMS] where the 
connection
between this nonlinear equation and
the invariant circles of the Semi-Standard Map is explained.

We restrict ourselves to $\La_0=1$ and inquire about the behaviour 
near that
``resonance'' of the solution~$F(q,z)$ of the equation
$$
F(q,qz) - 2 F(q,z) + F(q,q\ii z) = - z\,e^{F(q,z)}.
\Eqno{\eqSSM}
$$
There is an analytic solution~$F$ which, for each $q\in\D\cup\E$, is 
analytic
for $z$ close to the origin and which is characterized by $F(q,0)=0$. 
It is shown in~[BMS] that, as $q$ tends non-tangentially to~1, 
$F(q,(q-1)^2 z)$ tends to $-2\log(1+z/2)$ (in that paper the 
non-tangential
limit is computed for the other resonances as well).
We now claim that this limit is nothing but the beginning of a 
Gevrey-1
asymptotic expansion and give some indications about the 
corresponding Borel transform.

We define the moving singular half-lines to be the half-lines
$\pm\ze_b(z)[1,+\infty[$ for $b\in\Z$, with 
$$
\ze_b(z) = 2\pi(-i\log z+ i\log 2 + \pi + 2\pi b).
$$

\Proc{Theorem 5.4}
{There is an analytic function~$\hat F(\xi,z)$ which, for each 
$z\in\D_2$,
is holomorphic for $\xi$ in the complement of the half-lines
$\pm\ze_b(z)[1,+\infty[$
and has at most exponential growth on the lines passing through the 
origin and
avoiding the points~$\ze_b(z)$, such that
$$
F(1+t,t^2 z) = -2\log(1+z/2) + \int_0^{\pm\infty} \hat 
F(\xi,z)\,e^{-\xi/t}\,d\xi.
$$
In particular $F\in\cG_1(1,z H^\infty(\D_r))$ for $0<r<2$. 
}

The main difference with respect to the linear case is the necessity 
of
rescaling the variable~$z$ when $q$ approaches~1, instead of simply 
multiplying~$f(q,z)$
by some regularizing factor like~$t=q-1$, and this is precisely due 
to the nonlinear character of~\eqSSM. 
The analysis is of course more complicated, one needs to iterate a 
work which is 
analogous to that of Section~4.2, and this is why we restricted 
ourselves to the 
first resonance ($\La_0=1$) and to the holomorphic star of~$\hat F$ 
with respect to~$\xi$.
The case of the other resonances should be tractable.
We suspect that $F(1+t,t^2z)$ is resurgent with respect to~$t$, \ie that $\hat 
F(\xi,z)$
can be analytically continued with isolated singularities only, but 
this is
probably much more difficult to prove.


\vfill \eject

\beginsection{Appendix}

\beginsection{A.1 Hadamard's product}

\Def{Definition A1.1}{The {\sl Hadamard product} of two formal series
$$
A(z)=\sum_{j\ge 0} a_{j}z^j \,,\quad B(z)=\sum_{j\ge 0} b_{j}z^j 
$$
is the formal series
$$
(A\odot B)(z) = \sum_{j\ge 0} a_{j}b_{j} z^j.
$$
}

If $A$ and $B$ are convergent power series with radii of convergence 
$r_{A}$ and $r_{B}$ then $A\odot B$ converges on the disk of radius 
$r_{A}r_{B}$.

We refer to [Be] for a detailed study of Hadamard algebras, \ie 
algebras of formal power series in one variable with the product 
given by the Hadamard product.

The topological complex vector space~${\C}\{z\}$ 
with the product $\odot$
is a commutative  complex algebra with unit $\delta (z) = \sum_{j=0}^\infty 
z^j$. The Hadamard product is a convolution: 
if $A,B\in\C\{ z\}$ and~$\ga$ 
is a simple continuous curve around the origin, contained in the  
convergence domain of~$A$ and~$B$, one has 
$$
(A\odot B)(z) = {1\over 2\pi i} \int_{\ga}A(w)B\left({z\over 
w}\right){dw\over w} 
$$
for $|z|$ small enough.
The celebrated Hadamard Multiplication Theorem
states that $A\odot B$ has in all sheets of its Riemann surface 
singularities at most at points lying over $\alpha\cdot \beta$, where 
$\alpha$ is a nonregular point of $A$ and $B$ is a non regular 
point of $f$, and possibly at points lying over the origin 
[Sc]. A less general but more precise statement can be given as 
follows. 

Let $\Omega$ be an open subset of~${\C}$ and let 
${\cal O}(\Omega )$ denote the topological complex vector  
space of all functions which are 
holomorphic on~$\Omega$ with the usual locally convex topology given 
by uniform convergence on compact subsets of~$\Omega$. 

Let $\Omega_{1}$, $\Omega_{2}$ denote two open subset of ${\C}$ 
such that $0\in \Omega_{1}\cap\Omega_{2}$ and define 
$$
\Omega_{1}\odot \Omega_{2}= 
{\C}\setminus \{z\in {\C}\, \mid\, z=z_{1}z_{2} , \, 
z_{i}\notin \Omega_{i} , i=1,2\}\; . 
$$
Then Hadamard's Theorem can be stated as follows ([M\"u]):

\Proc{Theorem A1.1}{Let $\Omega_{1}$, $\Omega_{2}$ be as above, and let 
$L\in {\cal O}(\Omega_{1})$. There exists a unique continuous linear 
mapping $H_{L}$ from ${\cal O}(\Omega_{1})$ into ${\cal 
O}(\Omega_{1}\odot\Omega_{2})$ such that, for all $\ph\in {\cal 
O}(\Omega_{2})$ and for all $z$ with sufficiently small modulus, one 
has
$$
(H_{L}\ph)(z) = (L\odot \ph)(z). 
$$
}

In fact, we use mostly the case of functions analytic and bounded in disks, for which
we have the following easy result (with the notation $B_r=z H^\infty(\D_r)$ for all
$r>0$):

\Proc{Lemma A1.1}
{Let $0<\rho'<\rho$ and $L\in B_{\rho'/\rho}$. The Hadamard product defines a bounded
operator $\ph\in B_\rho \mapsto L\odot \ph\in B_{\rho'}$, whose operator norm is
$\le \Vert L \Vert_{B_{\rho'/\rho}}$.
}


\vfill\eject

\beginsection{A.2 Some elementary properties of the fundamental solution}


In this appendix we collect the statement and the proof 
of some elementary properties of the fundamental solution already used 
in the Introduction.

\Proc{Lemma A2.1}{Let $\de=z(1-z)^{-1}$.
If $q\in\C^{*}\setminus\cR$, the series 
$$
\dst
\sum_{\La\in\cR} ({q\over\La}-1)^{-1} \cL_{m(\La )} 
$$
converges to~$f_\de(q,\cdot)$ in~$\C[[z]]$.
}

We recall that if $J$ is a countable set and if  $(f_{j})_{j\in J}$
is a family of formal series, this family is summable if 
for all integer~$m$, the set $\{j\in J\,|\ens f_{j}\notin O(z^m) \}$
is finite. In this case the series  $\sum_{j\in J} f_{j}$ is called 
convergent in~$\C[[z]]$ and its sum is a formal series independent on 
the choice of an ordering on~$J$.
(This is the well-known notion of convergence associated to the 
$z$-adic valuation). 

\proof
The valuation of~$\cL_{m(\La )}$ is~$m(\La)$ and for each~$m$
the set~$\cR_{m}^{*}$ is finite. The series $f$ 
mentioned in the above lemma  converges thus formally,
and it can be rewritten as
$$
f = \sum_{m\ge1} \sum_{\La\in\cR_{m}^{*}} \sum_{j\ge1} 
{z^{jm}\over jm}({q\over\La}-1)^{-1} =
\sum_{(m,j)\in\hbox{$\N^{*}\!\times\!\N^{*}$}} 
\sum_{\La\in\cR_{m}^{*}} {z^{jm}\over jm}({q\over\La}-1)^{-1}.
$$
By reordering the terms of the summable family indexed by~$\N^{*}\times\N^{*}$, 
one finds
$$
f = \sum_{\ell\ge1} \sum_{m|\ell} \sum_{\La\in\cR_{m}^{*}} {z^\ell\over \ell} 
({q\over\La}-1)^{-1} = 
\sum_{\ell\ge1}{z^\ell\over \ell} \sum_{\La\in\cR_{\ell}} ({q\over\La}-1)^{-1}.
$$

In the coefficient of~$z^\ell$ one recognizes
the decomposition into simple elements 
of the corresponding coefficient in
$$
f_\de(q,z) = \sum_{\ell\ge1} {z^\ell \over q^\ell-1}.
$$
\qed

By means of the Hadamard 
product, the ``decomposition into simple elements'' just proved for the 
fundamental solution can be extended to the general solution~$f_{g}$
of (1.1):

\Proc{Proposition~A2.1}{Let $g\in z\C[[z]]$. 
If $q\in\C^{*}\setminus\cR$, the series
$\dst
\sum_{\La\in\cR} ({q\over\La}-1)^{-1} g\odot\cL_{m(\La )}
$
converges to~$f_g$ in~$\C[[z]]$.
}

\proof 
The identities
$$
g = g\odot\de\,,\quad 
f_g = g\odot f_\de
$$
are evident. On the other hand, for any summable family
$(f_{j})_{j\in J}$ de~$\C[[z]]$, the family $(g\odot f_{j})_{j\in J}$
is summable (because the Hadamard product with a formal series~$g$
does not decrease the  valuation), and
$$
g\odot\sum_{j\in J} f_{j} = \sum_{j\in J} g\odot f_{j}.
$$
The result follows then from Lemma~A2.1.
\qed




\Proc{Lemma A2.2}
{Let 
$$
S = \bigl\{\, q=e^{2\pi i x} \mid x\in {\R}\setminus {\Q}, \; 
\limsup_{k\rightarrow \infty}
{\log m_{k+1}\over m_{k}} = +\infty  \,\bigr\}. 
$$
where 
$(n_{k}/m_{k})_{k\ge 1}$ is the sequence of the 
convergents to~$x$ (see Appendix~A.3 for its definition and properties).
For each $q\in S$ the fundamental solution 
$
f_{\delta}(q,z) = \sum_{n\ge1} {z^n\over q^n -1}
$
diverges. $S$ is a $G_{\delta}$-dense 
subset of~$\S^{1}$ of measure zero. 
On the contrary, if $q=e^{2\pi i x}$ and 
$\limsup_{k\rightarrow\infty}{\log m_{k+1}\over m_{k}}\le M$, 
then $f_{\delta}(q,z)$ converges in the disk $|z|<e^{-M}$.
}

\proof
The divergence of $f_{\delta}$ when $q\in S$ is well-known
([HL], [Sim]), together with the convergence statement. 
$S$ is a $G_{\delta}$-dense
in~$\S^1$,
since it is 
immediate to check that 
$$
S = \bigcap_{j\ge 0} \bigcup_{n/m}
\bigl\{\, q=e^{2\pi i x} \mid |x-n/m| < {e^{-jm}\over m}\,\bigr\}. 
$$
This also shows that $S$ has measure zero.
\qed

\beginsection{A.3 Some arithmetical results. Continued fractions}

Let $[x]$ denote the usual integer 
part of a real number $x$, $\{x\}$ its fractional part: 
$\{x\} =x-[x]$. Let
$G = {\sqrt{5}+1\over 2}, \; 
g = G^{-1} = {\sqrt{5}-1\over 2}$.

To each $x \in \R \setminus \Q$ we associate its continued fraction
expansion as follows. Let
$$
x_0  = x - [x] , \quad a_0  = [x],
$$
then one obviously has
$x_0 = a_0 + x_{0}$, $a_0\in\Z$, $x_0\in\,]0,1[$. 
We will consider the iteration of the Gauss map 
$A\, :\,]0,1[\,\rightarrow [0,1[$, $A(x)=\{x^{-1}\}$:
we define inductively
$$
x_{k+1} = \biggl\{{1\over x_{k}}\biggr\},\quad
            a_{k+1} = \biggl[ {1 \over x_k} \biggr].
$$
This can be done for all $k\ge0$ since $x$ is irrational, thus 
$$
x_{k}^{-1} = a_{k+1} +  x_{k+1}, \qquad x_{k+1}\in\,]0,1[,\qquad a_{k+1}\in\N^*, 
$$
and we have
$$
x=a_0 + x_0=a_0+{1 \over a_1 + 
     x_1}= \ldots =a_0 + \displaystyle{1 \over a_1
     + \displaystyle{1  \over a_2 + \ddots +
     \displaystyle{1 \over a_k +  x_k}}}.
$$
We will write
$$
x=[a_0,a_1,\ldots ,a_k,
     \ldots]. 
$$
The integers $a_0,a_1,\ldots,a_k,\ldots$ are called the {\sl partial 
quotients} of $x$. 
The {\sl $k$th-convergent} is defined by
$$
{n_{k} \over m_{k}} = [a_0,a_1,\ldots ,a_k] =
                     a_0 + \displaystyle{1 \over a_1
     + \displaystyle{1  \over a_2 + \ddots +
     \displaystyle{1 \over a_k }}}
$$
and $\frac{n_k}{m_k}\to x$ as $k\to\infty$.
It is immediate to check that the numerators $n_{k}$ and denominators
$m_{k}$ are recursively determined by
$$
\displaylines
{n_{-2}=0, \qquad  n_{-1}=1, \qquad  n_k=a_k n_{k-1}+n_{k-2}, 
\qquad k \ge 0; \cr
 m_{-2}=1, \qquad  m_{-1}=0, \qquad  m_k=a_k m_{k-1}+ m_{k-2}, 
\qquad k \ge 0. \cr}
$$
Moreover
$$
\openup2\jot
\displaylines
{\hfill x = {n_{k}+n_{k-1}x_{k}\over m_{k}+m_{k-1}x_{k}} \hfill \LLap{\eqxnm}\cr
      x_k = - {m_{k} x -n_k \over m_{k-1} x - n_{k-1}} \cr 
 \hfill m_k n_{k-1} - n_k m_{k-1} = (-1)^k. \hfill \LLap{\eqAlt} \cr}
$$
Let
$$
\beta_k = \Pi_{i=0}^k x_i = (-1)^k (m_{k} x -n_k)\quad\hbox{for\ }k\ge 0,\quad
  \hbox{and\ }\be_{-1}=1.
$$
Then
$$
\eqalign{x_k &= {\beta_k \over \beta_{k-1}} \cr
            \beta_{k-2} &= a_k \beta_{k-1} + \beta_k \cr}
$$
>>From the definitions given one, easily proves by induction
the following proposition (we refer, for example, to~[MMY]
for its proof)

\Proc{Proposition A3.1}
{For all $x \in \R \setminus \Q$ and for all $k \ge 0$ one has
\item{(i)}\qquad $m_{k+2} > m_{k+1} > 0$;
\item{(ii)}\qquad $n_k > 0$ when $x>0$ and $n_k< 0$ when $x<0$;
\item{(iii)}\qquad $ \left|m_{k} x -n_k\right|
={\dst 1\over\dst m_{k+1}+m_kx_{k+1}}$,
so that ${\dst 1\over\dst 2}<\beta_km_{k+1}<1$;
\item{(iv)}\qquad  $\beta_k\le g^k$.
}


\remark{A3.1} 
{Note that from {\sl (iii)} and {\sl (iv)} one gets
$m_k\ge{\dst1\over\dst 2}G^{k-1}$.
}

\remark{A3.2} 
{From {\sl (iii)} one gets
$$
{1\over 2m_km_{k+1}}< {1 \over m_k (m_k + m_{k+1})}
< \left| x - {n_{k} \over m_{k}} \right| <
   {1 \over m_km_{k+1}} \Eqno\ineqfc
$$
Note also that  {\sl (iv)} remains valid for $x\in \Q$: 
in this case there exists $j\ge 0$ such that $x_{j}=0$ and the 
$x_{k}$ with $k\ge j$ are undefined; we set $\beta_{k}=0$ for all 
$k\ge j$. 
}

A partial converse of (iii), Proposition~A3.1, is provided by the 
following very useful Proposition (see [HW], Theorem 184, p.\ 153)

\Proc{Proposition A3.2}{Let $x\in\R\setminus\Q$. If 
$\left|{n\over m}-x\right|<{1\over 2m^{2}}$ then ${n\over m}$ is 
a convergent of $x$.}

The bound (iii), Proposition~A3.1, on the approximation provided by 
the convergents implies that $m_{k}|m_{k}x-n_{k}|<a_{k+1}^{-1}$. 
One can also prove the following ([HW], Theorem 193, p.\ 164)


\Proc{Proposition A3.3}{For each $x\in\R\setminus\Q$, there exist infinitely many 
rational numbers~${n\over m}$ such that
$\left|{n\over m}-x\right|<{1\over \sqrt{5}\, m^{2}}$.}


Among all rational approximations the convergents are the most 
accurate in a very precise sense: 

\Proc{Proposition A3.4 (The law of best approximation)}{If $1\le m\le 
m_{k}$, $(n,m)\not= (n_{k},m_{k})$ and $k\ge 1$, then 
$|mx-n|>|m_{k}x-n_{k}|$. Moreover, if $(n,m)\not= (n_{k-1},m_{k-1})$ and $k>1$,
then $|mx-n|>|m_{k-1}x-n_{k-1}|$.}
For a proof see [HW], Theorem 182, p.\ 151--52.

\beginsection{A.4 Proof of Lemma~3.3}

Let $\al\in\,]0,1[$ be a quadratic irrational number. 
Recall that $\N^*\times\Z$ has been partitioned into
$$
\cE^-=\ao (D,N)\in\N^*\times\Z \,|\ens N/D < \al \af
\quad\text{and}\quad
\cE^+=\ao (D,N)\in\N^*\times\Z \,|\ens N/D > \al \af.
$$
We define
$$
\nu_\pm = \ka_\pm^2 = \liminf_{(D,N)\in\cE^\pm} \ao D^2 |\frac{N}{D} - \al| \af.
$$
Thus $\nu_\pm = \nu_\pm(e^{2\pi i\al})$ with the notation of Definition~3.6.
Our aim is to find numbers~$\ka_+'$ and~$\ka_-'$, and decompositions
$$
\cE^+ = \cF^+ \cup \cE^+_* \cup \cA^+, 
\quad
\cE^- = \cF^- \cup \cE^-_* \cup \cA^-,
$$
with specific properties about the way the quantities $D^2 |\frac{N}{D} - \al|$
approach~$\ka_\pm^2$. 
\bigbreak

Let $P(X)$ be the polynomial of definition of~$\al$:
$$
P(X) = aX^2 + bX + c = a(X-\al)(X-\ov\al), \qquad a,b,c\in\Z,\qquad a\ge1,
$$
$$
   \al = \frac{-b+\eps\sqrt{\De}}{2a},\quad
\ov\al = \frac{-b-\eps\sqrt{\De}}{2a},\quad
                     \eps\in\{-1,+1\},\quad
                  \De = b^2 - 4ac \ge 2.
$$
The idea is simply to use the fact that, for all $(D,N)\in\N^*\times\Z$, the
expression
$$
a N^2 + b ND + c D^2 = a (\frac{N}{D}-\ov\al) (\frac{N}{D}-\al) D^2
$$
can assume only nonzero integral values and will allow to control 
the quantity~$|\frac{N}{D}-\al|\, D^2$ when it is small.

For $r\in\N^*$, we define the sets~$\cS^+_r$ and~$\cS^-_r$ by
$$
\cS_r^\pm = \ao (D,N)\in\cE^\pm  \,;\ens 
| a N^2 + b ND + c D^2 | = r \;\text{and}\; 
|\frac{N}{D} - \ov\al |\ge \frac{9}{10} |\al-\ov\al| \af.
$$
Let us denote by $\bigl\{\frac{n_k}{m_k}\bigr\}$ the sequence of convergents of~$\al$.
We know that 
$$
\forall p\ge0,\quad
(n_{2p},m_{2p})\in\cE^-
\ens\text{and}\ens
(n_{2p+1},m_{2p+1})\in\cE^+,
$$
moreover $\{|\frac{n_k}{m_k}-\al|m_k^2\}$ is bounded by~1.
>>From that we easily deduce that, at least for some values of~$r\in\N^*$, the first
projection of~$\cS_r^\pm$ is infinite (\ie there are inifinitely many possible
``denominators''~$D$ for which there exists $N\in\Z$ such that $(D,N)\in\cS_r^\pm$).
Therefore we can define
$$
r^\pm = \min\ao r\in\N^* \,|\ens \text{the first projection of~$\cS_r^\pm$ is
infinite}\af.
\Eqno{\defrpm}
$$

We are now ready to define the sets $\cA^\pm$, $\cE^\pm_*$, $\cF^\pm$ and the
numbers~$\ka_+'$ and~$\ka_-'$.
For the sake of simplicity we henceforth restrict ourselves to the case of the `plus'
sign.
\bigbreak

The set~$\cS^+_{r^+}$ has an infinite first projection,
whereas 
$$
\cS^+_{<r^+} = \bigcup_{1\le r < r^+} \cS^+_r
$$
has a finite first projection, and the inequality
$aN^2+bND+cD^2\ge r^++1$ holds for all~$(D,N)$ in
$$
\cS^+_{>r^+} = \bigcup_{r > r^+} \cS^+_r.
$$
The set $\cA^+$ will consist of all $(D,N)\in\cS_{r^+}^+$ with $D$ large enough;
for them we have the identity $aN^2+bND+cD^2 = r^+$.
We define the function
$\nu(\de) = \frac{r^++1}{a[|\al-\ov\al|+\de ]}$ and we pick some $\de_0>0$ such that
$\nu(\de_0)>\nu = \frac{r^+}{a|\al-\ov\al|}$. 
Notice that
$$
\nu = \frac{r^+}{\sqrt{\De}}.
$$
We also set
$$
\nu'=\nu(\de_0) > \nu,\quad \ka' = (\nu')^{1/2}.
$$
It will be checked that $\nu= \nu_+$ and $\ka_+'$ will be nothing but~$\ka'$.
The following lemma will be used in order to define progressively~$\cE^+_*$ 
and~$\cF^+$:

\Proc{Lemma~A4.1}
{For any $D_0\in\N^*$, the set $\ao(D,N)\in\cE^+\,|\ens D\le D_0\af$ 
admits a partition $\cF_{D_0}\cup\cE_{D_0}$ with
$$
\cF_{D_0}\,\text{finite and} \quad
\forall (D,N)\in\cE_{D_0},\ens
(\frac{N}{D}-\al)D^2 \ge \nu'.
$$
}

\proof
Take $\cF_{D_0} = \ao(D,N)\in\N^*\times\Z\,|\ens D\le D_0 \;\text{and}\;
[\al D]+1 \le N \le [\al D]+ \frac{\nu'}{D}
\af$ 
and $\cE_{D_0} = \ao(D,N)\in\N^*\times\Z\,|\ens D\le D_0 \;\text{and}\;
N > [\al D]+ \frac{\nu'}{D}
\af$,
where $[\,.\,]$ denotes the integer part of a real number.
\qed

We will apply this lemma and treat successively each term of the partition
$$
\cE^+ = \cT^+ \cup \cS^+_{<r^+} \cup \cS^+_{r^+} \cup \cS^+_{>r^+},
$$
where $\cT^+ = \ao (D,N)\in\cE^+  \,;\ens  
|\frac{N}{D} - \ov\al |< \frac{9}{10} |\al-\ov\al| \af$.

\bigbreak

\noindent
-- We begin with~$\cS^+_{<r^+}$. Since its first projection is finite, we decide 
to distribute its elements among~$\cF^+$ and~$\cE^+_*$ according to Lemma~A4.1.

\bigbreak

\noindent
-- Suppose $(D,N)\in\cT^+$. The inequality 
$|\frac{N}{D} - \ov\al |< \frac{9}{10} |\al-\ov\al|$
implies 
$|\frac{N}{D} - \al |\ge \frac{1}{10} |\al-\ov\al|$,
thus
$|\frac{N}{D} - \al |\,D^2\ge \nu'$ as soon as $D$ is large enough, say $D>D_0$.
Thus we put $\cT^+\cap\{D>D_0\}$ in~$\cE^+_*$ and
we distribute the elements of $\cT^+\cap\{D\le D_0\}$
among~$\cF^+$ and~$\cE^+_*$ according to Lemma~A4.1.

\bigbreak

\noindent
-- Suppose $(D,N)\in\cS^+_{>r^+}$. We know that
$$
|\frac{N}{D} - \al |\,D^2 = \frac{aN^2+bND+cD^2}{a|\frac{N}{D}-\al|}
\ge \frac{r^++1}{a|\frac{N}{D}-\ov\al|}.
$$
Either 
$0<\frac{N}{D}-\al\le\de_0$, and therefore $|\frac{N}{D}-\ov\al|\le|\al-\ov\al|+\de_0$
and
$$
(\frac{N}{D} - \al)\,D^2\ge  \frac{r^++1}{a[ |\al-\ov\al|+\de_0 ]} = \nu';
$$
or 
$\frac{N}{D}-\al>\de_0$, and 
$|\frac{N}{D} - \al |\,D^2\ge \nu'$ as soon as $D$ is large enough, say $D>D_0$.

Thus we distribute the elements of 
$\cS^+_{>r^+}\cap\{\frac{N}{D}-\al>\de_0\;\text{and}\;D\le D_0\}$
among~$\cF^+$ and~$\cE^+_*$ according to Lemma~A4.1, and the rest goes in~$\cE^+_*$.

\bigbreak

\noindent
-- Finally we suppose $(D,N)\in\cS^+_{r^+}$.
We observe that
$$
|\frac{N}{D} - \al |\,D^2= \frac{r^+}{a|\frac{N}{D} - \ov\al|}
\le \frac{r^+}{\frac{9}{10}a|\al-\ov\al|} = D_0,
$$
thus $0<N-\al D\le \frac{D_0}{D}$, and necessarily $N=[\al D]+1$ as soon as $D>D_0$.
We define 
$$
\cA^+ = \ao (D,N)\in\cS^+_{r^+} \,|\ens D>D_0 \af
$$
and apply once more Lemma~A4.1 in order to distribute the elements
of $\cA^+\cap\{D\le D_0\}$ among~$\cF^+$ and~$\cE^+_*$.
This way $\cA^+$ consists of a sequence $\{(D,N^+(D))\}_{D\in\cD^+}$, where $\cD^+$ is some
infinite subset of~$\N^*$ and $N^+(D) = [\al D]+1$.
Because of the inequalities $0<\frac{N}{D}-\al\le \frac{D_0}{D^2}\,(\forall (D,N)\in\cS^+_{r^+})$, 
we have
$$
\frac{N^+(D)}{D} \longrightarrow \al \quad \text{as}\;D\rightarrow\infty,\; D\in\cD^+
$$
$$
\text{and}\ens
|\frac{N^+(D)}{D} - \al |\,D^2= \frac{r^+}{a|\frac{N^+(D)}{D} - \ov\al|}
\longrightarrow \frac{r^+}{a|\al-\ov\al|}=\nu
\quad \text{as}\;D\rightarrow\infty,\; D\in\cD^+.
$$
\bigbreak

At this stage, we have obtained a partition of~$\cE^+$ as $\cF^+\cup\cE^+_*\cup\cA^+$
which shows that $\nu=\nu_+$. We choose $\ka_+'=\ka'$, so that $\cF^+$ and~$\cE^+_*$
satisfy the properties announced in Lemma~3.3. There only remains to study more
accurately the set~$\cA^+$.

For $D\in\cD^+$, we define
$$
\rho^+(D) = (\frac{N^+(D)}{D} - \al)\,D^2 - \nu_+ 
       = \frac{r^+}{a|\frac{N^+(D)}{D} - \ov\al|} - \frac{r^+}{a|\al - \ov\al|}.
$$
An easy computation shows that
$$
D|\rho^+(D)| = \frac{1}{|\al-\ov\al|}\Bigl( \frac{r^+}{a(\frac{N^+(D)}{D}-\al)} \Bigr)^2 \frac{1}{D}
          \sim \frac{(r^+)^2}{a^2 |\al-\ov\al|^3} \cdot \frac{1}{D}.
$$
This proves that $D\rho^+(D)$ tends to~0 as $D$ tends to infinity.

The convergence of the series $\sum_{D\in\cD^+} D^{-1/2}$ will be guaranteed by
the following

\Proc{Lemma~A4.2}
{$$
\exists \be>0 \,/\quad
\forall D,D'\in\cD^+, \ens D<D' \Rightarrow D'-D \ge \be D.
$$
}

\proof
Suppose $D,D'\in\cD^+$ with $D<D'$.
We introduce the notations
$$
\eqalign
{    D' &= D + x, \quad x\in\N^*, \cr
N^+(D')&=N^+(D)+y, \quad y\in\Z, \cr
\nu^+(D)&=(\frac{N^+(D)}{D}-\al)\,D^2, \quad z = y - \al x.
}
$$
In fact, in what follows, only $D\in\cD^+$ will be considered as a free
variable, and $x\in\N^*$ is considered as another variable subject to the
condition $D+x\in\cD^+$. The other quantitites are functions of~$D$ and~$x$, and 
we want to bound from below~$D\ii x$ for large~$D$.

An easy computation allows to rewrite the identity 
$a N^+(D')^2 + b N^+(D')D' + c D'^2 = a N^+(D)^2 + b N^+(D)D + c D^2$
as
$$\eqalign
{-(a y^2 + b xy + c x^2) 
    &= D \Bigl[ (b\frac{N^+(D)}{D}+2c) x + (b+2a\frac{N^+(D)}{D}) y \Bigr]\cr
    &= D \Bigl[ (b\al+2c) x + (b+2a\al) y + (bx+2ay)\nu^+(D)D^{-2} \Bigr]\cr
    &= D \Bigl[ z\eps\sqrt{\De} + x\nu^+(D)D^{-2}\eps\sqrt{\De} + 2az\nu^+(D)D^{-2} \Bigr]
}$$
(the last equality stems from the identities $b+2a\al = \eps\sqrt{\De}$ and $b\al+2c = -\al
\eps\sqrt{\De}$). Now the left-hand side is a nonzero integer, thus has absolute value
greater or equal to~1, and this allows us to bound from below at least one of the three
terms in the last right-hand side:
we retain that
$$
D|z| \ens\text{or}\ens x\nu^+(D)D\ii  \ens\text{or}\ens
\frac{2a\nu^+(D)}{\sqrt{\De}}D\ii|z| \ens \ge \, \frac{1}{3\sqrt{\De}}.
$$
For large~$D$ the third possibility will be excluded by the asymptotic analysis of~$z$.

The relations 
$$
N^+(D')=\al D'+\frac{\nu^+(D')}{D'}
\quad\text{and}\quad
N^+(D)=\al D+\frac{\nu^+(D)}{D}
$$
yield the formula
$$
-DD'z = D'\nu^+(D) - D \nu^+(D').
$$
We saw earlier that $\rho^+(D) \sim \text{const} D^{-2}$ as $D\rightarrow\infty$, thus
$$
\nu^+(D) = \nu_+ + \rho^+(D) = \nu_+ + \cO(D^{-2}).
$$
Similarly, since $D<D'$,
$$
\nu^+(D') = \nu_+ + \cO(D'^{-2}) = \nu_+ + \cO(D^{-2}).
$$
Here and below the symbol~$\cO$ involves a uniformness statement with respect to~$x$.
We can compute
$$
\eqalign
{-DD'z &= (D+x)(\nu^+(D)) - D(\nu^+(D')) 
       = x\nu^+(D)+ D(\nu^+(D)-\nu^+(D')) \cr
       &= x\nu^+(D) + \cO(D\ii).
}
$$
Thus 
$$
-D z = \frac{x}{D+x} \nu^+(D) + \cO(D^{-2}),
$$
in particular $D|z|$ is bounded from above and this eliminates the possibility that 
$\frac{2a\nu^+(D)}{\sqrt{\De}}D\ii|z| \ge \frac{1}{3\sqrt{\De}}$ as soon as $D$ is
large enough. 

Thus we are left with two cases:

\noindent
-- either $x\nu^+(D)D\ii \ge \frac{1}{3\sqrt{\De}}$,
thus $x\ge\be_1 D$ for $D$ large enough, with
$$
\be_1 = \frac{1}{6\nu_+\sqrt{\De}} = \frac{1}{6r^+} < 1,
$$
\noindent 
-- or $D|z| \ge \frac{1}{3\sqrt{\De}}$, and according to the above estimate of~$-Dz$, 
$$
\frac{x}{D+x} \nu_+ \ge \frac{1}{6\sqrt{\De}}
$$
for $D$ large enough, and $x\ge \frac{\be_1}{1-\be_1} D$ in that case.

Hence, in all cases, $x\ge\be_1 D$ as soon as $D>D_0$, therefore $x\ge\be D$ for 
all $D\in\cD^+$ with $\be=\min\{D_0\ii,\be_1\}$.
\qed

Thus, if we number the elements of~$\cD^+$ as an increasing
sequence~${\{D_p^+\}}_{p\ge0}$, we have $D^+_p \ge (1+\be)^p D^+_0$, and this
completes the proof of the statements relative to~$\cA^+$ which have their
counterpart for~$\cA^-$.

\bigbreak

Lastly we focus on the case of~$\cA^\eps$ with $\nu_\eps\le\nu_{-\eps}$.
In this situation we will prove that $\cA^\eps$ consists only of
couples~$(m_{w(p)},n_{w(p)})$, at least for $D$ large enough, \ie that
modulo a finite subset of it $\cA^\eps$
corresponds to convergents of~$\al$ only. This will allow us to obtain a better
control of the sequence~$\{D_p^\eps\}$.

\Proc{Lemma~A4.3}
{$$
r^\eps < \demi a |\al-\ov\al|.
$$
}

\proof
We know that there exists a subsequence ${\{n_{v(p)}/m_{v(p)}\}}_{p\ge0}$ of convergents
of~$\al$ such that 
$$
\forall p\ge0,\quad
|\frac{n_{v(p)}}{m_{v(p)}} - \al|æ\, m_{v(p)}^2 < \frac{1}{\sqrt{5}}.
$$
Suppose that $(D,N)=(m_{v(p)},n_{v(p)})$ for some~$p$.
We define $\eps_p$ to be the sign of~$(-1)^{v(p)-1}$, so that $(D,N)\in\cE^{\eps_p}$.
We have $|\frac{N}{D}-\ov\al|\ge \frac{9}{10}|\al-\ov\al|$ for $p$ large enough, and
$$
\eqalign
{|aN^2+bND+cD^2| &= a |\frac{N}{D}-\ov\al| \, |\frac{N}{D}-\al| \, D^2 \cr
                &< \frac{a}{\sqrt{5}} |\frac{N}{D}-\ov\al| \le \frac{9a}{20}|\al-\ov\al|
}
$$
for $p$ large enough,
\ie $(D,N)\in\cS^{\eps_p}_{r_p}$ with $r_p \le \frac{9a}{20}|\al-\ov\al|$.
Therefore, in view of our definition of~$r^\pm$ in the formula~\defrpm,
$r^+\le \frac{9a}{20}|\al-\ov\al|$ or $r^-\le \frac{9a}{20}|\al-\ov\al|$,
according to whether infinitely many~$(m_{v(p)},n_{v(p)})$ lie in~$\cE^+$ or in~$\cE^-$.
\qed

Now $\cA^\eps\subset\cS^\eps_{r^\eps}$. Thus, if $(D,N)\in\cA^\eps$, 
$$
|\frac{N}{D}-\al|\,D^2 = \frac{r^\eps}{a|\frac{N}{D}-\ov\al|} < \demi
$$
for $D$ large enough, hence $\frac{N}{D}$ belongs to the sequence of the convergents
of~$\al$. We can even conclude that $(D,N)=(m_{k},n_{k})$ for some $k\in\N$,
\ie that $D\wedge N=1$,
as soon as $D$ is large enough
(suppose indeed that $\forall D_0$, $\exists D>D_0$, $\exists N\in\Z$ such that
$(D,N)\in\cS^\eps_{r^\eps}$ and $D\wedge N\neq1$: 
the reduced forms $N'/D'$ of the fractions $N/D$ would yield infinitely many elements
of~$\cS^\eps_{r}$ with $1\le r<r^\eps$, and this would be in contradiction with the
definition of~$r^\eps$).

Thus there exist integers~$p_0,D_0$ and an increasing sequence~${\{k(p)\}}_{p\ge p_0}$ such
that 
$$
\cA^\eps \cap \{D>D_0\}  = \cS^\eps_{r^\eps} \cap \{D>D_0\}
                         = \ao (D^\eps_p,N^\eps(D^\eps_p))
                             = (m_{k(p)},n_{k(p)}),\; p\ge p_0
                         \af.
$$
Since $\al$ is a quadratic irrational number, by Lagrange's theorem its continued fraction
expansion is eventually periodic; we denote it by
$$
\al = [a_0,a_1,\ldots,a_{L-1},\ov{a_L,\ldots,a_{L+K-1}}],
$$
where $K\ge1$ is the period and~$L\in\N$.
The periodicity of the continued fraction expansion of~$\al$ will reflect somehow on the
structure of~$\cA^\eps$:

\Proc{Lemma~A4.4}
{Denote by $P(X)=a X^2 + b X + c$ the polynomial of definition of~$\al$ and let
$F(X,Y) = Y^2 P(\frac{X}{Y}) = a X^2 + b XY + c Y^2$.
The following identity holds for all $k\ge L$:
$$
F(n_k,m_k) = (-1)^K F(n_{K+k},m_{K+k}).
$$
}

\Proc{Corollary~A4.5}
{$$
\forall k\ge L,\quad
(m_k,n_k)\in\cA^\eps \cap \{D>D_0\} \Rightarrow 
(m_{k+2K},n_{k+2K})\in\cA^\eps \cap \{D>D_0\}.
$$
}

This corollary is sufficient to conclude the proof of Lemma~3.3.
Indeed, the sequence $\{\frac{m_{k+1}}{m_k}\}$ is bounded by some~$M>0$
(because the sequence~$\{a_k\}$ is bounded and
$\frac{m_{k+1}}{m_k}=\frac{a_k m_k + m_{k-1}}{m_k} < a_k + 1$),
and 
$$
\forall p\ge p_0,\qquad k(p)\ge L \Rightarrow
\frac{D^\eps_{p+1}}{D^\eps_p} = \frac{m_{k(p+1)}}{m_{k(p)}}
                            \le \frac{m_{k(p)+2K}}{m_{k(p)}}
                            \le M^{2K}.
$$

\medbreak \noindent {\sl Lemma~A4.4 implies Corollary~A4.5:}
Suppose $k\ge L$, $m_k>D_0$ and $(m_k,n_k)\in\cA^\eps$.
We have $m_{k+2K}>m_k>D_0$ and $(m_{k+2K},n_{k+2K})\in\cE^\eps$ because the two convergents
$\frac{n_k}{m_k}$ and~$\frac{n_{k+2K}}{m_{k+2K}}$ lie on the same side of~$\al$. 
In fact $\frac{n_{k+2K}}{m_{k+2K}}$ lies between $\frac{n_k}{m_k}$ and~$\al$, thus
$|\frac{n_{k+2K}}{m_{k+2K}}-\ov\al|>\frac{9}{10}|\al-\ov\al|$.
Therefore $\frac{n_{k+2K}}{m_{k+2K}}\in\cS^\eps_r\cap\{D>D_0\}$ with
$r=|F(m_{k+2K},n_{k+2K})|$, and Lemma~A4.4 shows that $r=|F(m_k,n_k)|=r^\eps$.
\qed

\Pf{Proof of Lemma~A4.4}
Let us first treat the case where $L=0$.
\smallbreak

We recall that $(n_{-2},m_{-2})=(0,1)$, $(n_{-1},m_{-1})=(1,0)$ and
$$
\forall k\ge0,\quad
(n_k,m_k) = a_k (n_{k-1},m_{k-1}) + (n_{k-2},m_{k-2}).
$$
The periodicity property $a_{K+k}=a_K$ allows one to check easily (by induction on~$k$) that
$$
\forall k\ge -2, \quad
(n_{K+k},m_{K+k}) = n_k(n_{K-1},m_{K-1}) + m_k(n_{K-2},m_{K-2}).
$$
On the other hand the identity
$$
\al = [a_0,a_1,\ldots,a_{K-1},\al] 
    = \frac{\al n_{K-1} + n_{K-2}}{\al m_{K-1} + m_{K-2}}
$$
shows that the polynomial
$$
P_1(X) = m_{K-1} X^2 + (m_{K-2}-n_{K-1}) X - n_{K-2}
$$
vanishes at $X=\al$, \ie
belongs to the ideal of~$\Q[X]$ generated by~$P(X)$:
$$
P_1(X) = \frac{m_{K-1}}{a} P(X).
$$
We can thus content ourselves with checking that
$$
\forall k\ge0,\quad
F_1(n_{K+k},m_{K+k}) = (-1)^K F_1(n_K,m_K),
$$
where $F_1(X,Y) = Y^2 P_1(\frac{X}{Y}) = m_{K-1} X^2 + (m_{K-2}-n_{K-1}) XY - n_{K-2} Y^2$.

This is a simple computation: for $k\ge2$,
$$
\eqalign
{F_1(n_{K+k},m_{K+k}) &= 
m_{K-1} (n_k n_{K-1} + m_k n_{K-2})^2 \cr
&\ens + (m_{K-2}-n_{K-1}) (n_k n_{K-1} + m_k n_{K-2})(n_k m_{K-1} + m_k m_{K-2}) \cr
&\ens - n_{K-2} (n_k m_{K-1} + m_k m_{K-2})^2 \cr
&= A n_k^2 + B n_k m_k + C m_k^2,
}
$$
$$
\eqalign{\text{with}\quad
A &= m_{K-1} (n_{K-1}m_{K-2} - m_{K-1}n_{K-2}) = (-1)^K m_{K-1}, \cr
B &= (m_{K-2}-n_{K-1})(n_{K-1}m_{K-2} - m_{K-1}n_{K-2}) = (-1)^K (m_{K-2}-n_{K-1}), \cr
C &= n_{K-2} (m_{K-1}n_{K-2} - n_{K-1}m_{K-2}) = (-1)^{K-1} n_{K-2}.
}
$$
This ends the proof of Lemma~A4.4 in the case $L=0$.
\bigbreak

We now proceed by induction on~$L$.
We suppose that $\al = [a_0,a_1,\ldots,a_{L-1},\ov{a_L,\ldots,a_{L+K-1}}]$ with $L\ge1$, and
that the convergents $\{n'_k/m'_k\}$ of
$$
\al' = [a_1,a_2,\ldots,a_{L-1},\ov{a_L,\ldots,a_{L+K-1}}]
$$
satisfy
$$
\forall k\ge L-1,\quad
G(n'_k,m'_k) = (-1)^K G(n'_{K+k},m'_{K+k}),
$$
where $G(X,Y) = Y^2 Q(\frac{X}{Y})$ and $Q(X)$ is the polynomial of definition of~$\al'$.

The identity
$$
\al = [a_0,\al'] = a_0 + \frac{1}{\al'}
$$
shows that the polynomial
$$
P_1(X) = (X-a_0)^2 Q(\frac{1}{X-a_0}) \in \Z[X]
$$
vanishes at $X=\al$, thus $P_1(X)$ is a rational multiple of the polynomial of definition
of~$\al$ and we can content ourselves with checking that
$$
\forall k\ge L,\quad
F_1(n_{K+k},m_{K+k}) = (-1)^K F_1(n_K,m_K),
$$
where $F_1(X,Y) = Y^2 P_1(\frac{X}{Y}) = (X-a_0 Y)^2 Q(\frac{Y}{X-a_0 Y}) = G(Y,X-a_0Y)$.

Let us express the convergents of~$\al$ in terms of those of~$\al'$:
if $k\ge1$,
$$
\frac{n_k}{m_k} = [a_0,a_1,\ldots,a_k] = a_0 + \frac{1}{[a_1,a_2,\ldots,a_k]} = a_0 +
\frac{m'_{k-1}}{n'_{k-1}},
$$
thus $n_k = a_0 n'_{k-1} + m'_{k-1}$, $m_k=n'_{k-1}$ and~$n_k-a_0 m_k=m'_{k-1}$.
Hence,
$$
\forall k\ge1, \quad
F_1(n_k,m_k) = G(n'_{k-1},m'_{k-1}),
$$
and by the inductive hypothesis
$$
\forall k\ge L,\quad
F_1(n_{K+k},m_{K+k}) = (-1)^K F_1(n_k,m_k).
$$
\qed


\beginsection{A.5 Reminder about Borel-Laplace summation}

{\em General notations and properties}
\ppar
Let $B$ a Banach algebra.
When dealing with formal series 
$\sum a_n Q^n \in B[[Q]]$,
it is convenient for us to use the variable $x=Q\ii$;
we first define the {\em formal Borel transform} (or {\em formal 
inverse Laplace
transform}) of formal series without constant term:
$$
\ti\cL\ii \,: \left\{ \eqalign
{%
x\ii B[[x\ii]]   	          &\ens\rightarrow \qquad B[[\xi]] \cr
\ti\phi = \smash{\sum_{n\ge0}}a_n x^{-n-1} &\ens\mapsto \ens \hat\phi 
= 
							\smash{\sum_{n\ge0}}
							 a_n \frac{\xi^n}{n!} .
} \right.
$$

\bigbreak
\noindent
Clearly, the Borel transform has nonzero radius of convergence if and 
only if we 
start with a formal Gevrey-1 series:
$\ti\phi\in x\ii{B[[x\ii]]}_1 \Leftrightarrow \ti\cL\ii \ti\phi \in 
B\{\xi\}$.
And starting with a convergent power-series we would obtain an entire 
function
of exponential type in all directions.

The multiplication of Gevrey-1 formal series is tranformed into convolution of
holomorphic germs:
$$
\ti\cL\ii (\ti\phi_1 \ti\phi_2) = \hat\phi_1*\hat\phi_2, \qquad
\hat\phi_i=\ti\cL\ii\ti\phi_i, \quad 
\hat\phi_1*\hat\phi_2(\xi) = \int_0^\xi 
\hat\phi_1(\xi_1)\hat\phi_2(\xi-\xi_1)\,d\xi_1.
$$
By extending the formal Borel transform to the constant series~1, we 
introduce a 
unit~$\de_0$ for the convolution:
$$
\ti\cL\ii\,: \; \ti\phi=\sum_{n\ge0}a_n x^{-n}\in {B[[x\ii]]}_1 
\;\mapsto\ens
a_0 \de_0 + \hat\phi \in B\de_0\oplus B\{\xi\},
\qquad
\hat\phi = \sum_{n\ge0} a_{n+1} \frac{\xi^n}{n!}.
$$
We will often refer to the plane of the complex variable~$\xi$ as to 
the Borel
plane, and to~$B\{\xi\}$ or~$B\de_0\oplus B\{\xi\}$ as to the 
convolutive model
in contrast with the formal model~${B[[x]]}_1$.

The counterpart of $\pa = \frac{d}{dx}$ in the convolutive model is
the multiplication by~$-\xi$:
$$
\ti\cL\ii (\pa\ti\phi) = \hat\pa (\ti\cL\ii\ti\phi),\qquad
\hat\pa\,: \left\{ \eqalign
{%
B\de_0\oplus B\{\xi\}  	 &\ens\rightarrow \ens	 B\{\xi\} \cr
a_0 \de_0 + \hat\phi\ens &\ens  \mapsto \ens\ens 	\hat\psi, 
						\qquad \hat\psi(\xi)=-\xi\hat\phi(\xi),
} \right. 
$$
while multiplication by~$x$ of a series without constant term amounts
essentially to differentiation with respect to~$\xi$: 
if $\ti\phi\in x\ii {B[[x\ii]]}_1$ and $\hat\phi = \ti\cL\ii\ti\phi$,
$$
\ti\cL\ii (x\ti\phi) = \hat\phi(0)\de_0 + \frac{d\hat\phi}{d\xi}.
$$

\bigbreak
\noindent
{\em Borel-Laplace summation}
\ppar
Let $\th\in[0,2\pi[$. Among all Gevrey-1 formal series, some of them 
have a
Borel transform $a_0\de_0+\hat\phi$ with a holomorphic 
germ~$\hat\phi$ which
extends analytically along the half-line $[0,e^{i\th}\infty[$ with at 
most
exponential growth. In such a case one can perform the{\em Laplace 
transform of
direction~$\th$}: 
$$
\hat\cL^\th \,: \ens
a_0\de_0+\hat\phi \;\mapsto\; \phi^\th,
\qquad
\phi^\th(x) = a_0 + \int_0^{e^{i\th}\infty} 
\hat\phi(\xi)\,e^{-x\xi}\,d\xi.
$$
The resulting function~$\phi^\th$ is holomorphic at least in a 
half-plane
bisected by the conjugate direction
(at least the half-plane $\RE (x\,e^{i\th}) > \de$ if we assume 
$e^{-\de|\xi|} \nor\hat\phi(\xi)\nor$ bounded).

If $\hat\phi$ extends analytically with at most exponential growth in 
a sector
$\{ \th_1 \le \arg\xi \le \th_2 \}$, by moving the direction of 
integration and
using the Cauchy Theorem we get a function analytic in a sectorial 
neighborhood
of infinity of aperture $\pi+\th_2-\th_1$.
But, according to Nevanlinna's Theorem, analyticity and exponential 
growth in a
half-strip  
$\{ \dist(\xi,[0,e^{i\th}\infty[) <\rho \}$ are sufficient to ensure 
that the initial
formal series~$\ti\phi$ is the Gevrey-1 asymptotic expansion at 
infinity in a
half-plane of~$\phi^\th$. 

The interest of this process is that $\hat\cL^\th \circ \ti\cL\ii$ 
preserves
multiplication, differentiation, etc., 
thus starting with the formal solution~$\ti\phi$ of some equation, 
studying the analytic continuation of~$\hat\phi$
and performing~$\hat\cL^\th$ for some direction~$\th$
may lead to an analytic solution of the equation (and even to distinct
solutions with the same asymptotics, if analytic continuation is 
possible in
several directions of the Borel plane with singularities in between).

\bigbreak
\noindent
{\em Effect of some changes of variable}
\ppar
Let $\ti\phi\in {B[[x\ii]]}_1$ and $\ti\cL\ii \ti\phi = a_0\de_0 + 
\hat\phi$.
Let us express the formal Borel transform of $\ti\psi(x) = 
\ti\phi(f(x))$ in
terms of that of~$\ti\phi$ for some elementary changes of 
variable~$f$.

\bigbreak
\noindent -- For $\ti\psi(x) = \ti\phi(\la x)$ with some $\la\in\C^*$,
$$
\ti\cL\ii\ti\psi = a_0 + \hat\psi, 
\qquad
\hat\psi(\xi) = \la\ii \hat\phi(\la\ii\xi).
$$

\bigbreak
\noindent -- For $\ti\psi(x) = \ti\phi(x+b)$ with some $b\in\C$,
$$
\ti\cL\ii\ti\psi = a_0 + \hat\psi, 
\qquad
\hat\psi(\xi) = e^{-b\xi} \, \hat\phi(\xi).
$$

\bigbreak
\noindent -- For $\ti\psi(x) = \ti\phi(x+\ti L(x))$ with some $\ti 
L\in x\ii{\C[[x\ii]]}_1$
and $\hat L = \ti\cL\ii\ti L$,
the Taylor formula yields
$$
\ti\cL\ii\ti\psi = a_0 + \hat\psi, 
\qquad
\hat\psi = \hat\phi + \sum_{r\ge1} \hat L^{*r} * 
\frac{\hat\pa^r\hat\phi}{r!},
\qquad
\hat L^{*r} = \underbrace{\hat L*\cdots*\hat L}_{r\rm\;times}.
$$
The above series is uniformly convergent in any closed 
disk which is contained in the disks of convergence of~$\hat\phi$ 
and~$\hat L$.
We say that $\hat\psi$ is obtained from~$\hat\phi$ by {\em
composition-convolution}, the counterpart of postcomposition 
by~$\id+L$, an
operation which may look more complicated but is in fact more 
regularizing than
postcomposition itself.

\bigbreak
\noindent
{\em Simple resurgent functions}
\ppar
In \'Ecalle's theory [E1], the holomorphic germ~$\hat\phi$ is called 
the {\em
minor} of~$\ti\phi$. The formal series~$\ti\phi$ is said to be a {\em 
simple
resurgent function} if its minor satisfies the following properties:

(i) on any broken line issuing from the origin, there is a finite set 
of points such 
that $\hat\phi$ may be continued analytically along any path that 
closely
follows the broken line in the forward direction, while circumventing 
(to the
left or to the right) those singular points;

(ii) any determination of~$\hat\phi$ in the vicinity of a singular
point~$\om$ has the form 
$$
\hat\phi(\om+\ze) = \frac{c}{2\pi i\ze} + 
\hat\psi(\ze)\frac{\log\ze}{2\pi i} +
\hat R(\ze),
\qquad
c\in B,\ens \hat\psi,\hat R\in B\{\ze\}.
$$

A nontrivial fact is the stability under convolution of this 
requirement: the
set of simple resurgent functions is a subalgebra of~${B[[x\ii]]}_1$.
We met in Section~4.2 an example of simple resurgent function where 
the minor
extended to a meromorphic function with simple poles only, thus a 
uniform
function. But since Resurgence theory is intended to deal with 
nonlinear
problems, and since convolution usually creates ramification, it is 
important
that condition~(i) authorise ramified and not only uniform analytic
continuation.\footnote{\noteRam}
{For us the source of ramification was only the 
composition-convolution induced
by some change of variable; but the fact that, when using the 
appropriate
variable, the minor was meromorphic was related to the linear 
character of the
problem under study.}

It is essential to be able to analyze the singularities which appear 
in the
convolutive model, since they are responsible for the divergence in 
the formal
model. This can be done by means of {\em alien calculus}, which 
relies on a
family of new derivations. For each $\om\in\C^*$,
there is a linear operator~$\De_\om$ of the algebra of simple 
resurgent 
functions which satisfies the Leibniz rule and measures the singular 
behaviour
of the analytic continuation at~$\om$ of the minor of the function on 
which it
is evaluated. 

For instance, if the minor~$\hat\phi$ is meromorphic, 
$\De_\om \ti\phi = 2\pi i \mathop{\hbox{Res}}(\hat\phi,\om)$. 
If the minor is not meromorphic but analytic on~$[0,\om[$ (the 
singular
point~$\om$ is ``viewed'' from the origin, and not hidden by other 
singular points),
$\De_\om \ti\phi = c + \ti\cL\hat\psi$
with notations as in~(ii).
The general formula is of the same kind but takes into account the 
singularities at~$\om$
of the various determinations of~$\hat\phi$ associated to paths which 
follow the 
segment~$[0,\om[$ while circumventing the intermediary singular 
points.

This operator~$\De_\om$ is called {\em alien derivation of 
index~$\om$} by
contrast with the natural derivation~$\pa$. There is a relation
$$
\De_\om \circ \pa = (\pa - \om) \circ \De_\om,
$$
but no relation between the alien derivations themselves: they 
generate a free
Lie algebra.
The point of view on Resurgence theory that we have indicated is 
rather
restrictive and we refer the interested reader to~[E1], [E2], [E3] 
for further
properties and more general definitions.

\vskip .5 truecm\noindent
{\bf Acknowledgements.} This research started with a visit of the first Author
to the research group ``Astronomie et Syst\`emes Dynamiques'' in 1996. It has
been supported by the CNR, the CNRS, the Institut de M\'ecanique C\'eleste, and
by CEE contract ERB-CHRX-CT94-0460.  The Authors are grateful to
J.-C.~Yoccoz for his interest and some useful discussions.  The results of this
study have been announced in various conferences between~1998 and~2000 (Bressanone,
La Rochelle, Cetraro, Aussois, Pisa, IHES and IMPA), whose organizers we wish to
thank.



\vfill \eject

\beginsection{References}

\frenchspacing


\item{[ALG]} J.-C. Archer and E. Le Gruyer, ``On the Whitney's 
extension theorem,'' {\sl Bull. Sci. Math. \bf 119} (1995), 235--266.
\item{[Ar]} V.~I. Arnold, ``On the mappings of the circumference onto 
itself,'' {\sl Translations of the Amer. Math. Soc. \bf 46} (1961), 2nd series, 
213--284.
\item{[BMS]} A. Berretti, S. Marmi and D. Sauzin,
``Limit at resonances of linearizations of some complex 
analytic dynamical systems,'' {\sl Ergodic Theory and Dynamical 
Systems \bf 20} (2000), 963--990.
\item{[Be]} B. Benzaghou, ``Alg\`ebres de Hadamard,''
{\sl Bull. de la Soc. Math. de France \bf 98} (1970), 209--252.
\item{[Be1]} A. Beurling, ``Sur les fonctions limites quasi-analytiques des
fractions rationelles,'' {\sl 8th Scandinavian Math. Congress}, 
Stockholm (1934), 199--210; in {\sl Collected Works of Arne Beurling}, 
Birkh\"auser, Boston Basel Berlin, Vol. I (1989), 109--120.
\item{[Be2]} A. Beurling, ``On quasianalyticity and general 
distributions,'' {\sl Multilithed Lecture Notes}, Summer School, 
Stanford University (1961); in {\sl Collected Works of Arne Beurling},
Birkh\"auser, Boston Basel Berlin, Vol.\  I (1989), 309--338.
\item{[Bo]} E. Borel, {\sl Le\c cons sur les fonctions monog\`enes
uniformes d'une variable complexe}, Gauthier-Villars, Paris
(1917).
\item{[Ca]} T. Carleman, {\sl Les fonctions quasi analytiques},
Gauthier-Villars, Paris (1926).
\item{[CCD]} B.~Candelpergher, M.-A.~Coppo and E.~Delabaere, ``La sommation de Ramanujan,''
{\sl L'Ensei\-gne\-ment math\'ematique \bf 43} (1997), 93--132.
\item{[De]} A. Denjoy, ``Sur les s\'eries de fractions 
rationelles,'' {\sl Bull. de la Soc. Math. de France \bf 52} (1924), 418--434.
\item{[Du]} D. Duverney, ``Explicit computation of Pad\'e-Hermite 
approximants,'' {\sl J. Approx. Th. \bf 88} (1997), 80--91.
\item{[E1]} J. \'Ecalle,
{\sl Les fonctions r\'esurgentes et leurs applications},
Publ. math. d'Orsay,
vol. I: {\bf 81--05}, vol. II: {\bf 81--06}, vol. III: {\bf 85--05}
(1981, 1985).
\item{[E2]} J. \'Ecalle, ``Singularit\'es non abordables par la 
g\'eom\'etrie,'' {\sl Ann. Inst. Fourier, Grenoble} {\bf 42} (1992), 73--194.
\item{[E3]} J. \'Ecalle, {\sl Introduction aux fonctions analysables
et preuve constructive de la conjecture de Dulac}, Hermann, Paris
(1992).
\item{[Ga]} T. W. Gamelin, {\sl Uniform Algebras}, Prentice-Hall, Engelwood Cliffs (1969).
\item{[Gam]} J. L. Gammel, ``Continuation of functions beyond natural
boundaries,'' {\sl Rocky Mountain J. Math. \bf 4} (1974), 203--206.
\item{[Gl]} G. Glaeser, ``\'Etude de quelques alg\`ebres 
tayloriennes,'' {\sl J. Anal. Math. Jerusalem \bf 6} (1958), 1--124.
\item{[GN]} J. L. Gammel and J. Nuttall, ``Convergence of Pad\'e 
Approximants to quasianalytic functions beyond natural boundaries,''
{\sl J. Math. Analys. Appl. \bf 43} (1973), 694--696.
\item{[Gou]} E. Goursat, ``Sur les fonctions \`a espaces 
lacunaires,'' {\sl Bull. Sci. Math. \bf 11} (1887), 109--114.
\item{[He]} M. R. Herman, ``Simple proofs of local 
conjugacy theorems for diffeomorphisms of the circle
with almost every rotation numbers,'' {\sl Bull. Soc. Bras. Mat. \bf 16} (1985), 45--83.
\item{[HL]} G.~H. Hardy and J.~E. Littlewood,
``Notes on the theory of series (XXIV): a curious power-series,''
{\sl Proc. Cambridge Phil. Soc. \bf 42} (1946), 85--88.
\item{[HW]} G.~H. Hardy and E.~M. Wright, {\sl An introduction to the 
theory of numbers}, Oxford University Press (1979).
\item{[Khi]} A. Ya. Khinchin, {\sl Continued Fractions}, 
The University of Chicago Press, Chicago London (1964).
\item{[Ko]} A. N. Kolmogorov, ``The General Theory of Dynamical 
Systems and Classical Mechanics'', address to the 1954 International 
Congress of Mathematicians, Amsterdam.
\item{[La]} S. Lang, {\sl Elliptic functions},
Graduate Texts in Math. {\bf 112}, Springer-Verlag (1987).
\item{[Ma]} B. Malgrange, ``Sommation des s\'eries divergentes,''
{\sl Expos. Math. \bf 13} (1995), 163--222.
\item{[MMY]} S. Marmi, P. Moussa and J.-C. Yoccoz,
``The Brjuno functions and their regularity properties,'' 
{\sl Comm. Math. Phys. \bf 186} (1997), 265--293.
\item{[M\"u]} J. M\"uller, ``The Hadamard Multiplication Theorem and 
Applications in Summability Theory,'' {\sl Complex Variables \bf 18} (1992), 155-166.
\item{[P1]} H. Poincar\'e, ``Sur les fonctions \`a espaces lacunaires,''
{\sl Amer. J. Math. \bf 14} (1892), 201--221.
\item{[P2]} H. Poincar\'e, ``Analyse de ses travaux sur la th\'eorie
g\'en\'erale des fonctions d'une variable,'' {\sl Acta Math. \bf 38} 
(1921), 65--70.
\item{[P\"o]} J. P\"oschel, ``Integrability of Hamiltonian systems on 
Cantor sets,'' {\sl Comm. Pure Appl. Math. \bf 35} (1982), 653--696. 
\item{[Ra]} J.-P. Ramis, {\sl S\'eries divergentes et th\'eories 
asymptotiques}, Panoramas et synth\'eses, Suppl. au Bull. de la Soc. 
Math. de France {\bf 121} (1993).
\item{[Re]} R. Remmert, {\sl Classical Topics in Complex Function Theory},
Graduate Texts in Math. {\bf 172}, Springer-Verlag (1998).
\item{[Ris]} E. Risler, {\sl Lin\'earisation des perturbations holomorphes 
des rotations et applications}, M\'emoires de la Soc. Math. de France {\bf 77} (1999).
\item{[Sc]} S. Schottlaender, ``Der Hadamardsche Multiplicationssatz
und weitere Kompositions s\"atze der Funktionentheorie,''
{\sl Math. Nachr. \bf 11} (1954), 239--294.
\item{[Si]} R. V. Sibiliev, ``Uniqueness theorems for Wolff-Denjoy 
series,'' {\sl St. Petersburg Math. J. \bf 7} (1996), 145--168.
\item{[Sim]} B. Simon, ``Almost periodic Schr\"odinger operators IV.
The Maryland model,'' {\sl Ann. Phys. \bf 159} (1985), 157--183.
\item{[St]} E. M. Stein, {\sl Singular integrals and differentiability 
properties of functions}, Princeton Math. Series {\bf 30} (1970).
\item{[Th]} V. Thilliez, 
``Quelques propri\'et\'es de quasi-analyticit\'e,''
{\sl Gazette des Mathematiciens \bf 70} (1996), 49--68.
\item{[Tr]} F. G. Tricomi, ``Determinazione del valore di un classico 
prodotto infinito,'' {\sl Rend. Acc. Naz. Lincei
Classe Sci. Fis. Mat. Nat. \bf XIV} (1953), 3--7.
\item{[Tj]} W. J. Trjitzinsky, ``On Quasi-Analytic Functions,''
{\sl Ann. of Math. \bf 30} (1930), 526--546.
\item{[We]} A. Weil, {\sl Elliptic Functions according to Eisenstein and 
Kronecker}, Springer-Verlag, Berlin Heidelberg New York (1976).
\item{[Wh]} H. Whitney, ``Analytic extensions of differentiable functions
defined in closed sets,'' {\sl Trans. Amer. Math. Soc. \bf 36} (1934),
63--89.
\item{[Wi]} A. Wintner, ``The linear difference equation of first 
order for angular variables,'' {\sl Duke Math. J. \bf 12} (1945), 445--449.
\item{[Wk]} J. Winkler, ``A uniqueness theorem for monogenic functions,''
{\sl Ann. Acad. Sci. Fennicae Ser. A. I. Math. \bf 18} (1993), 105--116.
\item{[Wo]} J. Wolff, ``Sur les s\'eries $\sum {A_{n}\over z-z_{n}}$,''
{\sl Comptes Rendus Acad. Sci. Paris \bf 173} (1921), 1327--1328.
\item{[Y1]} J.-C. Yoccoz, ``Conjugaison des diff\'eomorphismes analytiques du cercle,'' 
manuscript (1988).
\item{[Y2]} J.-C. Yoccoz, ``Petits diviseurs en dimension un,'' {\sl Ast\'erisque \bf 231} (1995).
\item{[Y3]} J.-C. Yoccoz, ``Analytic linearization of analytic circle 
diffeomorphisms,'' lectures given at the CIME school ``Dynamical 
Systems and Small Divisors'', Cetraro (Italy) 1998, to appear in the 
{\sl CIME series of Lecture Notes in Math.}
\item{[Za]} L. Zalcman, {\sl Analytic capacity and rational approximation},
Lecture Notes in Math. {\bf 50}, Springer (1968).


\bye